\RequirePackage{fix-cm}
\documentclass[twocolumn]{svjour3}          
\smartqed  
\usepackage{latexsym}
\usepackage[pdftex,colorlinks,breaklinks,bookmarksopen,bookmarksnumbered,citecolor=red,urlcolor=red]{hyperref}
\usepackage{graphicx,xcolor}
\usepackage{amsmath,amssymb,amsfonts}
\usepackage{fancybox}
\usepackage{subfigure}
\usepackage{float}
\usepackage{multirow}
\usepackage{bbdd} 
\usepackage{xspace}
\usepackage{pgfplots}
\pgfplotsset{compat=newest}
\usepgfplotslibrary{units}
\usetikzlibrary{pgfplots.units}
\usepackage{notation}
\def \mathbi#1{\textbf{\em #1}}
\def \und {\underline}

\def \Om {\Omega}
\def \dOm {\partial \Om}

\def \dE {\partial E}

\def \la {\lambda}

\def \Uu {\und{U}}
\def \uu {\und{u}}

\def \us {\uu_s}

\def \GM {\und{GM}}

\def \dO {\diff \Om}
\def \dS {\diff S}

\def \cont {\dsigma}
\def \defo {\depsil}

\def \K {\mathbi{K}}

\def \intO {\int_{\Om}}

\def \intOmi {\int_{\Om_i}}
\def \intOme {\int_{\Om_e}}

\def \intE {\int_{E}}
\def \intdE {\int_{\partial E}}

\def \cre {\mathrm{cre}}

\def \std {\mathrm{std}}

\def \loc {\mathrm{loc}}

\def \part {\mathrm{part}}

\def \Ec{\mathcal{E}}

\def \Hc {\mathcal{H}}
\def \Ic {\mathcal{I}}
\def \Jc {\mathcal{J}}

\def \Lc {\mathcal{L}}
\def \Mc {\mathcal{M}}
\def \Nc {\mathcal{N}}

\def \Pc{\mathcal{P}}
\def \Qc{\mathcal{Q}}

\def \Sc {\mathcal{S}}

\def \Uc {\mathcal{U}}
\def \Vc {\mathcal{V}}

\def \Ucb {\boldsymbol{\Uc}}
\def \Vcb {\boldsymbol{\Vc}}
\def \Scb {\boldsymbol{\Sc}}

\def \Abb {\mathbb{A}}

\def \Cbb {\mathbb{C}}

\def \Lbb {\mathbb{L}}

\def \Pbb {\mathbb{P}}

 \journalname{Computational Mechanics}
\begin{document}

\title{An enhanced method with local energy minimization for the robust a posteriori construction of equilibrated stress fields in finite element analyses
}

\titlerunning{An enhanced method with local energy minimization for the construction of equilibrated stress fields}        

\author{Florent Pled         \and
        Ludovic Chamoin         \and
        Pierre Ladev\`eze 
}


\institute{F. Pled \at
              LMT-Cachan (ENS-Cachan/CNRS/Paris 6 University), \\
              61 Avenue du Pr\'esident Wilson, 94235 CACHAN Cedex, France \\
              \email{pled@lmt.ens-cachan.fr}           
           \and
           L. Chamoin \at
           LMT-Cachan (ENS-Cachan/CNRS/Paris 6 University) \\
              \email{chamoin@lmt.ens-cachan.fr}
           \and
           P. Ladev\`eze \at
           LMT-Cachan (ENS-Cachan/CNRS/Paris 6 University), \\
           EADS Foundation Chair, Advanced Computational Structural Mechanics, France \\
              \email{ladeveze@lmt.ens-cachan.fr}
}

\date{Received: 29 July 2011 / Accepted: 21 August 2011}

\maketitle

\begin{abstract}

In the context of global/goal-oriented error estimation applied to computational mechanics, the need to obtain reliable and guaranteed bounds on the discretization error has motivated the use of residual error estimators. These estimators require the construction of admissible stress fields verifying the equilibrium exactly. This article focuses on a recent method, based on a flux-equilibration procedure and called the \textit{element equilibration + star-patch technique} (EESPT), that provides for such stress fields. The standard version relies on a strong prolongation condition in order to calculate equilibrated tractions along finite element boundaries. Here, we propose an enhanced version, which is based on a weak prolongation condition resulting in a local minimization of the complementary energy and leads to optimal tractions in selected regions. Geometric and error estimate criteria are introduced to select the relevant zones for optimizing the tractions. We demonstrate how this optimization procedure is important and relevant to produce sharper estimators at affordable computational cost, especially when the error estimate criterion is used. Two- and three-dimensional numerical experiments demonstrate the efficiency of the improved technique.
\keywords{Verification \and Finite element method \and Admissible stress field \and Non-intrusive techniques \and Strict error bounds}
\end{abstract}

\section{Introduction}\label{1}

In a wide variety of engineering disciplines, verification of the quality of the numerical modeling of physical systems has become an important issue at both industrial and research levels. Starting from an initial mathematical model, referred to as the \textit{reference model}, and coming from continuum mechanics, one usually constructs a discretized model suited to current numerical engineering tools. One of the most powerful and popular design tools is the finite element method (FEM); it is extensively used to obtain approximate numerical solutions. The mastering and control of the quality of a finite element analysis boomed about $40$ years ago \cite{Lad75,Lad83,Bab78a,Zie87}. Pioneering developments of effective methods concerning the assessment of the global error discretization provided a reliable mean to control the global quality of a FE simulation \cite{Ver96,Bab01,Ste03,Lad04}. Nowadays, research activities are turning to goal-oriented error estimation, \ie assessment of the error on local quantities providing local error bounds \cite{Par97,Ran97,Cir98,Per98,Lad99a,Pru99,Str00,Bec01,Wib06,Lad08,Cha07,Cha08,Pan10,Lad10}. One of the topical key issues concerns robust global/goal-oriented error estimation methods, \ie techniques providing strict and relevant bounds on the error. Such methods currently require the construction of an admissible stress field, \ie a stress tensor that verifies the equilibrium equations exactly.

Several techniques currently enable to construct an admissible stress field. The first approach, based on a equilibrium-type finite element method, is not realistic today despite its remarkable efficiency, as it requires another global solution of the reference problem through a dual formulation. The second approach, called the \textit{element equilibration technique} (EET) \cite{Lad83,Ain93bis,Lad96,Lad97,Flo02,Par10}, is built on an energy relation linking the FE data and the searched stress field, called the strong prolongation condition, and allowing the construction of equilibrated tractions along element boundaries. Starting from those balanced tractions and FE data, an approximate resolution of local problems defined at the element scale leads to the calculation of an admissible stress field at reasonable computational cost. The third approach, called the \textit{star-patch equilibration technique} (SPET) \cite{Mac00,Car00,Mor03,Pru04,Par06,Par09,Moi09,Cot09}, is built on the partition of unity concept allowing the resolution of local self-equilibrated problems defined at the patch of elements scale and leading to accurate admissible stress fields at higher computational cost. The last existing approach, called the \textit{element equilibration + star-patch technique} (EESPT) and recently introduced in \cite{Lad10bis,Ple11}, results of combining the former two approaches as it is built on both the strong prolongation condition and the partition of unity method leading to equilibrated tractions whose construction is easier to implement. This last method has attractive features, as it seems to be a good trade-off between performance, computation cost and simplicity of implementation \cite{Lad10bis,Ple11}. Nevertheless, the element-by-element calculation of an admissible stress field is similar to that of the EET. Besides, it is worth noticing that EET and EESPT methods are similar in the case of first-order FE interpolation degree. Only their practical implementations differ.

In this work, we first revisit the main features of this new hybrid technique, namely the EESPT. This paper is a continuation of previous papers \cite{Lad10bis,Ple11}, where the EESPT method is introduced. We then go one step further and focus on an enhanced version of the EESPT method, inspired from an idea introduced in \cite{Lad97} and subsequently developed in \cite{Flo02}. This improved version is based on a weak prolongation condition resulting in local minimization of the complementary energy that leads to optimal tractions. The main thrust resides in the introduction of geometric and error estimate criteria which allow to select the relevant zones for optimizing the tractions, namely the highly distorted and mostly concentrated error regions. This enables to reduce the error estimate without increasing significantly the corresponding computational cost.

The paper is divided into seven sections: after this introduction, \Sect{2} introduces the reference problem and the finite element one considered in this work in order to introduce the basic notion of admissible fields; \Sect{3} revisits the standard version of the EESPT me-thod based on a strong prolongation condition, while \Sect{4} presents the enhanced version for constructing fluxes based on a weak prolongation condition; \Sect{5} deals with the main technical features regarding its practical implementation; the capabilities of the proposed optimized approach are illustrated through several two- and three-dimensional numerical examples in \Sect{6}; eventually, \Sect{7} suggests recommendations and topics for future research.

\section{Basics on error estimation and admissible solutions}\label{2}

\subsection{Statement of the reference problem}\label{2.1}

\begin{figure}
\centering\includegraphics[scale=0.22]{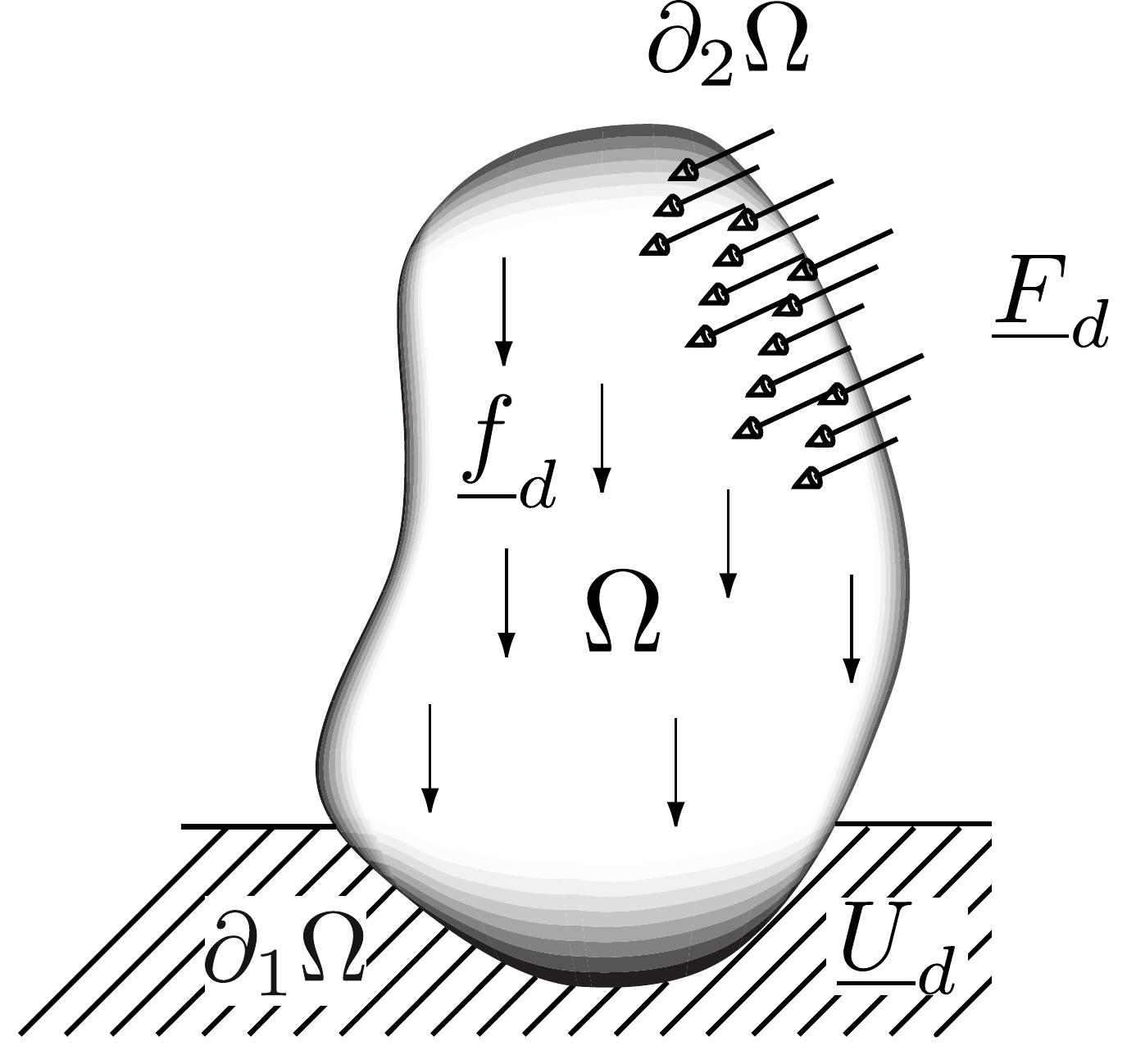}
\caption{Representation of the structure and its environment.}
\label{fig1:pbref_black_wo_line}
\end{figure}

Let us consider a mechanical structure defined in an open bounded domain $\Om$, with boundary $\dOm$ (see \Fig{fig1:pbref_black_wo_line}), and subjected to a prescribed mechanical loading: a displacement field $\Uu_d$ on part $\partial_1 \Om \neq \varnothing$; a traction force density $\und{F}_d$ on the complementary part $\partial_2 \Om$ of $\dOm$ such that $\overline{\partial_1 \Om \cup \partial_2 \Om} = \dOm$, $\partial_1 \Om \cap \partial_2 \Om = \varnothing$; a body force field $\und{f}_d$ within $\Om$. 

Besides, we consider a material with isotropic, homogeneous, linear and elastic behavior under the assumptions of small perturbations state, quasi-static loading and isothermal case. 

The reference problem to be solved reads as follows:
Find a displacement/stress pair $(\uu ,\cont)$ in the space domain $\Om$, which verifies:
\begin{itemize}
\item[$\bullet$]
the kinematic conditions:
\begin{equation}\label{eq1:CAref}
\uu \in \Ucb; \quad \uu_{\restrictto{\partial_1 \Om}} = \Uu_d; \quad \defo(\uu) = \frac{1}{2}\big(\nabla \uu + \nabla^T \uu \big);
\end{equation}
\item[$\bullet$]
the equilibrium equations:
\begin{equation}\label{eq1:SAref}
\begin{aligned}
\cont \in \Scb; \quad \forall \: \uu^{\ast} \in \Ucb_0, \intO \Tr\big[\cont \: \defo(\uu^{\ast})\big] \dO \\
= \intO \und{f}_d \cdot \uu^{\ast} \dO + \int_{\partial_2 \Om}\und{F}_d \cdot \uu^{\ast} \dS;
\end{aligned}
\end{equation}
\item[$\bullet$]
the constitutive relation:
\begin{equation}\label{eq1:RDCref}
\cont(M) = \K \: \defo\big(\uu(M)\big) \quad \forall \: M \in \Om,
\end{equation}
\end{itemize}

where $\defo(\uu)$ represents the classical linearized strain tensor associated to displacement field $\uu$, while operator $\K$ stands for Hooke's elasticity tensor. $\Ucb = \left\{\uu \in [\Hc^1(\Om)]^3 \right\}$ and $\Scb = \left\{ \cont \in \Mc_s(3) \cap [\Lc^2(\Om)]^6 \right\}$ are the functional regularity spaces ensuring the existence of finite-energy solutions, where $\Mc_s(n)$ denotes the space of symmetric square matrices of order $n$, $\Hc^1(\Om)$ represents the standard Sobolev space of square integrable functions and first derivatives, and $\Lc^2(\Om)$ refers to the space of square integrable functions. $\Ucb_0 \subset \Ucb$ represents the vectorial space associated to $\Ucb$, \ie the space of functions satisfying homogeneous kinematic (Dirichlet boundary) conditions over $\partial_1 \Om$: $\Ucb_0 = \left\{ \uu^{\ast} \in \Ucb, \uu^{\ast}_{\restrictto{\partial_1 \Om}} = \und{0} \right\}$.

In the following, the exact solution of the reference problem is denoted by $(\uu, \cont)$. As this one remains inaccessible in practice, one has recourse to approximate resolution methods in order to achieve an approximate solution of that reference problem.

\subsection{Statement of the finite element problem}\label{2.2}

The standard Galerkin Finite Element Method (FEM), which is a well-established computer-aided engineering tool, is one of the most commonly used methods. It furnishes a numerical solution $(\uu_h, \cont_h)$ lying in the finite-dimensional spaces $\Ucb_{h} \times \Scb_{h} \subset \Ucb \times \Scb$. These are defined from piecewise continuous polynomial displacement shape functions associated with a spatial discretization (finite element space mesh $\Mc_h$) of the domain $\Om$. It is assumed that the prescribed displacement field $\Uu_d$ is compatible with the interpolation chosen for the FE discretization. Thus, the finite element problem to be solved reads as follows: 
Find a displacement/stress pair $(\uu_h (M),\cont_h (M)), M \in \Om$, which verifies:
\begin{itemize}
\item[$\bullet$]
the kinematic conditions:
\begin{equation}\label{eq1:CAEF}
\uu_h \in \Ucb_{h}; \quad {\uu_h}_{\restrictto{\partial_1 \Om}} = \Uu_d; \quad \defo(\uu_h) = \frac{1}{2}\big(\nabla \uu_h + \nabla^T \uu_h \big);
\end{equation}
\item[$\bullet$]
the equilibrium equations:
\begin{equation}\label{eq1:SAEF}
\begin{aligned}
\cont_h \in \Scb_h; \quad \forall \: \uu_h^{\ast} \in \Ucb_{h,0}, \intO \Tr\big[\cont_h \: \defo(\uu_h^{\ast})\big] \dO \\
= \intO \und{f}_d \cdot \uu_h^{\ast} \dO + \int_{\partial_2 \Om}\und{F}_d \cdot \uu_h^{\ast} \dS;
\end{aligned}
\end{equation}
\item[$\bullet$]
the constitutive relation:
\begin{equation}\label{eq1:RDCEF}
\cont_h(M) = \K \: \defo\big(\uu_h(M)\big) \quad \forall \: M \in \Om,
\end{equation}
\end{itemize}

where $\Ucb_{h,0} = \Ucb_h \cap \Ucb_0$.

In the displacement-type finite element framework, the FE solution $(\uu_h, \cont_h)$ satisfies both kinematic conditions (\ref{eq1:CAref}) and constitutive relation (\ref{eq1:RDCref}) of the reference problem, but fails to verify equilibrium equations (\ref{eq1:SAref}). These equilibrium deficiencies are the main approximation in the displacement-type FEM.

First, let us define the discretization error $\und{e}_h = \uu - \uu_h$, also called the exact error or true error, corresponding to the difference between the exact displacement solution and the FE one; the assessment of this error enables to control the numerical quality of the FE solution $(\uu_h, \cont_h)$. Usually, it is measured in terms of a suitable norm, such as the energy norm $\lnorm{\bullet}_{u, \Om} = \left( \intO \Tr\big[\K \: \defo(\bullet) \: \defo(\bullet)\big] \dO \right)^{1/2}$, which leads to a global discretization error $\lnorm{\und{e}_h}_{u, \Om}$. Secondly, local errors $e^{\loc}_h = I(\uu) - I(\uu_h)$ can be defined if one seeks to evaluate and measure errors on quantities of interest, \ie functional outputs $I(\uu)$ of the solution.

\subsection{Construction of admissible fields in the standard FEM framework}\label{2.3}

The need of obtaining reliable and guaranteed bounds of the discretization error has motivated the development of methods for constructing an admissible solution; those are currently the only way to achieve strict bounds on the error \cite{Lad04,Lad83}. The admissible pair, denoted $(\hat{\uu}_h, \hat{\cont}_h)$, should verify the kinematic conditions (\ref{eq1:CAref}) and equilibrium equations (\ref{eq1:SAref}) of the reference problem. As the most widespread finite element methods use a classical displacement formulation providing a kinematically admissible displacement field $\uu_h$, one usually chooses $\hat{\uu}_h = \uu_h$ for the sake of simplicity, apart from the case of incompressible materials (incompressibility being considered as an additional kinematic admissibility constraint, see \cite{Lad92}). Therefore, one focuses on the construction of an admissible stress field $\hat{\cont}_h$, which is the key technical ingredient. An overall description of the different techniques used to reconstruct such stress field has been presented in \Sect{1}, and one of these techniques will be detailed in \Sect{3}.

Starting from an admissible solution $(\hat{\uu}_h, \hat{\cont}_h)$, the measure $e_{\cre}(\hat{\uu}_h, \hat{\cont}_h) = \lnorm{\hat{\cont}_h - \K \: \defo(\hat{\uu}_h)}_{\cont,\Om}$ of the constitutive relation error (\ref{eq1:RDCref}) enables one to assess the measure of the global discretization error $\lnorm{\und{e}_h}_{u, \Om} = \lnorm{\uu - \hat{\uu}_h}_{u, \Om}$ in the sense of the energy norms $\lnorm{\bullet}_{u, \Om}$ and $\lnorm{\bullet}_{\cont, \Om} = \left( \intO \Tr\big[ \bullet \: \K^{-1} \: \bullet\big] \dO \right)^{1/2}$ without knowing the exact solution $\uu$. Indeed, the error in constitutive relation is connected to the classical discretization error in solution by the popular Prager-Synge hypercircle theorem \cite{Pra47} which reads:
\begin{equation}\label{eq1:PragerSynge}
\lnorm{\uu - \hat{\uu}_h}^2_{u, \Om} + \lnorm{\cont - \hat{\cont}_h}^2_{\cont,\Om} = \lnorm{\hat{\cont}_h - \K \: \defo(\hat{\uu}_h)}^2_{\cont,\Om},
\end{equation}
and leads to the inequality:
\begin{equation}
\lnorm{\und{e}_h}_{u, \Om} \leqslant e_{\cre}(\hat{\uu}_h, \hat{\cont}_h).
\end{equation}
Thus, the constitutive relation error $e_{\cre}(\hat{\uu}_h, \hat{\cont}_h)$ is a reliable error measure as it is a strict upper bound of the measure $\lnorm{\und{e}_h}_{u, \Om}$ of the exact discretization error. The quality of the obtained error bound is strongly dependent on the quality of the corresponding admissible stress field. Assessment of the accuracy of an estimator is commonly expressed in terms of the usual effectivity index with respect to the energy norm of a reference error (obtained using an ``overkill solution'' resulting from the use of a very refined mesh):
\begin{equation}
\eta = \frac{\theta}{\lnorm{\und{e}_h}_{u, \Om}}, \notag
\end{equation}
where $\theta$ denotes the error estimate and $\lnorm{\und{e}_h}_{u, \Om}$ stands for the energy norm of the exact discretization error (if available) or that of the reference error. The more the effectivity index $\eta$ is close to $1$, the more the estimator is relevant and the more the corresponding admissible stress field is similar to the unknown exact stress field, see (\ref{eq1:PragerSynge}). Let us now examine the main points regarding the standard and enhanced versions of the EESPT technique.

\section{Principles of the original version of the EESPT technique}\label{3}

\subsection{Notations}\label{3.1}

Let us define $\Ec$, $\Nc$, $\Ic$, $\Nc \setminus \Ic$ and $\Jc$ the set of elements, nodes, vertices, non-vertex nodes and edges of the FE mesh $\Mc_h$, respectively. $\Ec_i^{\Nc} \subset \Ec$, $\Ec_i^{\Ic} \subset \Ec$, $\Ec_i^{\Nc \setminus \Ic} \subset \Ec$ and $\Ec_{\Gamma}^{\Jc} \subset \Ec$ represent the set of elements connected to node $i$, vertex $i$, non-vertex node $i$ and edge $\Gamma$, respectively. $\Jc_i^{\Ic} \subset \Jc$ represents the set of edges connected to vertex $i$. $\Nc_E^{\Ec} \subset \Nc$ and $\Nc_{\Gamma}^{\Jc} \subset \Nc$ stand for the set of nodes associated with element $E$ and edge $\Gamma$, respectively. $\Ic_E^{\Ec} \subset \Ic$ and $\Ic_{\Gamma}^{\Jc} \subset \Ic$ designate the set of vertices connected to element $E$ and edge $\Gamma$. $\Nc_E^{\Ec} \setminus \Ic_E^{\Ec} \subset \Nc \setminus \Ic$ and $\Nc_{\Gamma}^{\Jc} \setminus \Ic_{\Gamma}^{\Jc} \subset \Nc \setminus \Ic$ denote the set of non-vertex nodes connected to element $E$ and edge $\Gamma$. Finally, the FE displacement field $\uu_h$ is assumed to belong to $\Ucb^{p}_{h}$, where $\Ucb^{p}_{h}$ corresponds to the FE interpolation space of degree less than or equal to $p$, $p$ being the FE interpolation degree. $\Uc^{p}_{h}$ refers to its one-dimensional correspondent.

\subsection{The element equilibration and star-patch technique (EESPT): principle of the construction}\label{3.2}

This technique, developed in Ladev\`eze \textit{et al} \cite{Lad10bis}, is a hybrid method as it combines the advantages of both EET and SPET methods introduced in \Sect{1}. The procedure to construct an admissible stress field is carried out in two main steps:
\begin{enumerate}
\item[(i)]
construction of tractions $\hat{\und{F}}_h$ in equilibrium with the external loading $(\und{F}_d, \und{f}_d)$ on element edges $\dE$ of the spatial mesh $\Mc_h$;\\
\item[(ii)]
calculation of an admissible stress field $\hat{\cont}_h$ solution of a static local problem over each element $E \in \Ec$ where equilibrated tractions $\hat{\und{F}}_h$ act as Neumann boundary conditions.
\end{enumerate}

The key ingredient used to set up an admissible stress field is an energy condition, called the \textit{strong prolongation condition}. This one consists of seeking $\hat{\cont}_h$ as an extension (or prolongation) of the FE stress field $\cont_h$ in the following sense:
\begin{equation}\label{eq1:prolongfort}
\intE \left(\hat{\cont}_h - \cont_h \right) \: \und{\nabla} \varphi_i \dO = \und{0} \quad \forall \: E \in \Ec, \ \forall \: i \in \Nc_E^{\Ec},
\end{equation}
where $\varphi_i \in \Uc^{p}_{h}$ represents the FE shape function associated with node $i$.

Tractions $\hat{\und{F}}_h$ involved in the first step are designed to represent the stress vectors $\hat{\cont}_{h \restrictto{E}} \: \und{n}_E$ on sides $\dE$ of element $E \in \Ec$:
\begin{equation}\label{eq1:defdensites}
\hat{\cont}_{h \restrictto{E}} \: \und{n}_E = \eta_E \: \hat{\und{F}}_h \quad \textrm{on} \ \dE,
\end{equation}
where $\und{n}_E$ is the unit outward normal vector to element $E$ and $\eta_E = \pm 1$ are constant functions ensuring continuity of the stress vector in the $\und{n}_E$ direction across element boundaries $\dE$.

Besides, these tractions are constructed in equilibrium with the mechanical external loads $(\und{F}_d, \und{f}_d)$, that reads:
 \begin{align}
 & \eta_E \: \hat{\und{F}}_h = \und{F}_d \quad \textrm{on} \ \dE \subset \partial_2 \Om \label{eq1:eqFd} \\ 
 & \intdE \eta_E \: \hat{\und{F}}_h \cdot \us \dS + \intE \und{f}_d \cdot \us \dO = 0 \quad \forall \: \us \in \Ucb_{S \restrictto{E}}, \label{eq1:eqfd}
 \end{align}
where $\Ucb_{S \restrictto{E}}$ denotes the set of rigid body displacement fields over element $E$.

Once the set of equilibrated tractions $\hat{\und{F}}_h$ has been constructed on the element sides, the second step merely consists of searching the local restriction $\hat{\cont}_{h \restrictto{E}}$ of an admissible stress field $\hat{\cont}_h$ to each element $E \in \Ec$ as the solution of the following local problem $\Pc_E^{\Ec}$:
\begin{equation}\label{eq1:pblocal}
\hat{\cont}_{h \restrictto{E}} \in \Scb_{\hat{\und{F}}_h} \iff \left\{
\begin{aligned}
& \hat{\cont}_{h \restrictto{E}} \in \Scb \\
& \und{\diver} \: \hat{\cont}_{h \restrictto{E}} + \und{f}_d = \und{0} \quad \textrm{in} \ E \\
& \hat{\cont}_{h \restrictto{E}} \: \und{n}_E = \eta_E \: \hat{\und{F}}_h \quad \textrm{on} \ \dE
\end{aligned}
\right.
\end{equation}

One can obtain an approximation of $\hat{\cont}_{h \restrictto{E}}$ by using a standard displacement-type FEM in each element $E \in \Ec$, based on a dual formulation of local static problems $\Pc_E^{\Ec}$ (\ref{eq1:pblocal}). In practice, it is sufficient to consider a discretization of each element $E$ by a single element with an interpolation of degree $p+k$, where $p$ is the degree of interpolation of the initial FE analysis and $k$ an additional degree \cite{Bab94}. Numerical experiments tend to show that the use of an extra-degree $k=3$ enables one to obtain an approximate stress field of good quality \cite{Coo99}. Consequently, the practical resolution of local problems $\Pc_E^{\Ec}$ (\ref{eq1:pblocal}) leads to stress fields which are not exactly rigorous relative to the exact solution, but with respect to a refined solution.

\begin{remark}
Local problems $\Pc_E^{\Ec}$ (\ref{eq1:pblocal}) can be solved directly by searching an admissible stress field analytically in the polynomial form; the yielded stress field is strictly admissible only when the given body force field $\und{f}_d$ is a polynomial of degree compatible with that of $\hat{\cont}_{h \restrictto{E}}$; a step of decomposition of each element into subelements is required in order to ensure, at element vertices, the compatibility conditions resulting from the symmetry of the stress field $\hat{\cont}_{h \restrictto{E}}$ \cite{Lad04,Lad96}.
\end{remark}

The first step, which aims at constructing a set of equilibrated tractions, plays an important role in the quality of associated admissible stress fields and error estimates. 

\subsection{Original version of the construction of equilibrated tractions}\label{3.3}

Let us outline the main aspects related to the construction of equilibrated tractions. A detailed description and computational aspects of this method can be found in \cite{Lad10bis,Ple11}.

Starting from the strong prolongation condition (\ref{eq1:prolongfort}) rewritten in the global form:
\begin{equation}\label{eq2:prolongfort}
\begin{aligned}
& \intO \Tr\big[\left( \hat{\cont}_h - \cont_h \right) \: \defo(\und{v}_h^{\ast}) \big] \dO \\
& = \sum_{E \in \Ec} \intE \Tr\big[\left( \hat{\cont}_h - \cont_h \right) \: \defo(\und{v}_h^{\ast}) \big] \dO \\
& = 0 \quad \forall \: \und{v}_h^{\ast} \in \Vcb_h^p,
\end{aligned}
\end{equation}
where $\Vcb_h^p$ stands for the space of polynomial functions of degree $p$ which are continuous over each element $E \in \Ec$ and possibly discontinuous across inter-element edges, one can restrict (\ref{eq2:prolongfort}) to functions $\und{v}_h^{\ast} \in \Vcb_h^1$, since it is sufficient to satisfy equilibrium condition (\ref{eq1:eqfd}). Then, considering the weak form of the equilibrium equations verified by $\hat{\cont}_h$, tractions $\hat{\und{F}}_h$ satisfy: 
\begin{equation}\label{eq2:densites}
\begin{aligned}
& \sum_{E \in \Ec} \left[\intdE \eta_E \: \hat{\und{F}}_h \cdot \und{v}_h^{\ast} \dS \right. \\
& \left. - \int_E \left( \Tr\left[\cont_h \: \defo(\und{v}_h^{\ast})\right] - \und{f}_d \cdot \und{v}_h^{\ast} \right) \dO \right] = 0 \quad \forall \: \und{v}_h^{\ast} \in \Vcb_h^1.
\end{aligned}
\end{equation}

\begin{remark}
Global form (\ref{eq2:prolongfort}) of strong prolongation condition (\ref{eq1:prolongfort}) could have been modified by considering shape functions $\varphi_i \in \Uc_h^q$ ($1 \leqslant q \leqslant p$) in condition (\ref{eq1:prolongfort}), thus leading to condition (\ref{eq2:prolongfort}) with $\und{v}_h^{\ast} \in \Vcb_h^q$. For the sake of simplicity and practical purposes, one limits to the space $\Vcb_h^1$.
\end{remark}

Then, introduction of the partition of unity defined by the linear FE shape functions $\la_i \in \Uc_h^1$ into (\ref{eq2:densites}) yields the following system:
\begin{equation}\label{eq4:densites}
\begin{aligned}
& \sum_{E \in \Ec_i^{\Ic}} \left[ \intdE \eta_{E} \: \la_i \: \hat{\und{F}}_h \cdot \und{v}_h^{\ast} \dS \right. \\
& \left. - \int_E \left( \Tr\left[\cont_h \: \defo(\la_i \: \und{v}_h^{\ast})\right] - \und{f}_d \cdot \la_i \: \und{v}_h^{\ast} \right) \dO \right] = 0 \quad \forall \: \und{v}_h^{\ast} \in \Vcb_h^1.
\end{aligned}
\end{equation}

In order to confer more flexibility, we consider here the following set of local problems $\Pc_i^{\Ic}$ defined over the patch $\Om_i$ of elements $E \in \Ec_{i}^{\Ic}$ associated to each vertex $i \in \Ic$:

Find $\la_i \: \hat{\und{F}}_h^{(i)}$ such that:
\begin{equation}\label{eq2:pbdensitespatch}
\begin{aligned}
& \sum_{\Gamma \in \Jc_i^{\Ic}} \int_{\Gamma} \la_i \: \hat{\und{F}}_h^{(i)} \cdot \left( \sum_{E \in \Ec_{\Gamma}^{\Jc}} \eta_{E} \: {\und{v}_h^{\ast}}_{\restrictto{E}} \right) \dS \\
& = \Qc_{\Om_i}(\la_i \: \und{v}_h^{\ast}) \quad \forall \: \und{v}_h^{\ast} \in \Vcb_h^1,
\end{aligned}
\end{equation}
where
\begin{equation}
\displaystyle \Qc_{\Om_i}(\la_i \: \und{v}_h^{\ast}) = \intOmi \left( \Tr\left[\cont_h \: \defo(\la_i \: \und{v}_h^{\ast})\right] - \und{f}_d \cdot \la_i \: \und{v}_h^{\ast} \right) \dO.
\end{equation}

\begin{figure*}
\centering\includegraphics[width=0.9\textwidth]{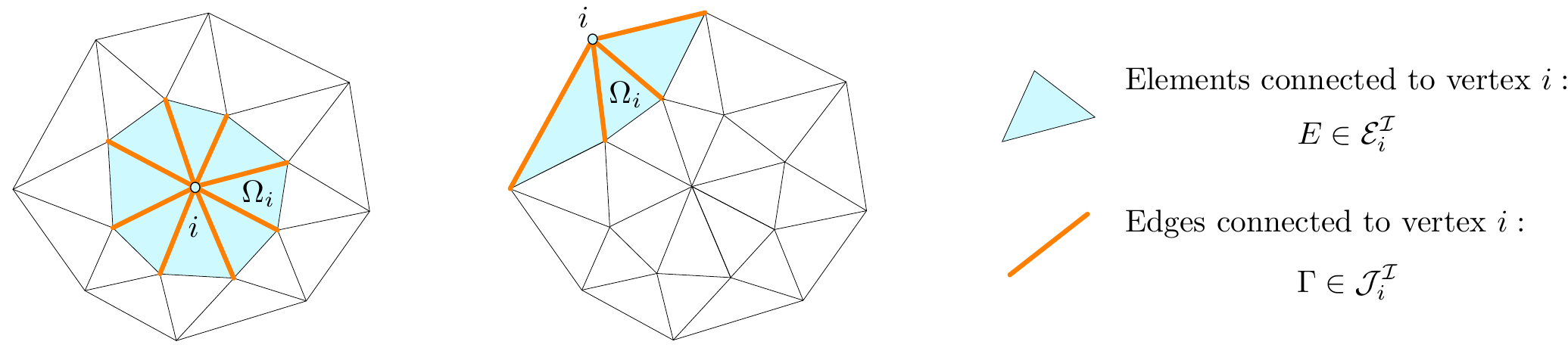}
\caption{Elements and edges connected to vertex $i$.}
\label{fig1:patchdouble_black_with_legend}
\end{figure*}

Indeed, it is worth recalling that quantity $\la_i \: \hat{\und{F}}_h^{(i)}$ is nonzero exclusively along edges $\Gamma \in \Jc_i^{\Ic}$. \Fig{fig1:patchdouble_black_with_legend} illustrates the sets of edges $\Gamma \in \Jc_i^{\Ic}$ and elements $\Ec \in \Ec_{i}^{\Ic}$ connected to vertex $i$.

The solvability and well-posedness of problems (\ref{eq2:pbdensitespatch}) is ensured for a FE interpolation degree $p \geqslant 2$ by considering the space ${\bar{\Ucb}_{h,0 \restrictto{\Om_i}}^1} = \left\{ \und{v}_h^{\ast} \in \bar{\Vcb}_h^1, {\und{v}_h^{\ast}}_{\restrictto{\Gamma \in \Jc_i^{\Ic} \cap \dOm}} = \und{0} \right\}$, where $\bar{\Vcb}_h^1$ defines the set of piecewise linear polynomial functions $\und{v}_h^{\ast} \in \Vcb_h^1$ which are continuous across edges $\Gamma \in \Jc_i^{\Ic}$, and therefore continuous over the whole patch $\Om_i$. A specific treatment is required in the case $p = 1$, for which problems (\ref{eq2:pbdensitespatch}) are substituted by:

Find $\la_i \: \hat{\und{F}}_h^{(i)}$ such that:
\begin{equation}\label{eq3:pbdensitespatch}
\begin{aligned}
& \sum_{\Gamma \in \Jc_i^{\Ic}} \int_{\Gamma} \la_i \: \hat{\und{F}}_h^{(i)} \cdot \left( \sum_{E \in \Ec_{\Gamma}^{\Jc}} \eta_{E} \: {\und{v}_h^{\ast}}_{\restrictto{E}} \right) \dS \\
& = \Qc_{\Om_i}(\la_i \: \und{v}_{h}^{\ast}(\und{x}_i)) \quad \forall \: \und{v}_h^{\ast} \in \Vcb_h^1,
\end{aligned}
\end{equation}
where
\begin{equation}
\begin{aligned}
& \displaystyle\Qc_{\Om_i}(\la_i \: \und{v}_{h}^{\ast}(\und{x}_i)) \\
& = \intOmi \left( \Tr\left[\cont_h \: \defo(\la_i \: \und{v}_h^{\ast}(\und{x}_i))\right] - \und{f}_d \cdot \la_i \: \und{v}_h^{\ast}(\und{x}_i) \right) \dO.
\end{aligned}
\end{equation}

One can demonstrate the existence of a solution to problems (\ref{eq3:pbdensitespatch}), by using both subspace ${\bar{\Ucb}_{h,0 \restrictto{\Om_i}}^1}$ and FE equilibrium. (see \cite{Lad10bis} for more details).

Uniqueness of the solution of such problems is guaranteed by the least-squares minimization of a cost function of the form \cite{Lad10bis}:
\begin{equation}\label{eq1:fonctioncoutEESPT}
J_{\Om_i}(\la_i \hat{\und{F}}_h^{(i)}) = \frac{1}{2} \sum_{\Gamma \in \Jc_i^{\Ic}} ( \la_i \: \hat{\und{F}}_h^{(i)} - \la_i \: \und{F}_h^{(i)} )_{\restrictto{\Gamma}}^2,
\end{equation}
which represents the gap between the searched solution $\la_i \: \hat{\und{F}}_{h \restrictto{\Gamma}}^{(i)}$ and the known quantity $\la_i \: \und{F}_{h \restrictto{\Gamma}}^{(i)}$ involving the projection of the FE stress field $\cont_h$ over the edge $\Gamma \in \Jc_i^{\Ic}$ and the traction force density $\und{F}_d$.

Eventually, one recovers tractions $\hat{\und{F}}_h$ along each edge $\Gamma \in \Jc$ directly from calculated quantities $\la_i \: \hat{\und{F}}_h^{(i)}$, which are sought in $\Ucb_{h \restrictto{\Gamma}}^p$, in such a way that $\hat{\und{F}}_h \in \Ucb_{h \restrictto{\Gamma}}^p$:
\begin{equation}\label{eq1:interpolationFh}
\hat{\und{F}}_{h \restrictto{\Gamma}} = \sum_{i \in \Ic_{\Gamma}^{\Jc}} ( \la_i \: \hat{\und{F}}_h^{(i)} )_{\restrictto{\Gamma}}.
\end{equation}

Besides, enforcement of conditions $\eta_E \: \hat{\und{F}}_h = \und{F}_d$ over edges $\Gamma \subset \partial_2 \Om$ can be achieved by adding these constraints in the constrained minimization problem.

In the following part, we describe the enhanced version in details with the purpose of emphasizing the key aspects of the method. 

\section{Principles of the enhanced version of the EESPT technique}\label{4}

The basic idea is to confer more flexibility and to give greater freedom in the construction of equilibrated tractions in order to improve the quality of associated admissible stress fields despite higher (but reasonable) computational cost. Indeed, in zones of high element aspect ratios or sharp gradients, the quality of admissible stress fields may be affected, thus resulting in large effectivity indices \cite{Lad97,Lad04}. This observation has spurred the development of an enhanced construction of equilibrated tractions. The principle, originally introduced in \cite{Lad97}, is to optimize the quality of the computed admissible stress field by improving the recovering strategy for the construction of equilibrated tractions. Thus, the construction of such tractions has been changed and is henceforth based on a weakened prolongation condition, which amounts to removing, from the strong prolongation condition, the shape functions associated with vertex nodes. Therefore, equilibrated stress field $\hat{\cont}_h$ is still constructed as an extension (or prolongation) of the FE stress field $\cont_h$, but the weak extension concerns only the non-vertex nodes of the FE mesh. From a practical point of view, this optimized construction of the tractions can be applied locally only in relevant zones in order to preserve an affordable computational cost. In this work, the main breakthrough is the use of sound criteria which enables to select these apposite regions. Let us define $\Ec_e \subset \Ec$ and $\Jc_e \subset \Jc$ the sets of elements and edges involved in the enhanced procedure, respectively. Several criteria have been considered to select part $\Om_e$, \ie the set $\Ec_e$ of elements involved in the global minimization step (see \Sect{4.2}).

\subsection{Enhanced version of the construction of equilibrated tractions}\label{4.1}

The prolongation condition needed for the construction of equilibrated tractions along $\Gamma \in \Jc_e$ is reduced to:
\begin{equation}\label{eq1:prolongfaible}
\intE \left(\hat{\cont}_h - \cont_h \right) \: \und{\nabla} \varphi_i \dO = \und{0} \quad \forall \: i \in \Nc_E^{\Ec} \setminus \Ic_E^{\Ec}, \ \forall \: E \in \Ec_e,
\end{equation}
where $\varphi_i$ is the FE shape function associated with non-vertex node $i$.

This relation is weaker than the strong prolongation condition as it involves only the FE shape functions of higher degree. In the following, the shape functions are described in hierarchical form in such a way that the linear part of the FE shape functions is associated with vertices.

As for the standard construction, densities $\hat{\und{F}}_h$ along edges $\Gamma \in \Jc_e$ are searched in a discretized space with the same interpolation degree as the FE displacement field $\uu_h$, \ie $\displaystyle\hat{\und{F}}_{h \restrictto{\Gamma}} \in \Ucb^{p}_{h \restrictto{\Gamma}}$.
 
The modification of the prolongation condition leads to a partition of the tractions:
\begin{equation}\label{eq1:decompF}
\hat{\und{F}}_h = \hat{\und{H}}_h + \hat{\und{R}}_h \quad \textrm{on} \ \Gamma \in \Jc_e,
\end{equation}
with
\begin{align}
& \int_{\Gamma} \hat{\und{H}}_h \: \varphi_i \dS = \und{0} \quad \forall \: i \in \Ic_{\Gamma}^{\Jc_e}, \label{eq1:decompH}\\
& \int_{\Gamma} \hat{\und{R}}_h \: \varphi_i \dS = \und{0} \quad \forall \: i \in \Nc_{\Gamma}^{\Jc_e} \setminus \Ic_{\Gamma}^{\Jc_e} \label{eq1:decompR}.
\end{align}

This decomposition is unique \cite{Lad04} and one can note that part $\hat{\und{H}}_h$ has zero resultant and moment on $\Gamma \in \Jc_e$. Now, let us focus on the determination of each part of tractions.\\
The determination of part $\hat{\und{H}}_h$ on $\Gamma \in \Jc_e$ is completely governed by the weak prolongation condition (\ref{eq1:prolongfaible}) and relation (\ref{eq1:decompH}), as the strong prolongation condition is no longer respected. As a result, part $\hat{\und{H}}_h$ on $\Gamma \in \Jc_e$ only depends on the data of the problem and the FE stress field. The construction of this part of the tractions is globally similar to the standard construction, thus it requires local calculations and therefore is not expensive in terms of computational time. Thus, part $\hat{\und{H}}_h$ can be determined explicitly.\\

The determination of part $\hat{\und{R}}_h$ on $\Gamma \in \Jc_e$ is performed by minimizing the complementary energy (or, equivalently, the constitutive relation error) locally on part $\Om_e \subset \Om$ containing elements $E \in \Ec_e$ under the following constraints:
\begin{itemize}
\item[$\bullet$] Neumann boundary conditions over edges $\Gamma \in \Jc_e \cap \partial_2 \Om$;
\item[$\bullet$] equilibrium conditions of tractions $\hat{\und{H}}_h + \hat{\und{R}}_h$ with body force field $\und{f}_d$ over each element $E \in \Ec_e$;
\item[$\bullet$] equilibrium conditions of tractions $\hat{\und{H}}_h + \hat{\und{R}}_h$ with body force field $\und{f}_d$ and standard tractions $\hat{\und{F}}^{\std}_h$ over each element $E \in \bar{\Ec}_e \setminus \Ec_e$,
\end{itemize}
where $\bar{\Ec}_e \subset \Ec$ denotes the set of elements $E$ connected to at least one edge $\Gamma \in \Jc_e$; therefore, $\bar{\Ec}_e \setminus \Ec_e$ contains all the elements connected to one and only one edge $\Gamma \in \Jc_e$; $\hat{\und{F}}^{\std}_h$ are pre-computed tractions over edges $\Gamma \in \dE \setminus \Jc_e$ coming from the standard construction over element $E \in \bar{\Ec}_e \setminus \Ec_e$. Let us recall that standard tractions are constructed in equilibrium with the external loading over edges $\Gamma \in \Jc \setminus \Jc_e$. Indeed, equilibrium conditions (\ref{eq1:eqFd}) and (\ref{eq1:eqfd}) are inherently enforced only over elements $E \in \Ec \setminus \bar{\Ec}_e$ and edges $\Gamma \in \Jc \setminus \Jc_e$ in the standard construction of equilibrated tractions.\\

Besides, one can reduce the computational cost resulting from the calculation of this global problem by introducing two local problems, referred to as $\Pc^H_{\vert E}$ and $\Pc^R_{\vert E}$, defined on each element $E \in \Ec_e$ and linked to parts $\hat{\und{H}}_h$ and $\hat{\und{R}}_h$, respectively \cite{Lad04}. Similarly, a local problem $\Pc^{H,F}_{\vert E}$ linked to $\hat{\und{H}}_h$ and $\hat{\und{F}}^{\std}_h$ is introduced on each element $E \in \bar{\Ec}_e \setminus \Ec_e$. Let us decompose $\hat{\cont}_h$ into two parts $\hat{\cont}^H_h$ and $\hat{\cont}^R_h$ such that:
\begin{equation}\label{eq1:decompsigma}
\hat{\cont}_h = \hat{\cont}^H_h + \hat{\cont}^R_h \quad \textrm{on} \ E \in \Ec_e,
\end{equation}
with
\begin{align}
& \hat{\cont}^H_{h \restrictto{E}} \: \und{n}_E = \eta_E \: \hat{\und{H}}_h \quad \textrm{on} \ \dE, \label{eq1:decompsigmaH}\\
& \hat{\cont}^R_{h \restrictto{E}} \: \und{n}_E = \eta_E \: \hat{\und{R}}_h \quad \textrm{on} \ \dE \label{eq1:decompsigmaR}.
\end{align}

Let us now consider the following local problems:
\begin{itemize}
\item[$\bullet$] Problem $\Pc^H_{\vert E}$:
\begin{equation}\label{eq1:pblocalH}
\hat{\cont}^H_{h \restrictto{E}} \in \Scb_{\hat{\und{H}}_h} \iff \left\{
\begin{aligned}
& \hat{\cont}^H_{h \restrictto{E}} \in \Scb \\
& \und{\diver} \: \hat{\cont}^H_{h \restrictto{E}} + \und{f}^H_E = \und{0} \quad \textrm{in} \ E \\
& \hat{\cont}^H_{h \restrictto{E}} \: \und{n}_E = \eta_E \: \hat{\und{H}}_h \quad \textrm{on} \ \dE,
\end{aligned}
\right.
\end{equation}
with
\begin{equation}
\begin{aligned}
\circ \ \und{f}^H_E & = - \displaystyle\frac{1}{\labs{E}} \intdE \eta_E \: \hat{\und{H}}_h \dS \\
& - \left( I^{-1}_G \intdE \left( \GM \wedge \eta_E \: \hat{\und{H}}_h \right) \dS \cdot \und{N} \right) \und{N} \wedge \GM
\end{aligned}
\end{equation}
in two dimensions, where $\labs{E}$ represents the measure of element $E$ and $\und{N}$ denotes the axis normal to the two-dimensional plane considered; $I_G$ is the scalar mass moment of inertia around axis $\und{N}$ passing through the center of mass $G$;
\begin{equation}
\begin{aligned}
\circ \ \und{f}^H_E & = - \displaystyle\frac{1}{\labs{E}} \intdE \eta_E \: \hat{\und{H}}_h \dS \\
& - \left( \boldsymbol{I}^{-1}_G \intdE \left( \GM \wedge \eta_E \: \hat{\und{H}}_h \right) \dS \right) \wedge \GM
\end{aligned}
\end{equation}
in three dimensions, where $\boldsymbol{I}_G$ is the mass moment of inertia tensor with respect to the center of mass $G$. 

$\und{f}^H_E$ is constructed in such a way that $\und{f}^H_E$ is in equilibrium with part $\hat{\und{H}}_h$. The weak form of problem $\Pc^H_{\vert E}$ reads:

Find $\hat{\cont}^H_{h \restrictto{E}} \in \Scb$ such that:
\begin{equation}\label{PH}
\begin{aligned}
\forall \: \uu^{\ast} \in \Ucb_{0 \restrictto{E}}, \quad \intE \Tr\big[\hat{\cont}^H_{h \restrictto{E}} \: \defo(\uu^{\ast})\big] \dO \\
= \intE \und{f}^H_E \cdot \uu^{\ast} \dO + \intdE \eta_E \: \hat{\und{H}}_h \cdot \uu^{\ast} \dS;
\end{aligned}
\end{equation}
Solving problems $\Pc^H_{\vert E}$ over each element $E \in \Ec_e$ leads to part $\hat{\cont}^H_{h \restrictto{E}}$, which is determined explicitly from $\hat{\und{H}}_h$.\\
\item[$\bullet$] Problem $\Pc^R_{\vert E}$:
\begin{equation}\label{eq1:pblocalR}
\hat{\cont}^R_{h \restrictto{E}} \in \Scb_{\hat{\und{R}}_h} \iff \left\{
\begin{aligned}
& \hat{\cont}^R_{h \restrictto{E}} \in \Scb \\
& \und{\diver} \: \hat{\cont}^R_{h \restrictto{E}} + \und{f}^R_E = \und{0} \quad \textrm{in} \ E \\
& \hat{\cont}^R_{h \restrictto{E}} \: \und{n}_E = \eta_E \: \hat{\und{R}}_h \quad \textrm{on} \ \dE,
\end{aligned}
\right.
\end{equation}
with
\begin{equation}
\begin{aligned}
\circ \ \und{f}^R_E & = - \displaystyle\frac{1}{\labs{E}} \intdE \eta_E \: \hat{\und{R}}_h \dS \\
& - \left( I^{-1}_G \intdE \left( \GM \wedge \eta_E \: \hat{\und{R}}_h \right) \dS \cdot \und{N} \right) \und{N} \wedge \GM
\end{aligned}
\end{equation}
in two dimensions;
\begin{equation}
\begin{aligned}
\circ \ \und{f}^R_E & = - \displaystyle\frac{1}{\labs{E}} \intdE \eta_E \: \hat{\und{R}}_h \dS \\
& - \left( \boldsymbol{I}^{-1}_G \intdE \left( \GM \wedge \eta_E \: \hat{\und{R}}_h \right) \dS \right) \wedge \GM
\end{aligned}
\end{equation}
in three dimensions.

As a result, $\und{f}^R_E$ is in equilibrium with part $\hat{\und{R}}_h$. The weak form of problem $\Pc^R_{\vert E}$ reads:

Find $\hat{\cont}^R_{h \restrictto{E}} \in \Scb$ such that:
\begin{equation}\label{PR}
\begin{aligned}
\forall \: \uu^{\ast} \in \Ucb_{0 \restrictto{E}}, \quad \intE \Tr\big[\hat{\cont}^R_{h \restrictto{E}} \: \defo(\uu^{\ast})\big] \dO \\
= \intE \und{f}^R_E \cdot \uu^{\ast} \dO + \intdE \eta_E \: \hat{\und{R}}_h \cdot \uu^{\ast} \dS;
\end{aligned}
\end{equation}
Solving problems $\Pc^R_{\vert E}$ over each element $E \in \Ec_e$ leads to a linear relation $\hat{\cont}^R_{h \restrictto{E}} (\hat{\und{R}}_{h \restrictto{\dE}})$, as $\und{f}^R_E$ is a linear function with respect to $\hat{\und{R}}_{h \restrictto{\dE}}$.
\end{itemize}

Subsequently, the global constrained minimization problem defined over $\Om_e$ consists of minimizing either the complementary energy on part $\Om_e$:
\begin{equation}\label{eq1:Ec}
\begin{aligned}
\frac{1}{2} & \intOme \Tr\big[ \hat{\cont}^R_h(\hat{\und{R}}_h) \: \K^{-1} \: \hat{\cont}^R_h(\hat{\und{R}}_h)\big] \dO \\
+ & \intOme \Tr\big[ \hat{\cont}^R_h(\hat{\und{R}}_h) \: \K^{-1} \: \hat{\cont}^H_h\big] \dO \\
- & \int_{\partial_1 \Om \cap \Jc_e} \hat{\cont}^R_h(\hat{\und{R}}_h) \: \und{n} \cdot \Uu_d \dS,
\end{aligned}
\end{equation}
or, equivalently, the constitutive relation error:
\begin{equation}\label{eq1:Erdc}
\begin{aligned}
\frac{1}{2} & \intOme \Tr\big[ \hat{\cont}^R_h(\hat{\und{R}}_h) \: \K^{-1} \: \hat{\cont}^R_h(\hat{\und{R}}_h)\big] \dO \\
+ & \intOme \Tr\big[ \hat{\cont}^R_h(\hat{\und{R}}_h) \: \K^{-1} \: \left( \hat{\cont}^H_h - \cont_h \right)\big] \dO,
\end{aligned}
\end{equation}
under the following constraints:
\begin{align}
 \bullet \ & \eta_E \: \hat{\und{R}}_h = \und{F}_d - \eta_E \: \hat{\und{H}}_h \quad \textrm{on edges} \ \Gamma \in \partial_2 \Om \cap \Jc_e, \label{eq2:eqFd} \\ 
 \bullet \ & \forall \: \us \in \Ucb_{S \restrictto{E}}, \ \forall E \in \Ec_e, \quad \intdE \eta_E \: \hat{\und{R}}_h \cdot \us \dS \nonumber \\
 = & - \intdE \eta_E \: \hat{\und{H}}_h \cdot \us \dS - \intE \und{f}_d \cdot \us \dO \label{eq2:eqfd} \\
 \bullet \ & \forall \: \us \in \Ucb_{S \restrictto{E}}, \ \forall E \in \bar{\Ec}_e \setminus \Ec_e, \quad \int_{\dE \cap \Jc_e} \eta_E \: \hat{\und{R}}_h \cdot \us \dS \nonumber \\
 = & - \int_{\dE \setminus \Jc_e} \eta_E \: \hat{\und{F}}^{\std}_h \cdot \us \dS \nonumber \\
 & - \int_{\dE \cap \Jc_e} \eta_E \: \hat{\und{H}}_h \cdot \us \dS - \intE \und{f}_d \cdot \us \dO \label{eq3:eqfd}
\end{align}

Constraint (\ref{eq2:eqFd}) enforces the equilibrium of tractions $\hat{\und{H}}_h + \hat{\und{R}}_h$ with the external traction force density $\und{F}_d$ on element edges $\Gamma \in \Jc_e \cap \partial_2 \Om$. Constraint (\ref{eq2:eqfd}) enforces the equilibrium of tractions $\hat{\und{H}}_h + \hat{\und{R}}_h$ on element edges $\Gamma \in \Jc_e$ with the external body force field $\und{f}_d$ over elements $E \in \Ec_e$. Eventually, constraint (\ref{eq3:eqfd}) enforces the equilibrium of tractions $\hat{\und{H}}_h + \hat{\und{R}}_h$ on element edge $\Gamma \in \Jc_e$ with $\und{f}_d$ and the standard tractions $\hat{\und{F}}^{\std}_h$ over elements $E \in \bar{\Ec}_e \setminus \Ec_e$.\\

The resulting stress field $\hat{\cont}_{h \restrictto{E}} = \hat{\cont}^H_{h \restrictto{E}} + \hat{\cont}^R_{h \restrictto{E}} (\hat{\und{R}}_{h \restrictto{\dE}})$, $E \in \Ec_e$ is statically admissible provided that constraints (\ref{eq2:eqFd}), (\ref{eq2:eqfd}) and (\ref{eq3:eqfd}) hold. Since $\und{f}^R_E$ and $\und{f}^H_E$ are built in equilibrium with parts $\hat{\und{R}}_{h \restrictto{\dE}}$ and $\hat{\und{H}}_{h \restrictto{\dE}}$, respectively, equilibrium conditions (\ref{eq2:eqfd}) can be rewritten in the form:
\begin{align}
& \intE \und{f}^R_E \dO = \intE\left( \und{f}_d - \und{f}^H_E \right) \dO \quad \forall E \in \Ec_e; \label{eq4:eqfd}\\
& \intE \GM \wedge \und{f}^R_E \dO = \intE \GM \wedge \left( \und{f}_d - \und{f}^H_E \right) \dO \quad \forall E \in \Ec_e. \label{eq5:eqfd}
\end{align}
Indeed, (\ref{eq4:eqfd}) and (\ref{eq5:eqfd}) convey the equilibrium between tractions $\hat{\und{H}}_h + \hat{\und{R}}_h$ (or body force field $\und{f}^R_E + \und{f}^H_E$, equivalently) and $\und{f}_d$ over elements $E \in \Ec_e$. Likewise, equilibrium conditions (\ref{eq3:eqfd}) over $E \in \bar{\Ec}_e \setminus \Ec_e$ can be recast in a form similar to relations (\ref{eq4:eqfd}) and (\ref{eq5:eqfd}).

Eventually, part $\hat{\und{R}}_h$ is recovered over edges $\Gamma \in \Jc_e$ by solving two local problems $\Pc^H_{\vert E}$ and $\Pc^R_{\vert E}$ defined at the element scale and one global problem defined at the part $\Om_e$ scale.

\subsection{Criteria used in the enhanced procedure}\label{4.2}

Three criteria have been considered for the selection of relevant zones involved in the optimization procedure:
\begin{itemize}
\item[$\bullet$] the radius ratio, which is the radius of the inscribed circle tangent to triangle's three edges divided by the radius of the circle circumscribed by triangle's three vertices for two-dimensional cases whose FE mesh is made of linear triangular elements; besides, it is the ratio between the radius of the inscribed sphere tangent to tetrahedron's four faces divided by the radius of the sphere circumscribed by tetrahedron's four vertices for three-dimensional cases whose FE mesh is made of linear tetrahedral elements;
\item[$\bullet$] the edge ratio, which is the ratio between the shortest element edge length and the longest element edge length for two-dimensional cases; similarly, the area ratio, which is the ratio between the smallest element face area and the largest element face area for three-dimensional cases; 
\item[$\bullet$] the estimate ratio, which is the ratio between the squared element-by-element contribution to the global estimate and the squared maximal local contribution to the global estimate.
\end{itemize}

The two geometric criteria, namely the radius ratio and edge ratio, are intended to take the most distorted elements into account in the optimization procedure, while the error estimate criterion allows the elements which contribute to the estimate $\theta$ the most to be selected. In \Sect{6}, we limit ourselves to alone the three aforementioned criteria.

Implementation issues regarding the enhanced procedure are discussed in the next section. The interested reader can refer to \cite{Lad10bis,Ple11} for more details about practical implementation for the standard versions of EET and EESPT methods. 

\section{Practical Implementation of the optimized procedure}\label{5}

\subsection{Discretized quantities}\label{5.1}

Over an edge $\Gamma \in \Jc_e$, searched quantities $\hat{\und{R}}_{h \restrictto{\Gamma}}$ and $\hat{\und{H}}_{h \restrictto{\Gamma}}$ are discretized in the form:
\begin{align}
& \hat{\und{R}}_{h \restrictto{\Gamma}} = [\varphi_{\restrictto{\Gamma}}]^{T} \: \hat{\und{r}}_{h, \Gamma}; \label{eq1:discretisation_Rh} \\
& \hat{\und{H}}_{h \restrictto{\Gamma}} = [\varphi_{\restrictto{\Gamma}}]^{T} \: \hat{\und{h}}_{h, \Gamma}, \label{eq1:discretisation_Hh}
\end{align}
where $[\varphi_{\restrictto{\Gamma}}]$ is the matrix of FE shape functions of degree $p$ over $\Gamma$ and $\hat{\und{r}}_{h, \Gamma}$ (respectively $\hat{\und{h}}_{h, \Gamma}$) is the unknown vector of components of $\hat{\und{R}}_{h \restrictto{\Gamma}}$ (respectively $\hat{\und{H}}_{h \restrictto{\Gamma}}$). Let $\hat{\und{r}}_{h}$ (respectively $\hat{\und{h}}_h$) be the generalized vector corresponding to the combination of unknown vectors $\hat{\und{r}}_{h, \Gamma}$ ($\hat{\und{h}}_{h, \Gamma}$ respectively) for every edge $\Gamma \in \Jc_e$.

As already pointed out in \Sect{4.1}, $\hat{\und{h}}_h$ is explicitly calculated using weak prolongation condition (\ref{eq1:prolongfaible}) and relation (\ref{eq1:decompH}). From now on, $\hat{\und{h}}_h$ is assumed to be known.

\subsection{Practical resolution}\label{5.2}

In practice, local problems (\ref{PH}) and (\ref{PR}) are solved numerically in the same way as problem (\ref{eq1:pblocal}), that is, by duality on element $E$, a displacement-type FEM by considering the original FE mesh $\Mc_h$ with a $p + 3$ discretization over each element $E$.

Minimization function (\ref{eq1:Erdc}) (or (\ref{eq1:Ec})) takes the following matrix form:
\begin{equation}\label{eq1:Erdcdiscr}
\frac{1}{2} \: \hat{\und{r}}_{h}^{T} \: \Abb \: \hat{\und{r}}_{h} - \hat{\und{r}}_{h}^{T} \: \und{B},
\end{equation}
where $\Abb$ is a symmetric matrix.

Constraints (\ref{eq2:eqFd}), (\ref{eq4:eqfd}), (\ref{eq5:eqfd}) and (\ref{eq1:decompR}) are back imposed through the definition of extra sets of unknowns, the so-called Lagrange multipliers. First, equilibrium conditions (\ref{eq2:eqFd}) over edges $\Gamma \in \Jc_e \cap \partial_2 \Om$ read:
\begin{equation}\label{eq1:eqimpC}
\Cbb \: \hat{\und{r}}_h = \und{q}.
\end{equation}

Second, equilibrium conditions (\ref{eq4:eqfd}) and (\ref{eq5:eqfd}) over elements $E \in \Ec_e$ and similar equilibrium conditions over elements $E \in \bar{\Ec}_e \setminus \Ec_e$ read:
\begin{equation}\label{eq1:eqimpL}
\Lbb \: \hat{\und{r}}_h = \und{b}.
\end{equation}
Let us note that solvability of (\ref{eq1:eqimpL}) requires the verification of compatibility conditions resulting from the FE equilibrium. One way to overcome this problem is to properly eliminate the redundant equations. A zero or small pivots detection procedure has been performed to handle and discard redundant equations into (\ref{eq1:eqimpL}) and yields a reduced system which reads:
\begin{equation}\label{eq2:eqimpL}
\tilde{\Lbb} \: \hat{\und{r}}_h = \tilde{\und{b}}.
\end{equation}

Third, conditions (\ref{eq1:decompR}) over edges $\Gamma \in \Jc_e$ read:
\begin{equation}\label{eq1:eqimpP}
\Pbb \: \hat{\und{r}}_h = \und{0}.
\end{equation}
Let us notice that conditions (\ref{eq1:decompR}) (or system (\ref{eq1:eqimpP})) vanish in the case $p = 1$, $p$ being the FE interpolation degree.

Therefore, introducing the Lagrangian:
\begin{equation}\label{eq1:lagenhanced}
\begin{aligned}
L(\hat{\und{r}}_h,\und{\Lambda}_{\Cbb},\und{\Lambda}_{\Pbb},\und{\Lambda}_{\Lbb}) = & \frac{1}{2} \: \hat{\und{r}}_{h}^{T} \: \Abb \: \hat{\und{r}}_{h} - \hat{\und{r}}_{h}^{T} \: \und{B} + ( \Cbb \: \hat{\und{r}}_h - \und{q} )^T \: \und{\Lambda}_{\Cbb} \\
& + ( \Pbb \: \hat{\und{r}}_h )^T \: \und{\Lambda}_{\Pbb} + ( \tilde{\Lbb} \: \hat{\und{r}}_h - \tilde{\und{b}} )^T \: \und{\Lambda}_{\Lbb},
\end{aligned}
\end{equation}
the system to be solved takes the matrix form:
\begin{equation}\label{Pbminenhanced}
\begin{pmatrix}
\Abb & \Cbb^T & \Pbb^T & \tilde{\Lbb}^T \\
\Cbb & 0 & 0 & 0 \\
\Pbb & 0 & 0 & 0 \\
\tilde{\Lbb} & 0 & 0 & 0 \\
\end{pmatrix}
\begin{bmatrix}
\hat{\und{r}}_h \\
\und{\Lambda}_{\Cbb} \\
\und{\Lambda}_{\Pbb} \\
\und{\Lambda}_{\Lbb} \\
\end{bmatrix}
=
\begin{bmatrix}
\und{B} \\
\und{q} \\
\und{0} \\
\tilde{\und{b}} \\
\end{bmatrix},
\end{equation}
where $\und{\Lambda}_{\Cbb}$, $\und{\Lambda}_{\Pbb}$, $\und{\Lambda}_{\Lbb}$ represent the vectors of Lagrange multipliers associated with constraints (\ref{eq1:eqimpC}), (\ref{eq1:eqimpP}) and (\ref{eq2:eqimpL}), respectively.

\section{Numerical results}\label{6}

All the two- and three-dimensional structures considered are made of an isotropic linear elastic material. The two-dimensional cases are plane-stress problems. Values for Young's modulus and Poisson's ratio are set to $1$ and $0.3$, respectively. All the calculations of local solutions are performed using an interpolation of degree $p + k$ with $k = 3$ ($p$-refinement). The performance of the two error estimates, namely EET and EESPT, is analyzed with respect to the geometric and error estimate criteria in challenging real industrial applications. The broad applicability and relevance of the proposed techniques are investigated through the following two- and three-dimensional model problems:
\begin{itemize}
\item[$\bullet$] a two-dimensional cracked structure, already studied in \cite{Par06,Lad10bis}, which contains two different round cavities;
\item[$\bullet$] a two-dimensional weight sensor under bending already considered in a previous paper \cite{Ple11};
\item[$\bullet$] a three-dimensional open hole specimen under tension.
\end{itemize}

All those numerical experiments are performed by using FE meshes made of linear elements. Let us recall that both EET and EESPT methods give similar results in this particular case, as only their implementations are different.

\subsection{Cracked structure}\label{6.1}
\begin{figure}
\centering\includegraphics[scale = 0.36]{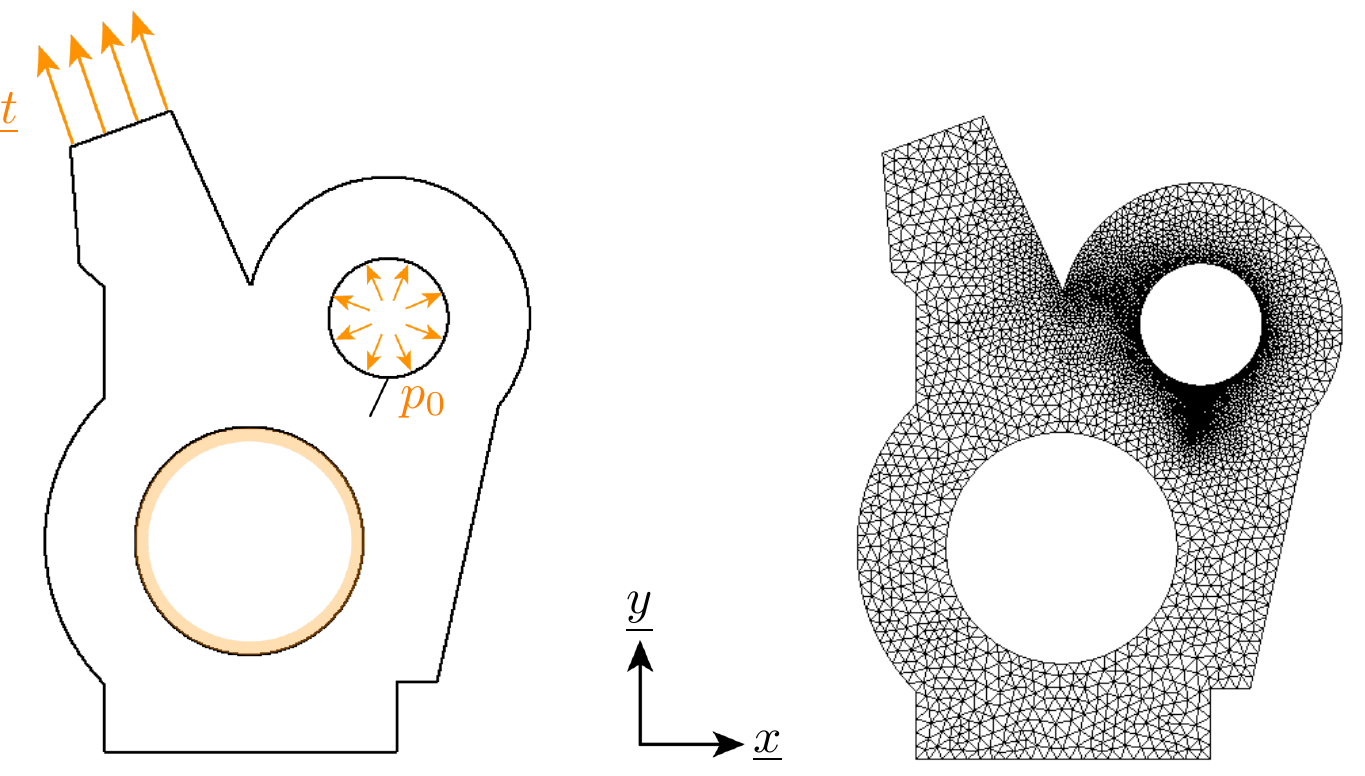}
\caption{Cracked structure model problem (left) and associated finite element mesh (right).}
\label{fig1:structure_fissuree_2D_geometry_fine_mesh}
\end{figure}

Let us consider a two-dimensional cracked structure, represented in \Fig{fig1:structure_fissuree_2D_geometry_fine_mesh}, which contains two circular holes. Under the plane stress assumption, the central round cavity is clamped, while the smaller one is subjected to an internal constant pressure $p_0$. The structure is also subjected to a unit normal traction force density $\und{t} = + \und{n}$ acting on the top-left side. Besides, a crack is located at the bottom of the second cavity. Therefore, a singularity is located at the crack tip where a strong stress concentration occurs. Traction-free boundary conditions are prescribed over the two lips of the crack and over the other sides. The mesh density increases toward the singularity introduced by the crack tip, which is the highest stress zone. The FE mesh containing $7 \, 751$ linear triangular elements and $4 \, 122$ nodes (\ie $8 \, 244$ d.o.f.) is given in \Fig{fig1:structure_fissuree_2D_geometry_fine_mesh}. The reference solution, also called \textquotedblleft quasi-exact\textquotedblright \, solution, is obtained through an overkill FE calculation with a reference mesh made of $1 \, 082 \, 011$ linear triangular elements and $543 \, 744$ nodes (\ie $1 \, 087 \, 488$ d.o.f.).
\begin{figure}
\centering\includegraphics[scale = 0.33]{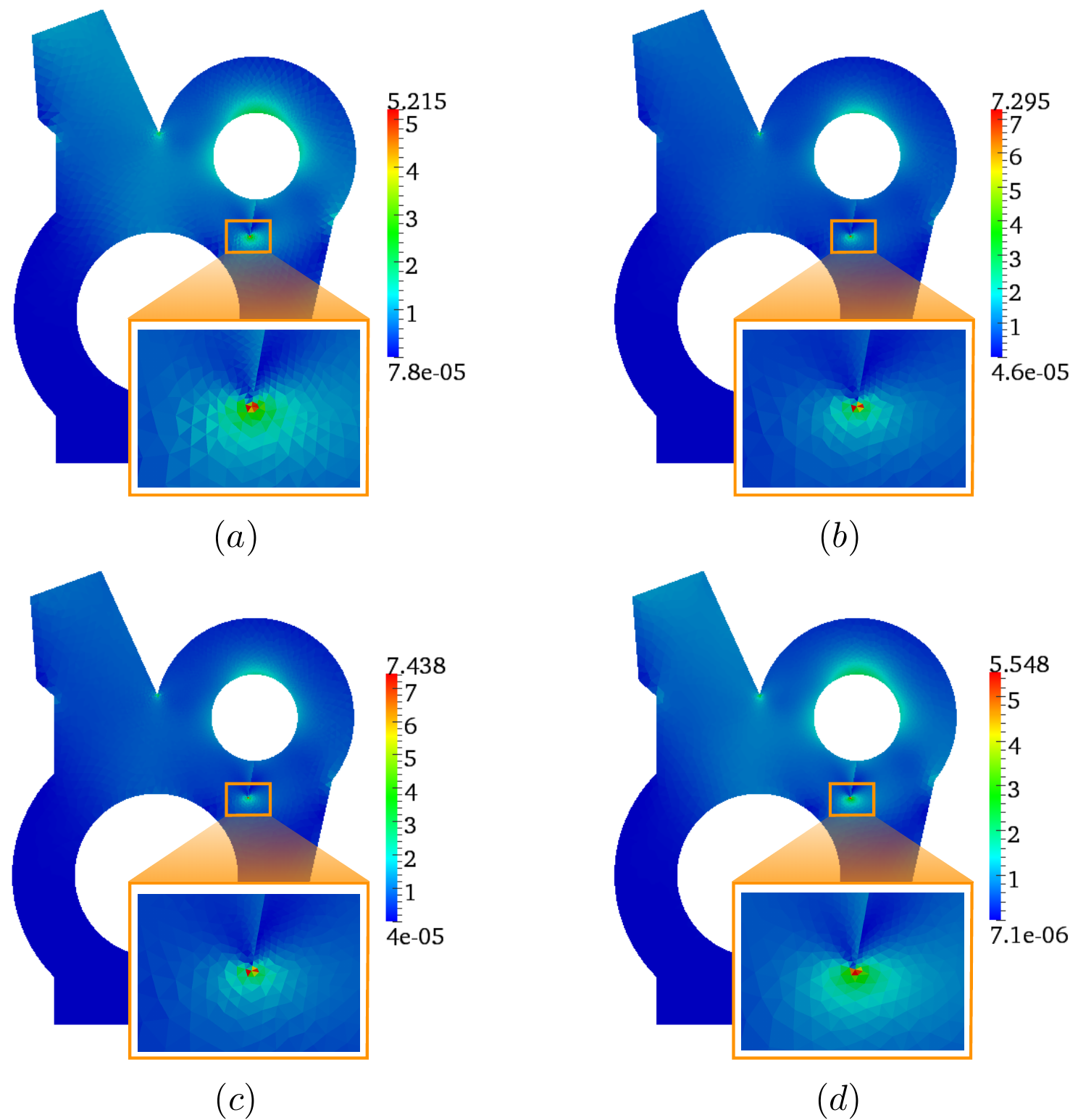}
\caption{Magnitude of the FE stress field (a) and the admissible stress field calculated using the standard versions of the EET (b) and the EESPT (c), and the full enhanced version (d).}
\label{fig1:sigma_EF_hat_structure_fissuree_2D}
\end{figure}
\begin{figure}
\centering\includegraphics[scale = 0.33]{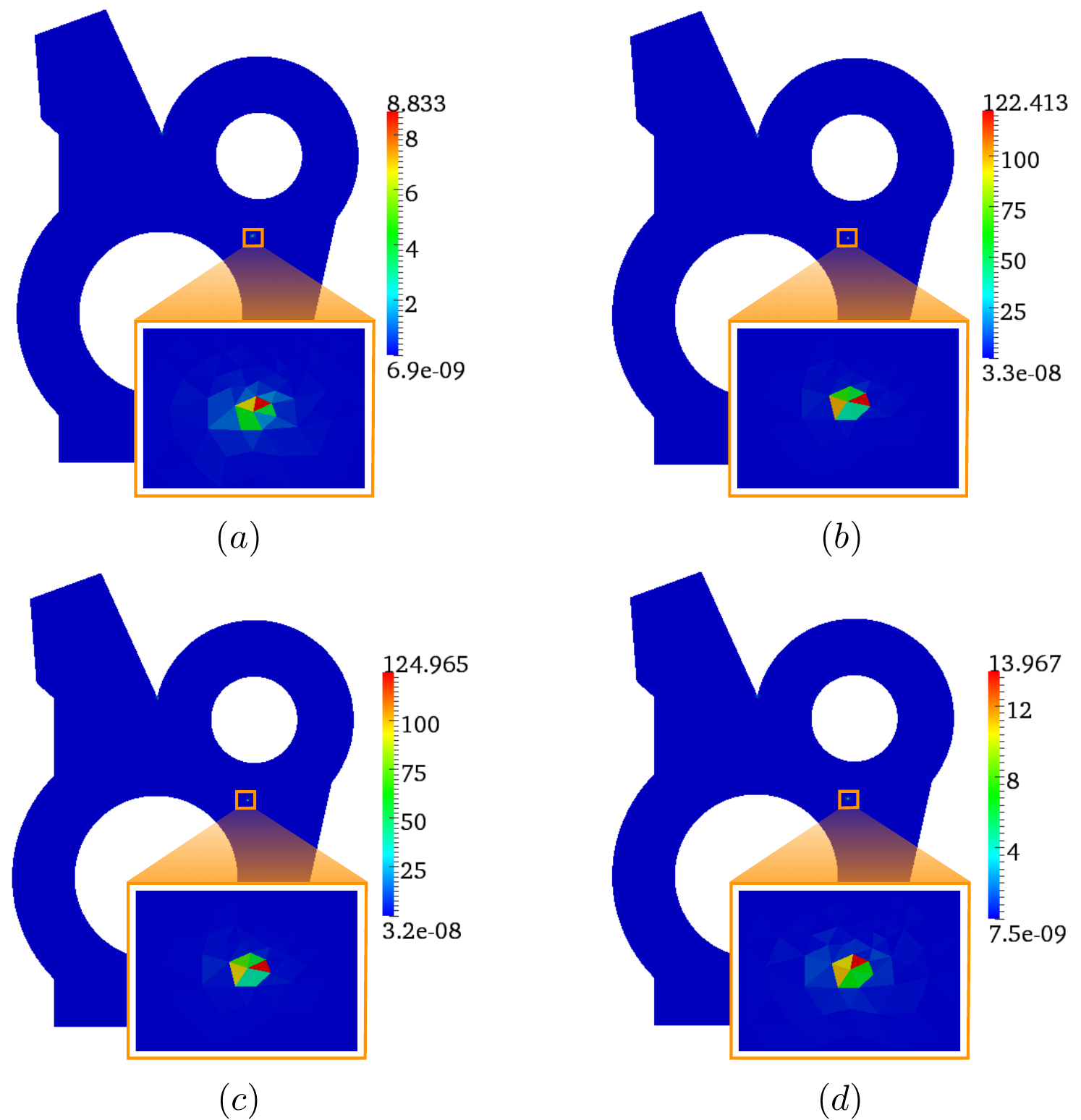}
\caption{Spatial distribution of local contributions to the density of the energy norm of the reference error (a) and that of the contributions to the density of the error estimates calculated using the standard versions of the EET (b) and the EESPT (c), and the full enhanced version (d).}
\label{fig1:distribution_discretization_error_estimators_structure_fissuree_2D}
\end{figure}

\Fig{fig1:sigma_EF_hat_structure_fissuree_2D} shows the magnitude $\displaystyle \sqrt{\Tr\big[\cont_h \: \cont_h \big]}$ of the FE stress field and that $\displaystyle \sqrt{\Tr\big[\hat{\cont}_h \: \hat{\cont}_h \big]}$ of the admissible stress fields obtained from the standard versions of EET and EESPT methods and the enhanced version in the case $\Ec_e = \Ec$ (hereafter referred to as the full enhanced version), \ie where all the elements of the FE mesh are implicated in the enhanced procedure. The highest stress region is located near the crack tip. As the FE mesh is refined adaptively in the vicinity of the crack tip, the spatial distribution of the reference and estimated error is expressed as a density which is the ratio between the squared elementary contribution to the reference (or estimated) error and the size of the element. \Fig{fig1:distribution_discretization_error_estimators_structure_fissuree_2D} displays the spatial distribution of the local contributions to the density of the energy norm of the reference error and that of the contributions to the density of the error estimates for the standard versions of EET and EESPT methods and the full enhanced version. These maps show that the higher contributions to the density of the reference and estimated error are concentrated in the elements surrounding the crack tip. Furthermore, the full enhanced version provides a better indicator of the elementary contributions to the energy norm of the reference error compared to the standard versions of EET and EESPT methods.

Given that the structure is subjected to homogeneous Dirichlet boundary conditions, the quasi-exact value of the energy norm of the discretization error (\ie the energy norm of the reference error) has been directly calculated:
\begin{equation}
\begin{aligned}
\lnorm{\und{e}_h}_{u, \Om} & = \sqrt{\lnorm{\uu}^2_{u, \Om} - \lnorm{\uu_h}^2_{u, \Om}} \\
& \simeq \sqrt{\lnorm{\uu_{ref}}^2_{u, \Om} - \lnorm{\uu_h}^2_{u, \Om}} \simeq 6.3302,
\end{aligned}
\end{equation}
and required a computational cost of about $1$ s, while that needed for calculating the local contributions to $\lnorm{\und{e}_h}_{u, \Om}$ reaches $3$ hours and $40$ minutes.

The graphs of the effectivity indices and corresponding normalized CPU time calculated using each criterion, namely the radius ratio, the edge ratio and the estimate ratio, are shown in \Figs{fig1:structure_fissuree_2D_radius_ratio}, \ref{fig1:structure_fissuree_2D_edge_ratio} and \ref{fig1:structure_fissuree_2D_estimate_ratio}, respectively, for both EET and EESPT methods. The r.h.s. graphs represent the evolution of the effectivity indices and corresponding normalized CPU time as functions of the number of elements involved in the enhanced procedure for both estimators. The normalized computational times obtained using the EET and EESPT methods have been computed with respect to the standard versions of each estimator. The radius ratio, edge ratio and estimate ratio range between $0.2923$ to $0.5000$, $0.4885$ to $0.9981$ and $9.5 \, 10^{-9}$ to $1.0$, respectively. Values of the effectivity indices and normalized CPU times corresponding to values $0.29$, $0.48$ and $1.1$ for the radius ratio, edge ratio and estimate ratio, respectively, are that calculated using the standard versions of EET and EESPT methods. First, the amount of computing time increases quasi-linearly with the number of elements involved in the minimization step. Then, one can observe that effectivity indices exhibit a downward behavior as the number of elements involved in the optimized procedure increases whatever the criterion we used. Nevertheless, the effectivity indices for both estimators drop more strongly in the case of the error estimate criterion than in the case of both geometric criteria. These reports are relevant as regards the choice of the error estimate criterion, as this one enables us to achieve sharper upper bounds while keeping an affordable computational cost. Besides, the use of other geometric criteria would probably lead to a trend similar to that obtained with the radius ratio or edge ratio.

\begin{figure*}
\centering
\begin{tikzpicture}[baseline]
\pgfplotsset{
xlabel near ticks,
ylabel near ticks,
tick label style={font=\footnotesize},
label style={font=\small},
legend style={font=\small},
try min ticks=7
}
\begin{axis}[
	width=0.36\textwidth,
	scale only axis,
	axis y line*=left,
	axis on top,
	xlabel=radius ratio,
	ylabel=effectivity index $\eta$,
	legend style={at={(0.04,0.48)},anchor=south west,legend columns=1},
	legend entries={$\eta_{EET}$,$\eta_{EESPT}$}
]
\addplot+[sharp plot,blue,solid,mark=*,mark options={blue}] table[x=rad_r,y=eta_EET] {structure_fissuree_2D_radius_ratio.txt};
\addplot+[sharp plot,red,solid,mark=triangle*,mark options={red}] table[x=rad_r,y=eta_EESPT] {structure_fissuree_2D_radius_ratio.txt};
\end{axis}
\begin{axis}[
	width=0.36\textwidth,
	scale only axis,
	ymax=240,
	axis x line=none,
	axis y line*=right,
	ylabel=normalized CPU time,
	legend style={at={(0.04,0.44)},anchor=north west,legend columns=1},
	legend entries={$t_{EET}$,$t_{EESPT}$}
]
\addplot+[sharp plot,teal,dashed,mark=*,mark options={teal}] table[x=rad_r,y expr=\thisrow{CPU_EET}/8] {structure_fissuree_2D_radius_ratio.txt};
\addplot+[sharp plot,orange,dashed,mark=triangle*,mark options={orange}] table[x=rad_r,y expr=\thisrow{CPU_EESPT}/9] {structure_fissuree_2D_radius_ratio.txt};
\end{axis}
\end{tikzpicture}
\begin{tikzpicture}[baseline]
\pgfplotsset{
xlabel near ticks,
ylabel near ticks,
tick label style={font=\footnotesize},
label style={font=\small},
legend style={font=\small},
try min ticks=7
}
\begin{axis}[
	width=0.36\textwidth,
	scaled x ticks=base 10:-3,
	scale only axis,
	axis y line*=left,
	axis on top,
	xlabel=number of elements,
	ylabel=effectivity index $\eta$,
	legend style={at={(0.96,0.96)},anchor=north east,legend columns=1},
]
\addplot+[sharp plot,blue,solid,mark=*,mark options={blue}] table[x=nb_elem,y=eta_EET] {structure_fissuree_2D_radius_ratio.txt};
\addplot+[sharp plot,red,solid,mark=triangle*,mark options={red}] table[x=nb_elem,y=eta_EESPT] {structure_fissuree_2D_radius_ratio.txt};
\end{axis}
\begin{axis}[
	width=0.36\textwidth,
	scale only axis,
	ymax=240,
	axis x line=none,
	axis y line*=right,
	ylabel=normalized CPU time,
	legend style={at={(0.04,0.48)},anchor=north west,legend columns=1},
]
\addplot+[sharp plot,teal,dashed,mark=*,mark options={teal}] table[x=nb_elem,y expr=\thisrow{CPU_EET}/8] {structure_fissuree_2D_radius_ratio.txt};
\addplot+[sharp plot,orange,dashed,mark=triangle*,mark options={orange}] table[x=nb_elem,y expr=\thisrow{CPU_EESPT}/9] {structure_fissuree_2D_radius_ratio.txt};
\end{axis}
\end{tikzpicture}
\caption{Effectivity indices and normalized CPU time for the error estimates EET and EESPT with respect to the radius ratio (left) and to the number of elements involved in the optimized procedure (right).}\label{fig1:structure_fissuree_2D_radius_ratio}
\end{figure*}
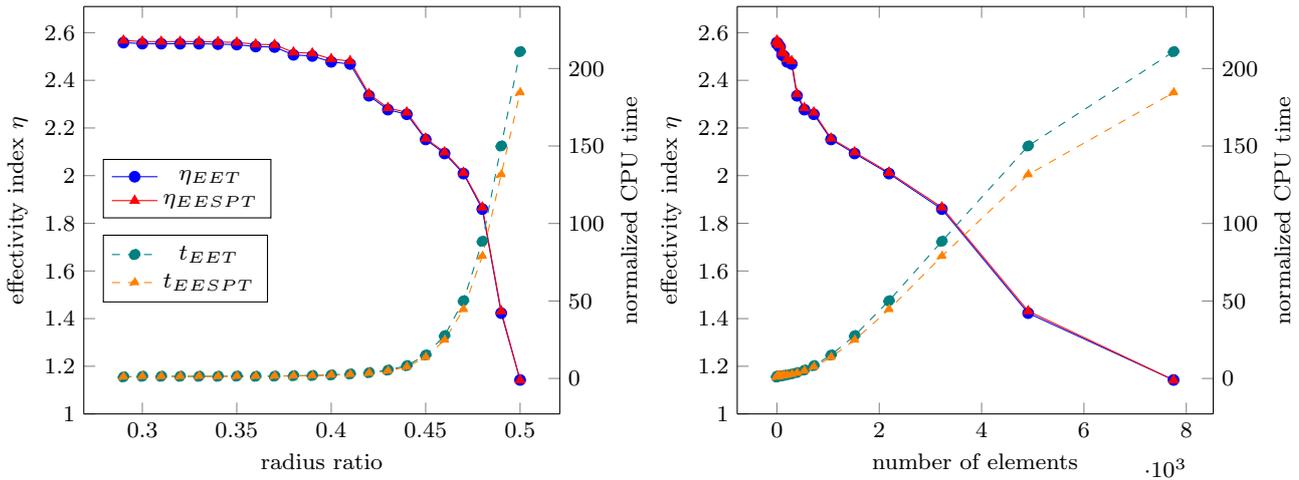

\begin{figure*}
\centering
\begin{tikzpicture}[baseline]
\pgfplotsset{
xlabel near ticks,
ylabel near ticks,
tick label style={font=\footnotesize},
label style={font=\small},
legend style={font=\small},
try min ticks=7
}
\begin{axis}[
	width=0.36\textwidth,
	scale only axis,
	axis y line*=left,
	axis on top,
	xlabel=edge ratio,
	ylabel=effectivity index $\eta$,
	legend style={at={(0.04,0.48)},anchor=south west,legend columns=1},
	legend entries={$\eta_{EET}$,$\eta_{EESPT}$}
]
\addplot+[sharp plot,blue,solid,mark=*,mark options={blue}] table[x=edg_r,y=eta_EET] {structure_fissuree_2D_edge_ratio.txt};
\addplot+[sharp plot,red,solid,mark=triangle*,mark options={red}] table[x=edg_r,y=eta_EESPT] {structure_fissuree_2D_edge_ratio.txt};
\end{axis}
\begin{axis}[
	width=0.36\textwidth,
	scale only axis,
	ymax=240,
	axis x line=none,
	axis y line*=right,
	ylabel=normalized CPU time,
	legend style={at={(0.04,0.44)},anchor=north west,legend columns=1},
	legend entries={$t_{EET}$,$t_{EESPT}$}
]
\addplot+[sharp plot,teal,dashed,mark=*,mark options={teal}] table[x=edg_r,y expr=\thisrow{CPU_EET}/8] {structure_fissuree_2D_edge_ratio.txt};
\addplot+[sharp plot,orange,dashed,mark=triangle*,mark options={orange}] table[x=edg_r,y expr=\thisrow{CPU_EESPT}/9] {structure_fissuree_2D_edge_ratio.txt};
\end{axis}
\end{tikzpicture}
\begin{tikzpicture}[baseline]
\pgfplotsset{
xlabel near ticks,
ylabel near ticks,
tick label style={font=\footnotesize},
label style={font=\small},
legend style={font=\small},
try min ticks=7
}
\begin{axis}[
	width=0.36\textwidth,
	scaled x ticks=base 10:-3,
	scale only axis,
	axis y line*=left,
	axis on top,
	xlabel=number of elements,
	ylabel=effectivity index $\eta$,
	legend style={at={(0.96,0.96)},anchor=north east,legend columns=1},
]
\addplot+[sharp plot,blue,solid,mark=*,mark options={blue}] table[x=nb_elem,y=eta_EET] {structure_fissuree_2D_edge_ratio.txt};
\addplot+[sharp plot,red,solid,mark=triangle*,mark options={red}] table[x=nb_elem,y=eta_EESPT] {structure_fissuree_2D_edge_ratio.txt};
\end{axis}
\begin{axis}[
	width=0.36\textwidth,
	scale only axis,
	ymax=240,
	axis x line=none,
	axis y line*=right,
	ylabel=normalized CPU time,
	legend style={at={(0.04,0.48)},anchor=north west,legend columns=1},
]
\addplot+[sharp plot,teal,dashed,mark=*,mark options={teal}] table[x=nb_elem,y expr=\thisrow{CPU_EET}/8] {structure_fissuree_2D_edge_ratio.txt};
\addplot+[sharp plot,orange,dashed,mark=triangle*,mark options={orange}] table[x=nb_elem,y expr=\thisrow{CPU_EESPT}/9] {structure_fissuree_2D_edge_ratio.txt};
\end{axis}
\end{tikzpicture}
\caption{Effectivity indices and normalized CPU time for the error estimates EET and EESPT with respect to the edge ratio (left) and to the number of elements involved in the optimized procedure (right).}\label{fig1:structure_fissuree_2D_edge_ratio}
\end{figure*}
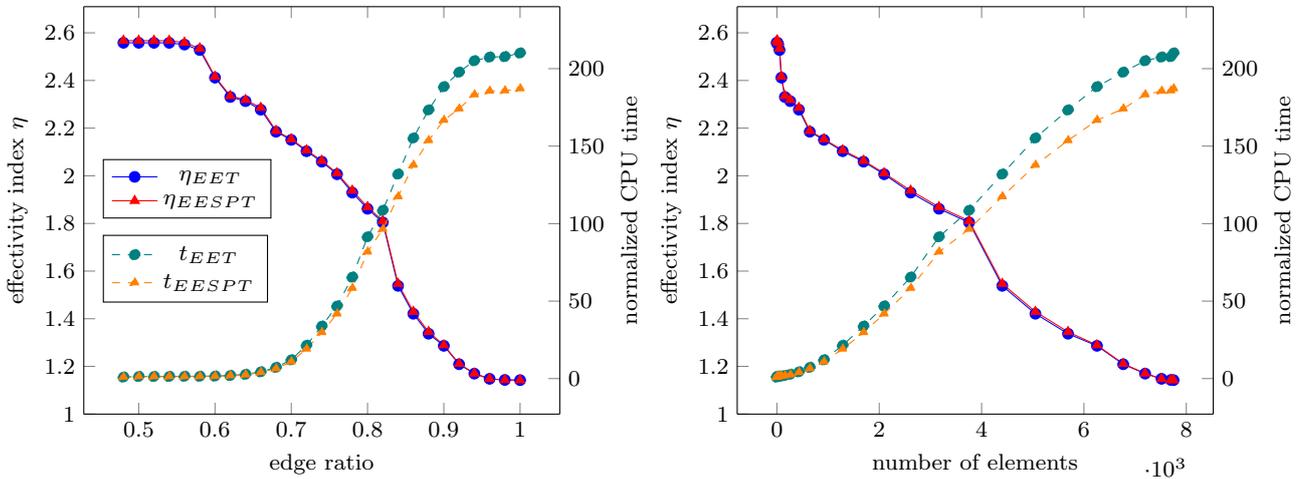

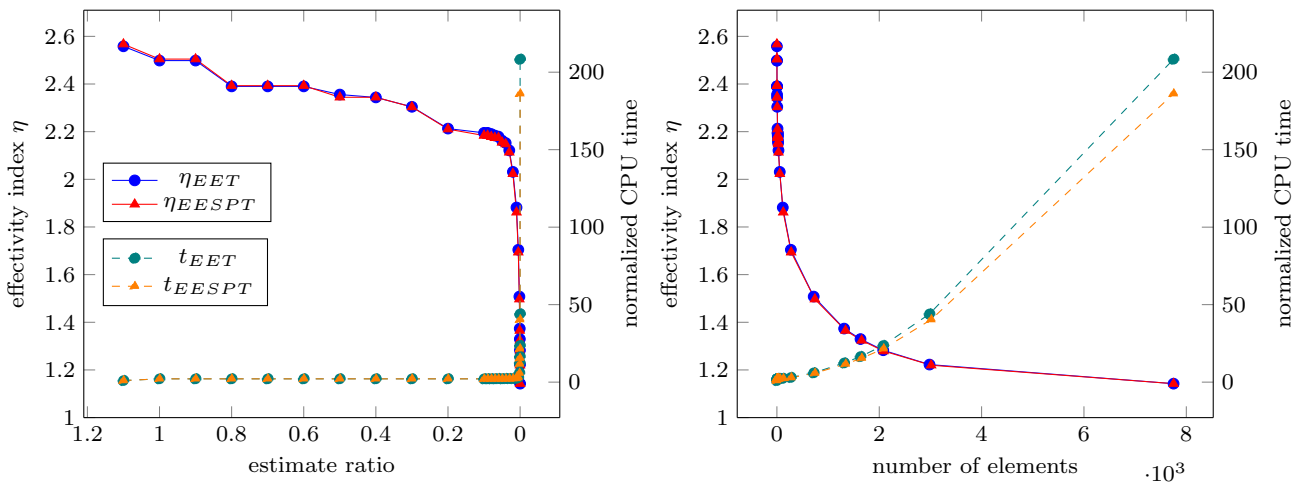
\begin{figure*}
\centering
\begin{tikzpicture}[baseline]
\pgfplotsset{
xlabel near ticks,
ylabel near ticks,
tick label style={font=\footnotesize},
label style={font=\small},
legend style={font=\small},
try min ticks=7
}
\begin{axis}[
	width=0.36\textwidth,
	scale only axis,
	axis y line*=left,
	axis on top,
	xlabel=estimate ratio,
	x dir=reverse,
	ylabel=effectivity index $\eta$,
	legend style={at={(0.04,0.48)},anchor=south west,legend columns=1},
	legend entries={$\eta_{EET}$,$\eta_{EESPT}$}
]
\addplot+[sharp plot,blue,solid,mark=*,mark options={blue}] table[x=est_r,y=eta_EET] {structure_fissuree_2D_estimate_ratio.txt};
\addplot+[sharp plot,red,solid,mark=triangle*,mark options={red}] table[x=est_r,y=eta_EESPT] {structure_fissuree_2D_estimate_ratio.txt};
\end{axis}
\begin{axis}[
	width=0.36\textwidth,
	scale only axis,
	ymax=240,
	axis x line=none,
	axis y line*=right,
	ylabel=normalized CPU time,
	x dir=reverse,
	legend style={at={(0.04,0.44)},anchor=north west,legend columns=1},
	legend entries={$t_{EET}$,$t_{EESPT}$}
]
\addplot+[sharp plot,teal,dashed,mark=*,mark options={teal}] table[x=est_r,y expr=\thisrow{CPU_EET}/8] {structure_fissuree_2D_estimate_ratio.txt};
\addplot+[sharp plot,orange,dashed,mark=triangle*,mark options={orange}] table[x=est_r,y expr=\thisrow{CPU_EESPT}/9] {structure_fissuree_2D_estimate_ratio.txt};
\end{axis}
\end{tikzpicture}
\begin{tikzpicture}[baseline]
\pgfplotsset{
xlabel near ticks,
ylabel near ticks,
tick label style={font=\footnotesize},
label style={font=\small},
legend style={font=\small},
try min ticks=7
}
\begin{axis}[
	width=0.36\textwidth,
	scaled x ticks=base 10:-3,
	scale only axis,
	axis y line*=left,
	axis on top,
	xlabel=number of elements,
	ylabel=effectivity index $\eta$,
	legend style={at={(0.96,0.96)},anchor=north east,legend columns=1},
]
\addplot+[sharp plot,blue,solid,mark=*,mark options={blue}] table[x=nb_elem_EET,y=eta_EET] {structure_fissuree_2D_estimate_ratio.txt};
\addplot+[sharp plot,red,solid,mark=triangle*,mark options={red}] table[x=nb_elem_EESPT,y=eta_EESPT] {structure_fissuree_2D_estimate_ratio.txt};
\end{axis}
\begin{axis}[
	width=0.36\textwidth,
	scale only axis,
	ymax=240,
	axis x line=none,
	axis y line*=right,
	ylabel=normalized CPU time,
	legend style={at={(0.04,0.48)},anchor=north west,legend columns=1},
]
\addplot+[sharp plot,teal,dashed,mark=*,mark options={teal}] table[x=nb_elem_EET,y expr=\thisrow{CPU_EET}/8] {structure_fissuree_2D_estimate_ratio.txt};
\addplot+[sharp plot,orange,dashed,mark=triangle*,mark options={orange}] table[x=nb_elem_EESPT,y expr=\thisrow{CPU_EESPT}/9] {structure_fissuree_2D_estimate_ratio.txt};
\end{axis}
\end{tikzpicture}
\caption{Effectivity indices and normalized CPU time for the error estimates EET and EESPT with respect to the estimate ratio (left) and to the number of elements involved in the optimized procedure (right).}\label{fig1:structure_fissuree_2D_estimate_ratio}
\end{figure*}

A sound efficiency factor is the ratio between the gain $g_{\eta}$ in effectivity index and the loss $l_t$ of computational cost with respect to the standard technique. The gain $g_{\eta}$ in effectivity index and the loss $l_t$ of computational cost are defined as 
\begin{equation}\label{eq1:efficiencyfactor}
g_{\eta} = \labs{\frac{\eta - \eta^{\std}}{\eta^{\std}}} \quad \text{and} \quad l_t = \labs{\frac{t - t^{\std}}{t^{\std}}}, \quad \text{respectively}, 
\end{equation}
where $\eta$ and $\eta^{\std}$ denote the effectivity indices obtained using the enhanced and standard procedures, respectively, while $t$ and $t^{\std}$ represent the corresponding normalized CPU times.
 
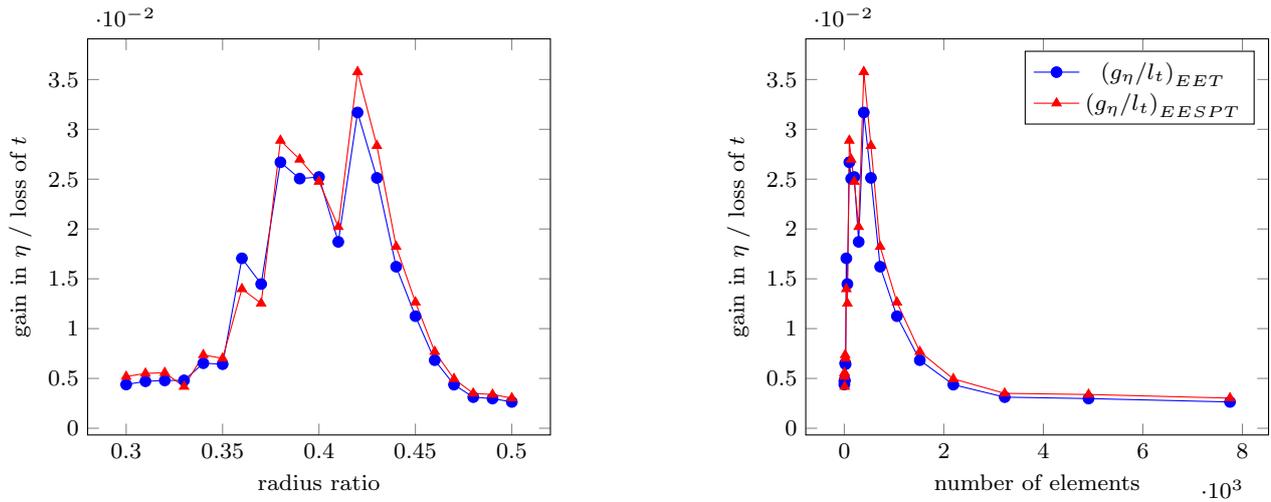
\begin{figure*}
\centering
\begin{tikzpicture}[baseline]
\pgfplotsset{
xlabel near ticks,
ylabel near ticks,
tick label style={font=\footnotesize},
label style={font=\small},
legend style={font=\small},
try min ticks=7
}
\begin{axis}[
	width=0.35\textwidth,
	scaled y ticks=base 10:2,
	scale only axis,
	xlabel=radius ratio,
	ylabel=gain in $\eta$ / loss of $t$,
	x filter/.code={
		\ifnum\coordindex<1
			\def\pgfmathresult{}
		\fi
	}
]
\addplot+[sharp plot,blue,solid,mark=*,mark options={blue}] table[x=rad_r,y expr=((2.558170-\thisrow{eta_EET})/2.558170)/((\thisrow{CPU_EET}-8)/8)] {structure_fissuree_2D_radius_ratio.txt};
\addplot+[sharp plot,red,solid,mark=triangle*,mark options={red}] table[x=rad_r,y expr=((2.567230-\thisrow{eta_EESPT})/2.567230)/((\thisrow{CPU_EESPT}-9)/9)] {structure_fissuree_2D_radius_ratio.txt};
\end{axis}
\end{tikzpicture}
\hspace{2cm}
\begin{tikzpicture}[baseline]
\pgfplotsset{
xlabel near ticks,
ylabel near ticks,
tick label style={font=\footnotesize},
label style={font=\small},
legend style={font=\small},
try min ticks=7
}
\begin{axis}[
	width=0.35\textwidth,
	scaled x ticks=base 10:-3,
	scaled y ticks=base 10:2,
	scale only axis,
	xlabel=number of elements,
	ylabel=gain in $\eta$ / loss of $t$,
	legend pos=north east,
	legend entries={$\displaystyle{(g_{\eta} / l_t )}_{EET}$, $\displaystyle{(g_{\eta} / l_t )}_{EESPT}$},
	x filter/.code={
		\ifnum\coordindex<1
			\def\pgfmathresult{}
		\fi
	}
]
\addplot+[sharp plot,blue,solid,mark=*,mark options={blue}] table[x=nb_elem,y expr=((2.558170-\thisrow{eta_EET})/2.558170)/((\thisrow{CPU_EET}-8)/8)] {structure_fissuree_2D_radius_ratio.txt};
\addplot+[sharp plot,red,solid,mark=triangle*,mark options={red}] table[x=nb_elem,y expr=((2.567230-\thisrow{eta_EESPT})/2.567230)/((\thisrow{CPU_EESPT}-9)/9)] {structure_fissuree_2D_radius_ratio.txt};
\end{axis}
\end{tikzpicture}
\caption{Ratio between gain in effectivity index and loss of normalized CPU time for the error estimates EET and EESPT with respect to the radius ratio (left) and to the number of elements involved in the optimized procedure (right).}\label{fig1:structure_fissuree_2D_radius_ratio_bis}
\end{figure*}

\begin{figure*}
\centering
\begin{tikzpicture}[baseline]
\pgfplotsset{
xlabel near ticks,
ylabel near ticks,
tick label style={font=\footnotesize},
label style={font=\small},
legend style={font=\small},
try min ticks=7
}
\begin{axis}[
	width=0.35\textwidth,
	scaled y ticks=base 10:2,
	scale only axis,
	xlabel=edge ratio,
	ylabel=gain in $\eta$ / loss of $t$,
	x filter/.code={
		\ifnum\coordindex<1
			\def\pgfmathresult{}
		\fi
	}
]
\addplot+[sharp plot,blue,solid,mark=*,mark options={blue}] table[x=edg_r,y expr=((2.558170-\thisrow{eta_EET})/2.558170)/((\thisrow{CPU_EET}-8)/8)] {structure_fissuree_2D_edge_ratio.txt};
\addplot+[sharp plot,red,solid,mark=triangle*,mark options={red}] table[x=edg_r,y expr=((2.567230-\thisrow{eta_EESPT})/2.567230)/((\thisrow{CPU_EESPT}-9)/9)] {structure_fissuree_2D_edge_ratio.txt};
\end{axis}
\end{tikzpicture}
\hspace{2cm}
\begin{tikzpicture}[baseline]
\pgfplotsset{
xlabel near ticks,
ylabel near ticks,
tick label style={font=\footnotesize},
label style={font=\small},
legend style={font=\small},
try min ticks=7
}
\begin{axis}[
	width=0.35\textwidth,
	scaled x ticks=base 10:-3,
	scaled y ticks=base 10:2,
	scale only axis,
	xlabel=number of elements,
	ylabel=gain in $\eta$ / loss of $t$,
	legend pos=north east,
	legend entries={$\displaystyle{(g_{\eta} / l_t )}_{EET}$, $\displaystyle{(g_{\eta} / l_t )}_{EESPT}$},
	x filter/.code={
		\ifnum\coordindex<1
			\def\pgfmathresult{}
		\fi
	}
]
\addplot+[sharp plot,blue,solid,mark=*,mark options={blue}] table[x=nb_elem,y expr=((2.558170-\thisrow{eta_EET})/2.558170)/((\thisrow{CPU_EET}-8)/8)] {structure_fissuree_2D_edge_ratio.txt};
\addplot+[sharp plot,red,solid,mark=triangle*,mark options={red}] table[x=nb_elem,y expr=((2.567230-\thisrow{eta_EESPT})/2.567230)/((\thisrow{CPU_EESPT}-9)/9)] {structure_fissuree_2D_edge_ratio.txt};
\end{axis}
\end{tikzpicture}
\caption{Ratio between gain in effectivity index and loss of normalized CPU time for the error estimates EET and EESPT with respect to the edge ratio (left) and to the number of elements involved in the optimized procedure (right).}\label{fig1:structure_fissuree_2D_edge_ratio_bis}
\end{figure*}
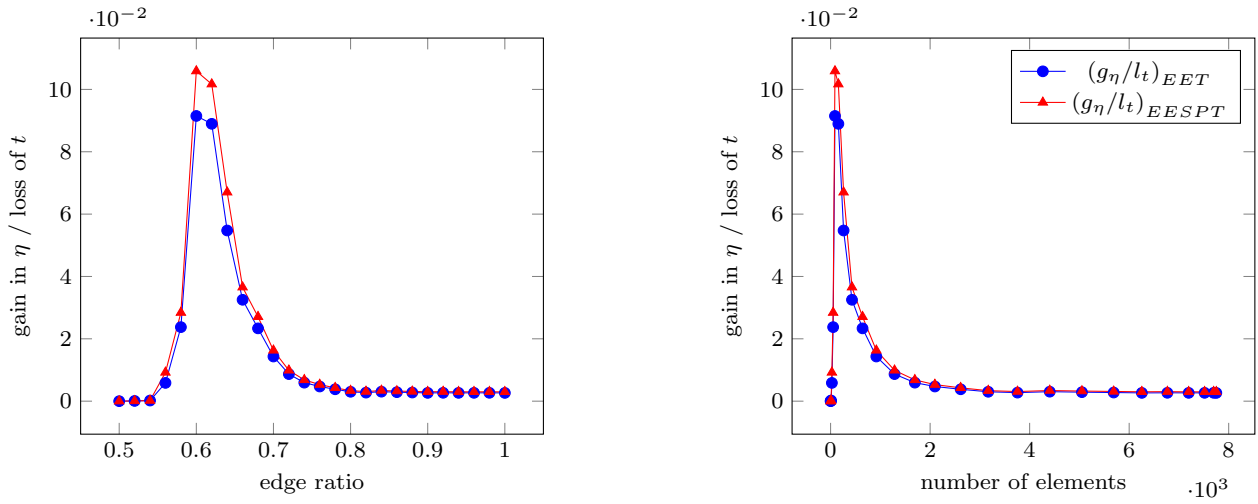

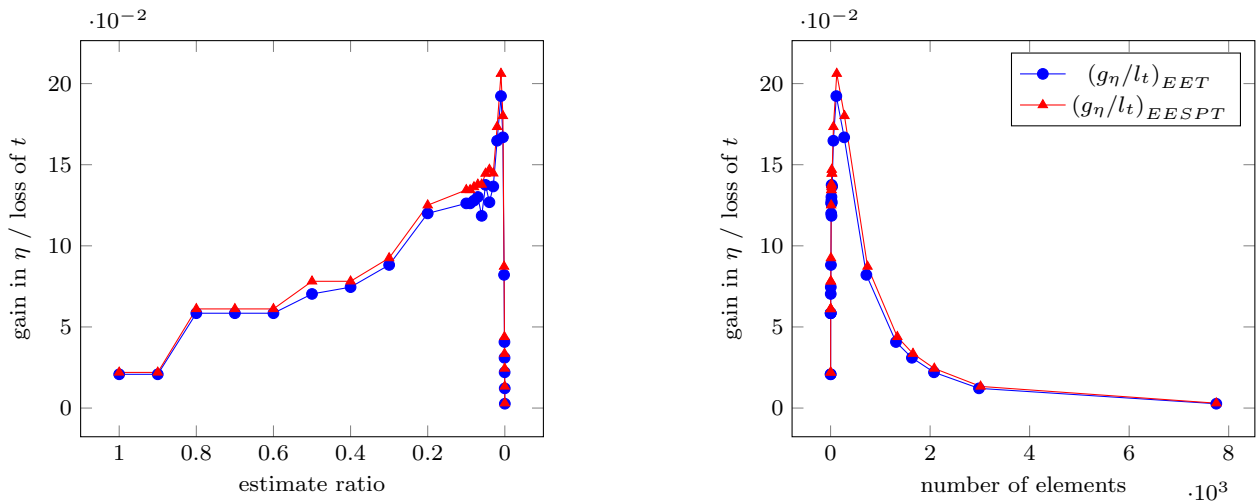
\begin{figure*}
\centering
\begin{tikzpicture}[baseline]
\pgfplotsset{
xlabel near ticks,
ylabel near ticks,
tick label style={font=\footnotesize},
label style={font=\small},
legend style={font=\small},
try min ticks=7
}
\begin{axis}[
	width=0.35\textwidth,
	scaled y ticks=base 10:2,
	scale only axis,
	xlabel=estimate ratio,
	x dir=reverse,
	ylabel=gain in $\eta$ / loss of $t$,
	x filter/.code={
		\ifnum\coordindex<1
			\def\pgfmathresult{}
		\fi
	}
]
\addplot+[sharp plot,blue,solid,mark=*,mark options={blue}] table[x=est_r,y expr=((2.558170-\thisrow{eta_EET})/2.558170)/((\thisrow{CPU_EET}-8)/8)] {structure_fissuree_2D_estimate_ratio.txt};
\addplot+[sharp plot,red,solid,mark=triangle*,mark options={red}] table[x=est_r,y expr=((2.567230-\thisrow{eta_EESPT})/2.567230)/((\thisrow{CPU_EESPT}-9)/9)] {structure_fissuree_2D_estimate_ratio.txt};
\end{axis}
\end{tikzpicture}
\hspace{2cm}
\begin{tikzpicture}[baseline]
\pgfplotsset{
xlabel near ticks,
ylabel near ticks,
tick label style={font=\footnotesize},
label style={font=\small},
legend style={font=\small},
try min ticks=7
}
\begin{axis}[
	width=0.35\textwidth,
	scaled x ticks=base 10:-3,
	scaled y ticks=base 10:2,
	scale only axis,
	xlabel=number of elements,
	ylabel=gain in $\eta$ / loss of $t$,
	legend pos=north east,
	legend entries={$\displaystyle{(g_{\eta} / l_t )}_{EET}$, $\displaystyle{(g_{\eta} / l_t )}_{EESPT}$},
	x filter/.code={
		\ifnum\coordindex<1
			\def\pgfmathresult{}
		\fi
	}
]
\addplot+[sharp plot,blue,solid,mark=*,mark options={blue}] table[x=nb_elem_EET,y expr=((2.558170-\thisrow{eta_EET})/2.558170)/((\thisrow{CPU_EET}-8)/8)] {structure_fissuree_2D_estimate_ratio.txt};
\addplot+[sharp plot,red,solid,mark=triangle*,mark options={red}] table[x=nb_elem_EESPT,y expr=((2.567230-\thisrow{eta_EESPT})/2.567230)/((\thisrow{CPU_EESPT}-9)/9)] {structure_fissuree_2D_estimate_ratio.txt};
\end{axis}
\end{tikzpicture}
\caption{Ratio between gain in effectivity index and loss of normalized CPU time for the error estimates EET and EESPT with respect to the estimate ratio (left) and to the number of elements involved in the optimized procedure (right).}\label{fig1:structure_fissuree_2D_estimate_ratio_bis}
\end{figure*}

\Figs{fig1:structure_fissuree_2D_radius_ratio_bis}, \ref{fig1:structure_fissuree_2D_edge_ratio_bis} and \ref{fig1:structure_fissuree_2D_estimate_ratio_bis} represent the evolution of this efficiency factor for both estimators as a function of the radius ratio, the edge ratio and the estimate ratio, respectively, (l.h.s. graphs) and as function of the number of elements involved in the enhanced procedure (r.h.s. graphs). One can observe that the efficiency factor $g_{\eta} / l_t$ becomes meaningful when a small percentage of the elements are implicated in the enhanced procedure regardless of the criterion. Nevertheless this efficiency factor obtained using the error estimate criterion reaches a value at least twice as high as the one given by the use of geometric criteria.

\subsection{Weight sensor}\label{6.2}

\begin{figure*}
\centering\includegraphics[width=0.9\textwidth]{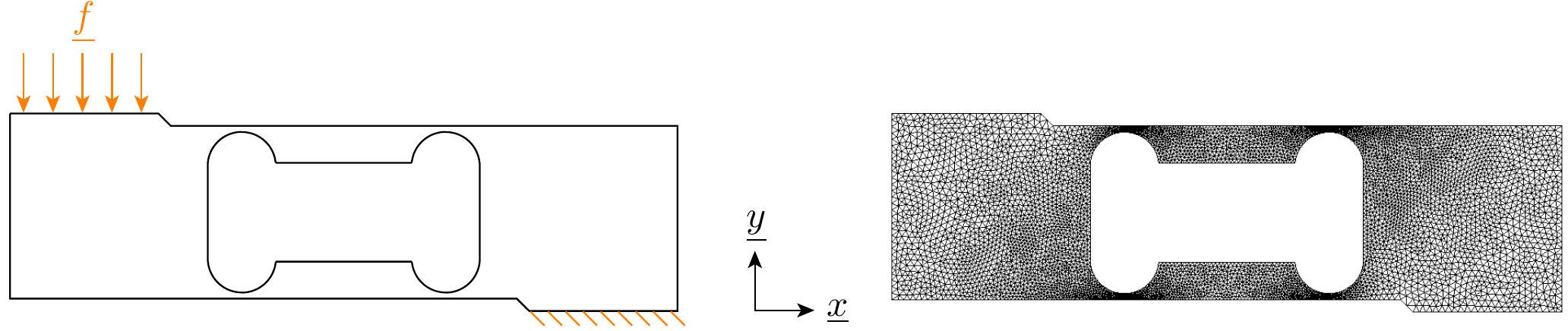}
\caption{Weight sensor model problem (left) and associated finite element mesh (right).}
\label{fig1:capteur_effort_2D_geometry_fine_mesh}
\end{figure*}

Now, let us consider a two-dimensional model of a weight sensor with two symmetric holes represented in \Fig{fig1:capteur_effort_2D_geometry_fine_mesh}. The structure, already studied in \cite{Ple11}, is loaded with a unit force density $\und{f} = -\und{y}$ along the top-left horizontal edge. Displacements are set to zero along the bottom-right horizontal edge. All the remaining edges are traction-free boundaries. The geometry and mesh considered, made of $11 \, 766$ linear triangular elements and $6 \, 299$ nodes (\ie $12 \, 598$ d.o.f.), are given in \Fig{fig1:capteur_effort_2D_geometry_fine_mesh}. The mesh is uniformly densified around the top and bottom regions of the two round cavities, which are the highest stress zones. The reference solution is performed using a very fine mesh consisting of $3 \, 331 \, 632$ linear triangular elements and $1 \, 671 \, 043$ nodes (\ie $3 \, 334 \, 086$ d.o.f.). 

\begin{figure}
\centering\includegraphics[scale = 0.32]{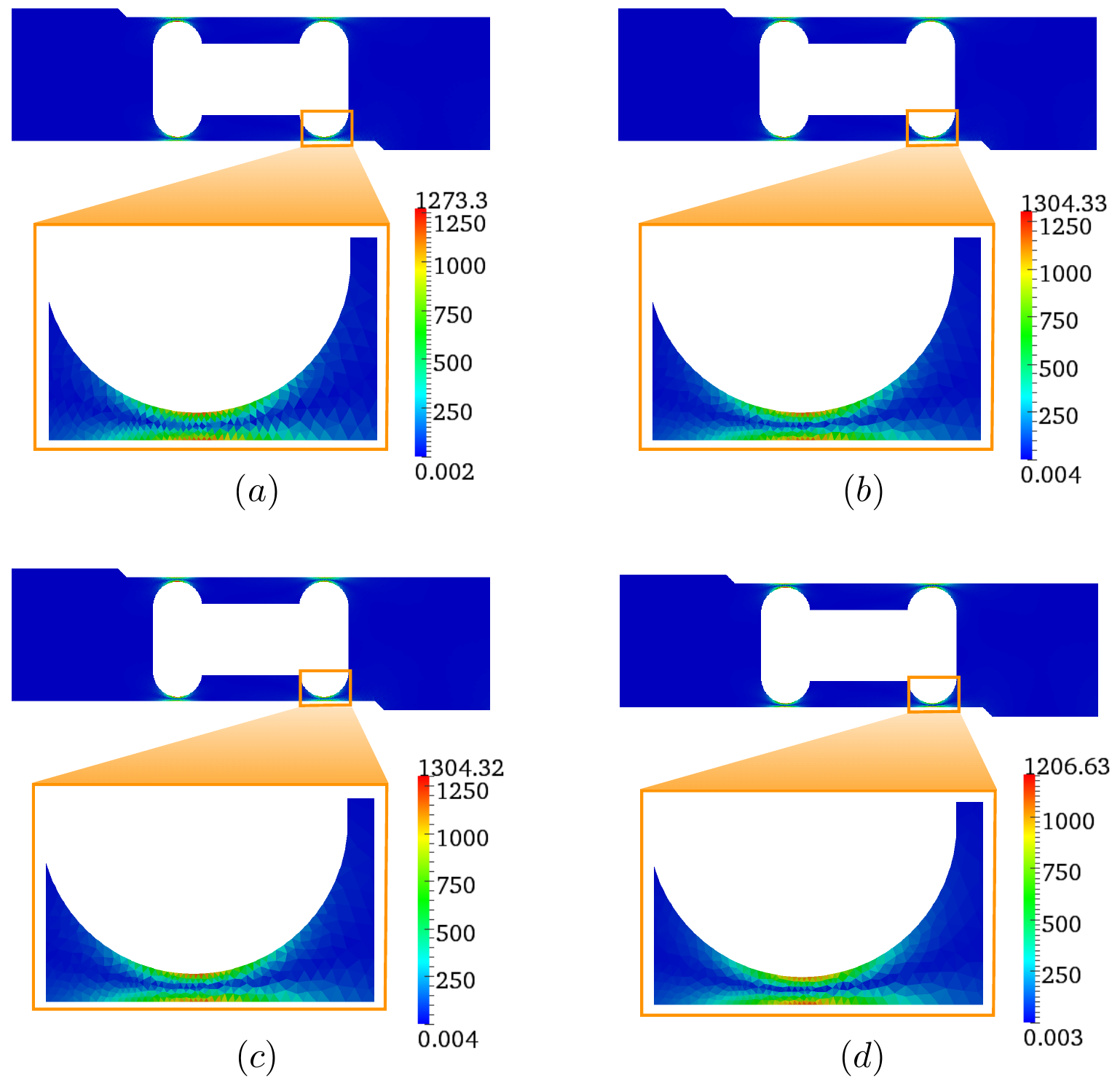}
\caption{Magnitude of the FE stress field (a) and the admissible stress field calculated using the standard versions of EET (b) and EESPT (c), the full enhanced version (d).}
\label{fig1:sigma_EF_hat_capteur_effort_2D}
\end{figure}

\begin{figure}
\centering\includegraphics[scale = 0.32]{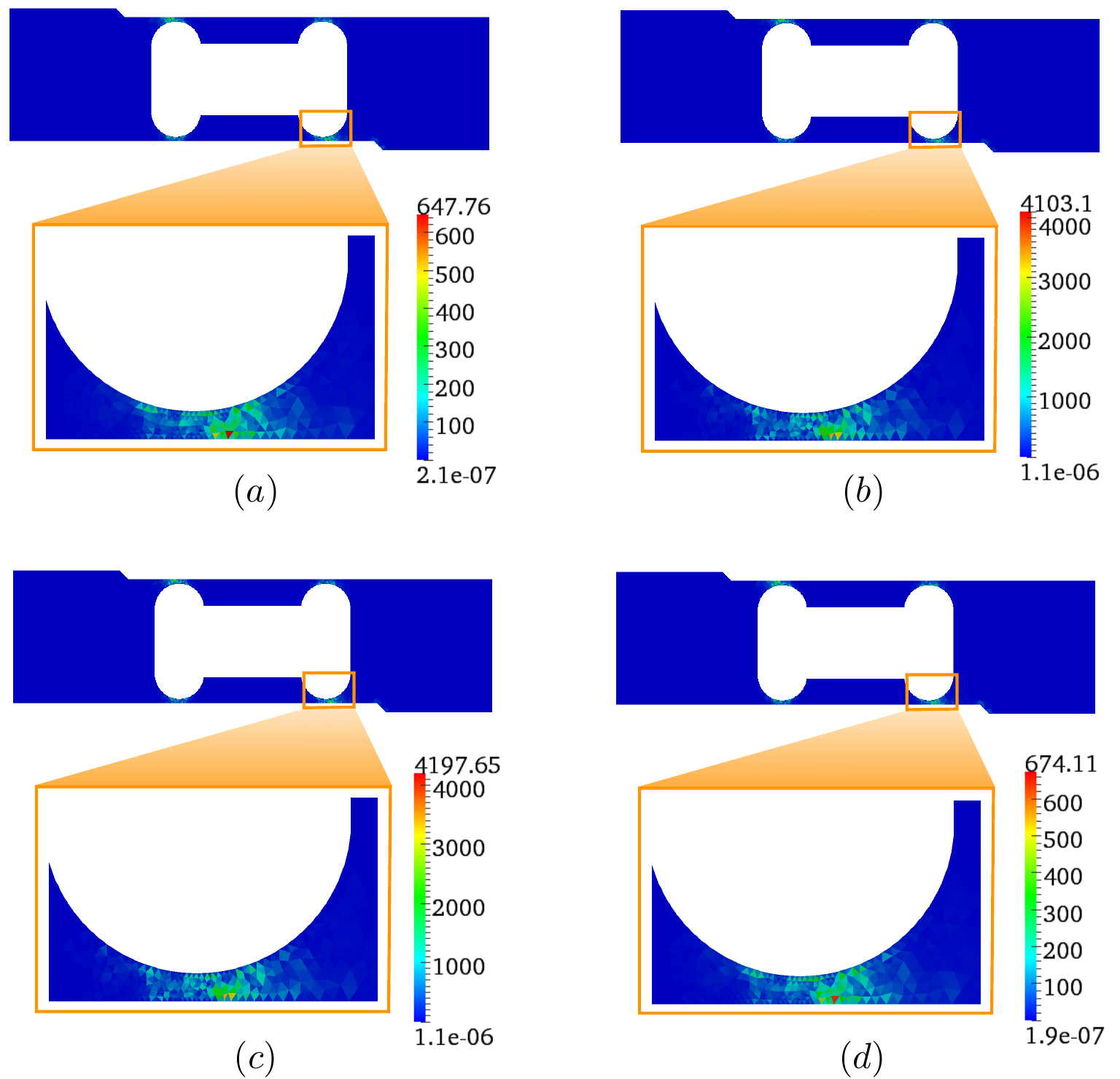}
\caption{Spatial distribution of local contributions to the energy norm of the reference error (a) and that of the contributions to the error estimates calculated using the standard versions of the EET (b) and the EESPT (c), and the full enhanced version (d).}
\label{fig1:distribution_discretization_error_estimators_capteur_effort_2D}
\end{figure}

The FE stress field and the admissible stress fields obtained from the standard versions of EET and EESPT methods and the full enhanced version are displayed in \Fig{fig1:sigma_EF_hat_capteur_effort_2D}. Maps of the local contributions to the energy norm of the reference error and that of contributions to the error estimates obtained from the standard versions of EET and EESPT methods are shown in \Fig{fig1:distribution_discretization_error_estimators_capteur_effort_2D}. One can observe that the error is by a majority localized around the top and bottom regions of the two cavities. Besides, maps $(a)$ and $(d)$ are strikingly close. More precisely, values of the elementary contributions to the estimated error computed using the full enhanced version do not overestimate that of the contributions to $\lnorm{\und{e}_h}_{u, \Om}$, contrary to the standard versions of EET and EESPT methods.

The calculation of the value of the energy norm of the reference error leads to:
\begin{equation}
\begin{aligned}
\lnorm{\und{e}_h}_{u, \Om} & = \sqrt{\lnorm{\uu}^2_{u, \Om} - \lnorm{\uu_h}^2_{u, \Om}} \\
& \simeq \sqrt{\lnorm{\uu_{ref}}^2_{u, \Om} - \lnorm{\uu_h}^2_{u, \Om}} \simeq 335.005,
\end{aligned}
\end{equation}
and the CPU time is about $2$ s, whereas that needed for the calculation of the positive contributions of each element of the mesh providing local indicator of the local energy norm of the reference error reaches about $15$ hours.

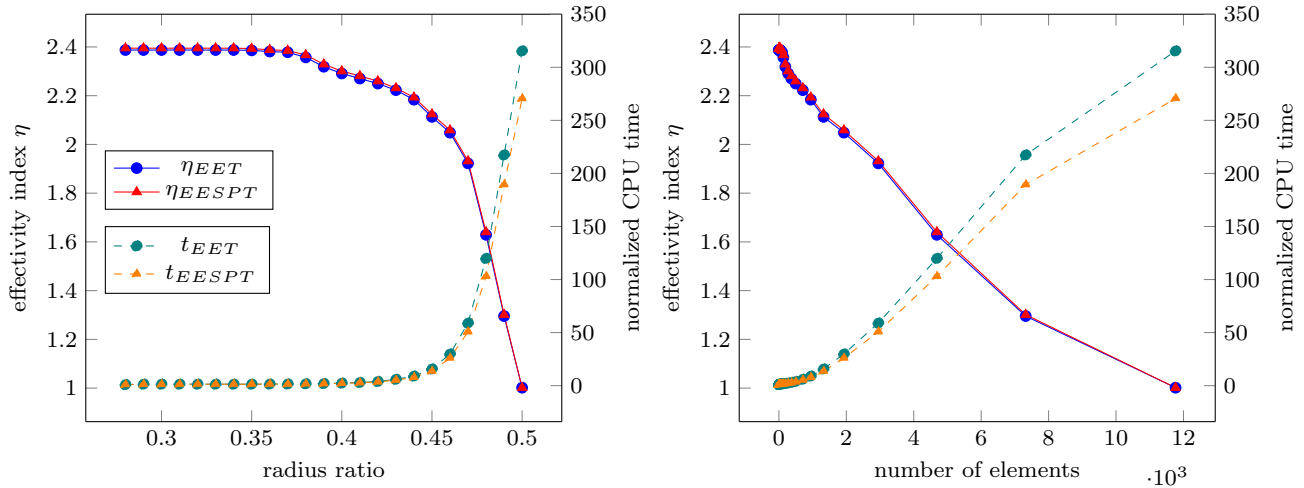
\begin{figure*}
\centering
\begin{tikzpicture}[baseline]
\pgfplotsset{
xlabel near ticks,
ylabel near ticks,
tick label style={font=\footnotesize},
label style={font=\small},
legend style={font=\small},
try min ticks=7
}
\begin{axis}[
	width=0.36\textwidth,
	scale only axis,
	axis y line*=left,
	axis on top,
	xlabel=radius ratio,
	ylabel=effectivity index $\eta$,
	legend style={at={(0.04,0.52)},anchor=south west,legend columns=1},
	legend entries={$\eta_{EET}$,$\eta_{EESPT}$}
]
\addplot+[sharp plot,blue,solid,mark=*,mark options={blue}] table[x=rad_r,y=eta_EET] {capteur_effort_2D_radius_ratio.txt};
\addplot+[sharp plot,red,solid,mark=triangle*,mark options={red}] table[x=rad_r,y=eta_EESPT] {capteur_effort_2D_radius_ratio.txt};
\end{axis}
\begin{axis}[
	width=0.36\textwidth,
	scale only axis,
	ymax=350,
	axis x line=none,
	axis y line*=right,
	ylabel=normalized CPU time,
	legend style={at={(0.04,0.48)},anchor=north west,legend columns=1},
	legend entries={$t_{EET}$,$t_{EESPT}$}
]
\addplot+[sharp plot,teal,dashed,mark=*,mark options={teal}] table[x=rad_r,y expr=\thisrow{CPU_EET}/12] {capteur_effort_2D_radius_ratio.txt};
\addplot+[sharp plot,orange,dashed,mark=triangle*,mark options={orange}] table[x=rad_r,y expr=\thisrow{CPU_EESPT}/14] {capteur_effort_2D_radius_ratio.txt};
\end{axis}
\end{tikzpicture}
\begin{tikzpicture}[baseline]
\pgfplotsset{
xlabel near ticks,
ylabel near ticks,
tick label style={font=\footnotesize},
label style={font=\small},
legend style={font=\small},
try min ticks=7
}
\begin{axis}[
	width=0.36\textwidth,
	scaled x ticks=base 10:-3,
	scale only axis,
	axis y line*=left,
	axis on top,
	xlabel=number of elements,
	ylabel=effectivity index $\eta$,
	legend style={at={(0.96,0.96)},anchor=north east,legend columns=1},
]
\addplot+[sharp plot,blue,solid,mark=*,mark options={blue}] table[x=nb_elem,y=eta_EET] {capteur_effort_2D_radius_ratio.txt};
\addplot+[sharp plot,red,solid,mark=triangle*,mark options={red}] table[x=nb_elem,y=eta_EESPT] {capteur_effort_2D_radius_ratio.txt};
\end{axis}
\begin{axis}[
	width=0.36\textwidth,
	scale only axis,
	ymax=350,
	axis x line=none,
	axis y line*=right,
	ylabel=normalized CPU time,
	legend style={at={(0.04,0.48)},anchor=north west,legend columns=1},
]
\addplot+[sharp plot,teal,dashed,mark=*,mark options={teal}] table[x=nb_elem,y expr=\thisrow{CPU_EET}/12] {capteur_effort_2D_radius_ratio.txt};
\addplot+[sharp plot,orange,dashed,mark=triangle*,mark options={orange}] table[x=nb_elem,y expr=\thisrow{CPU_EESPT}/14] {capteur_effort_2D_radius_ratio.txt};
\end{axis}
\end{tikzpicture}
\caption{Effectivity indices and normalized CPU time for the error estimates EET and EESPT with respect to the radius ratio (left) and to the number of elements involved in the optimized procedure (right).}\label{fig1:capteur_effort_2D_radius_ratio}
\end{figure*}

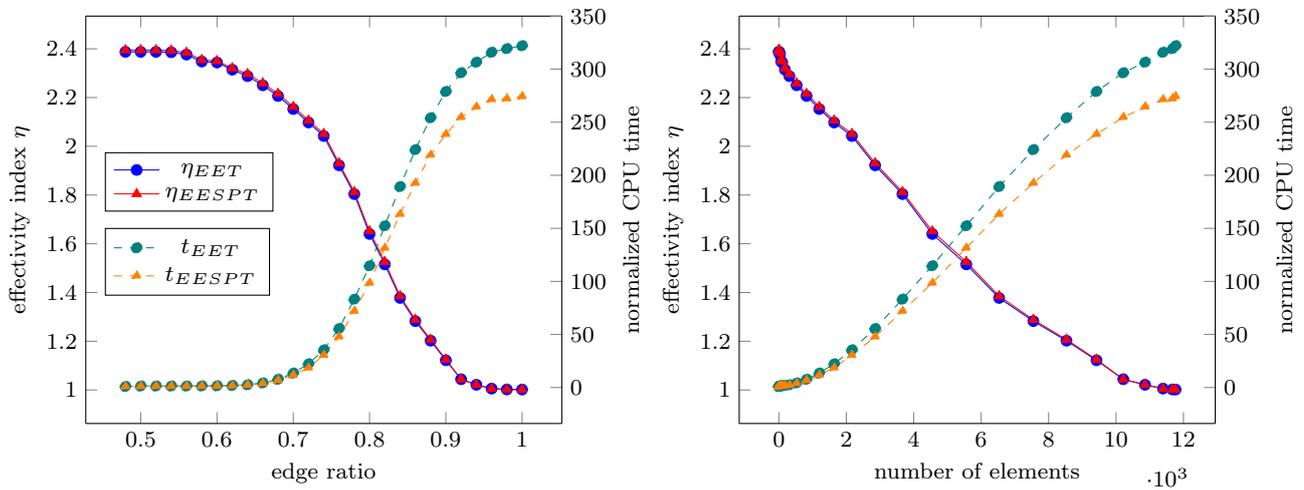
\begin{figure*}
\centering
\begin{tikzpicture}[baseline]
\pgfplotsset{
xlabel near ticks,
ylabel near ticks,
tick label style={font=\footnotesize},
label style={font=\small},
legend style={font=\small},
try min ticks=7
}
\begin{axis}[
	width=0.36\textwidth,
	scale only axis,
	axis y line*=left,
	axis on top,
	xlabel=edge ratio,
	ylabel=effectivity index $\eta$,
	legend style={at={(0.04,0.52)},anchor=south west,legend columns=1},
	legend entries={$\eta_{EET}$,$\eta_{EESPT}$}
]
\addplot+[sharp plot,blue,solid,mark=*,mark options={blue}] table[x=edg_r,y=eta_EET] {capteur_effort_2D_edge_ratio.txt};
\addplot+[sharp plot,red,solid,mark=triangle*,mark options={red}] table[x=edg_r,y=eta_EESPT] {capteur_effort_2D_edge_ratio.txt};
\end{axis}
\begin{axis}[
	width=0.36\textwidth,
	scale only axis,
	ymax=350,
	axis x line=none,
	axis y line*=right,
	ylabel=normalized CPU time,
	legend style={at={(0.04,0.48)},anchor=north west,legend columns=1},
	legend entries={$t_{EET}$,$t_{EESPT}$}
]
\addplot+[sharp plot,teal,dashed,mark=*,mark options={teal}] table[x=edg_r,y expr=\thisrow{CPU_EET}/12] {capteur_effort_2D_edge_ratio.txt};
\addplot+[sharp plot,orange,dashed,mark=triangle*,mark options={orange}] table[x=edg_r,y expr=\thisrow{CPU_EESPT}/14] {capteur_effort_2D_edge_ratio.txt};
\end{axis}
\end{tikzpicture}
\begin{tikzpicture}[baseline]
\pgfplotsset{
xlabel near ticks,
ylabel near ticks,
tick label style={font=\footnotesize},
label style={font=\small},
legend style={font=\small},
try min ticks=7
}
\begin{axis}[
	width=0.36\textwidth,
	scaled x ticks=base 10:-3,
	scale only axis,
	axis y line*=left,
	axis on top,
	xlabel=number of elements,
	ylabel=effectivity index $\eta$,
	legend style={at={(0.96,0.96)},anchor=north east,legend columns=1},
]
\addplot+[sharp plot,blue,solid,mark=*,mark options={blue}] table[x=nb_elem,y=eta_EET] {capteur_effort_2D_edge_ratio.txt};
\addplot+[sharp plot,red,solid,mark=triangle*,mark options={red}] table[x=nb_elem,y=eta_EESPT] {capteur_effort_2D_edge_ratio.txt};
\end{axis}
\begin{axis}[
	width=0.36\textwidth,
	scale only axis,
	ymax=350,
	axis x line=none,
	axis y line*=right,
	ylabel=normalized CPU time,
	legend style={at={(0.04,0.48)},anchor=north west,legend columns=1},
]
\addplot+[sharp plot,teal,dashed,mark=*,mark options={teal}] table[x=nb_elem,y expr=\thisrow{CPU_EET}/12] {capteur_effort_2D_edge_ratio.txt};
\addplot+[sharp plot,orange,dashed,mark=triangle*,mark options={orange}] table[x=nb_elem,y expr=\thisrow{CPU_EESPT}/14] {capteur_effort_2D_edge_ratio.txt};
\end{axis}
\end{tikzpicture}
\caption{Effectivity indices and normalized CPU time for the error estimates EET and EESPT with respect to the edge ratio (left) and to the number of elements involved in the optimized procedure (right).}\label{fig1:capteur_effort_2D_edge_ratio}
\end{figure*}

\begin{figure*}
\centering
\begin{tikzpicture}[baseline]
\pgfplotsset{
xlabel near ticks,
ylabel near ticks,
tick label style={font=\footnotesize},
label style={font=\small},
legend style={font=\small},
try min ticks=7
}
\begin{axis}[
	width=0.36\textwidth,
	scale only axis,
	axis y line*=left,
	axis on top,
	xlabel=estimate ratio,
	x dir=reverse,
	ylabel=effectivity index $\eta$,
	legend style={at={(0.04,0.52)},anchor=south west,legend columns=1},
	legend entries={$\eta_{EET}$,$\eta_{EESPT}$}
]
\addplot+[sharp plot,blue,solid,mark=*,mark options={blue}] table[x=est_r,y=eta_EET] {capteur_effort_2D_estimate_ratio.txt};
\addplot+[sharp plot,red,solid,mark=triangle*,mark options={red}] table[x=est_r,y=eta_EESPT] {capteur_effort_2D_estimate_ratio.txt};
\end{axis}
\begin{axis}[
	width=0.36\textwidth,
	scale only axis,
	ymax=350,
	axis x line=none,
	axis y line*=right,
	ylabel=normalized CPU time,
	x dir=reverse,
	legend style={at={(0.04,0.48)},anchor=north west,legend columns=1},
	legend entries={$t_{EET}$,$t_{EESPT}$}
]
\addplot+[sharp plot,teal,dashed,mark=*,mark options={teal}] table[x=est_r,y expr=\thisrow{CPU_EET}/12] {capteur_effort_2D_estimate_ratio.txt};
\addplot+[sharp plot,orange,dashed,mark=triangle*,mark options={orange}] table[x=est_r,y expr=\thisrow{CPU_EESPT}/14] {capteur_effort_2D_estimate_ratio.txt};
\end{axis}
\end{tikzpicture}
\begin{tikzpicture}[baseline]
\pgfplotsset{
xlabel near ticks,
ylabel near ticks,
tick label style={font=\footnotesize},
label style={font=\small},
legend style={font=\small},
try min ticks=7
}
\begin{axis}[
	width=0.36\textwidth,
	scaled x ticks=base 10:-3,
	scale only axis,
	axis y line*=left,
	axis on top,
	xlabel=number of elements,
	ylabel=effectivity index $\eta$,
	legend style={at={(0.96,0.96)},anchor=north east,legend columns=1},
]
\addplot+[sharp plot,blue,solid,mark=*,mark options={blue}] table[x=nb_elem_EET,y=eta_EET] {capteur_effort_2D_estimate_ratio.txt};
\addplot+[sharp plot,red,solid,mark=triangle*,mark options={red}] table[x=nb_elem_EESPT,y=eta_EESPT] {capteur_effort_2D_estimate_ratio.txt};
\end{axis}
\begin{axis}[
	width=0.36\textwidth,
	scale only axis,
	ymax=350,
	axis x line=none,
	axis y line*=right,
	ylabel=normalized CPU time,
	legend style={at={(0.04,0.48)},anchor=north west,legend columns=1},
]
\addplot+[sharp plot,teal,dashed,mark=*,mark options={teal}] table[x=nb_elem_EET,y expr=\thisrow{CPU_EET}/12] {capteur_effort_2D_estimate_ratio.txt};
\addplot+[sharp plot,orange,dashed,mark=triangle*,mark options={orange}] table[x=nb_elem_EESPT,y expr=\thisrow{CPU_EESPT}/14] {capteur_effort_2D_estimate_ratio.txt};
\end{axis}
\end{tikzpicture}
\caption{Effectivity indices and normalized CPU time for the error estimates EET and EESPT with respect to the estimate ratio (left) and to the number of elements involved in the optimized procedure (right).}\label{fig1:capteur_effort_2D_estimate_ratio}
\end{figure*}
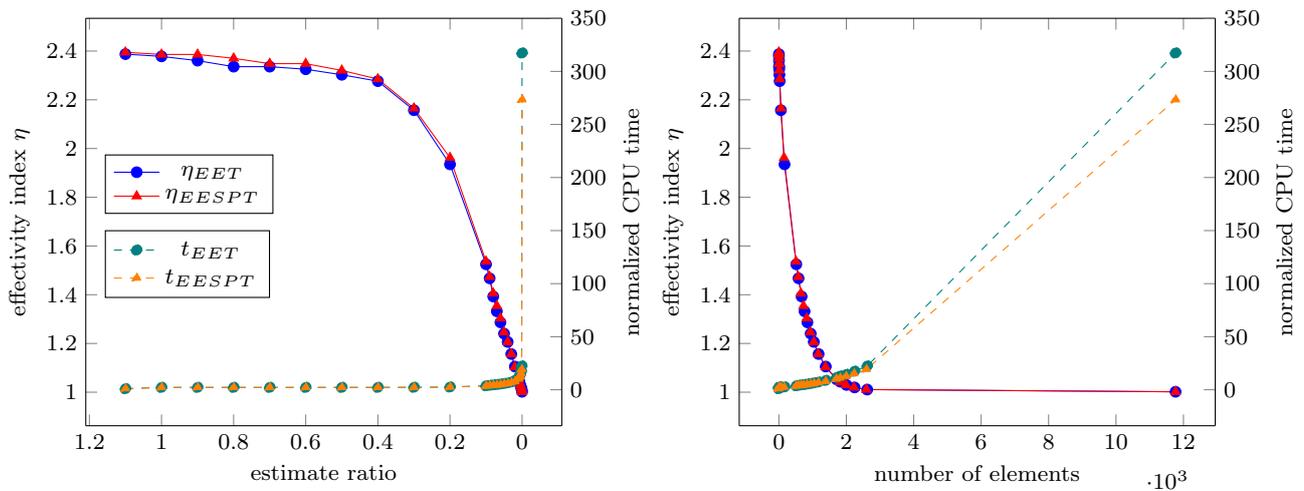

The effectivity indices and the corresponding normalized times of simulation computed with respect to the EET and EESPT methods are plotted versus the radius ratio, the edge ratio and the estimate ratio, respectively, in \Figs{fig1:capteur_effort_2D_radius_ratio}, \ref{fig1:capteur_effort_2D_edge_ratio} and \ref{fig1:capteur_effort_2D_estimate_ratio}. The radius ratio, edge ratio and estimate ratio range between $0.2863$ to $0.5000$, $0.4909$ to $0.9989$ and $3.10^{-10}$ to $1.0$, respectively. Values of the effectivity indices and normalized CPU times corresponding to values $0.28$, $0.48$ and $1.1$ for the radius ratio, edge ratio and estimate ratio, respectively, are that calculated using the standard versions of EET and EESPT methods. Furthermore, those figures show the evolutions of the effectivity indices and corresponding CPU times as functions of the number of elements involved in the enhanced procedure for both estimators. For both EET and EESPT estimators, the upper bounds of $\lnorm{\und{e}_h}_{u, \Om}$ become gradually more precise as the number of elements involved in the minimization procedure increases for all the criteria we considered. It confirms that quality of the balanced tractions has a strong influence on resulting error bounds. More precisely, the analysis of the curves reveals that the effectivity indices for both estimators plummet in the case of the error estimate criterion and reaches a value very close to one with only a small part of the elements involved in the optimization. On the contrary, in the cases of both geometric criteria (radius ratio and edge ratio), the evolution of the global effectivity indices remains more or less linear with respect to the number of elements involved in the enhanced procedure. 

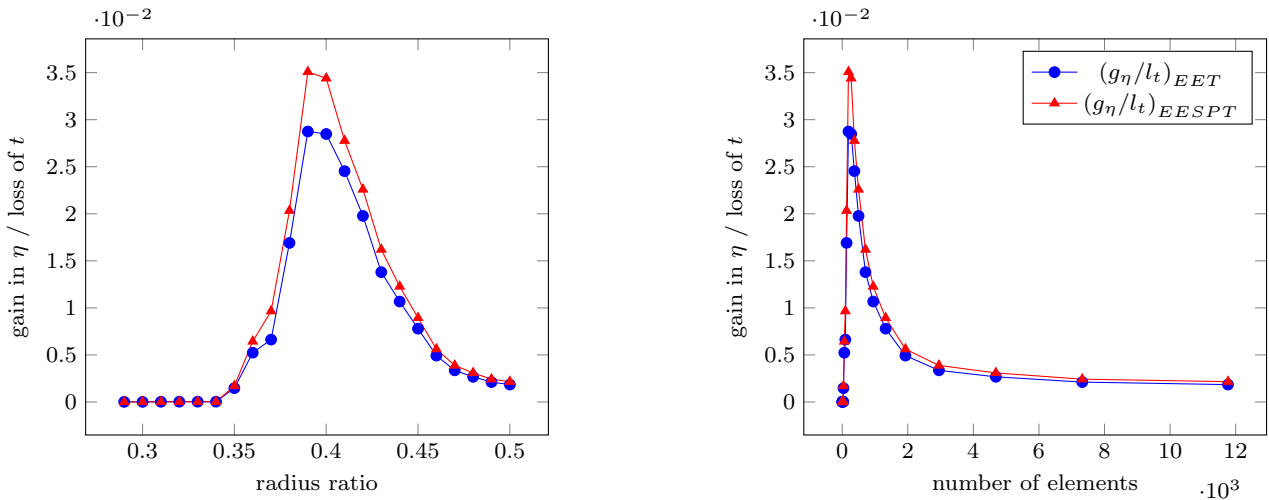
\begin{figure*}
\centering
\begin{tikzpicture}[baseline]
\pgfplotsset{
xlabel near ticks,
ylabel near ticks,
tick label style={font=\footnotesize},
label style={font=\small},
legend style={font=\small},
try min ticks=7
}
\begin{axis}[
	width=0.35\textwidth,
	scaled y ticks=base 10:2,
	scale only axis,
	xlabel=radius ratio,
	ylabel=gain in $\eta$ / loss of $t$,
	x filter/.code={
		\ifnum\coordindex<1
			\def\pgfmathresult{}
		\fi
	}
]
\addplot+[sharp plot,blue,solid,mark=*,mark options={blue}] table[x=rad_r,y expr=((2.387400-\thisrow{eta_EET})/2.387400)/((\thisrow{CPU_EET}-12)/12)] {capteur_effort_2D_radius_ratio.txt};
\addplot+[sharp plot,red,solid,mark=triangle*,mark options={red}] table[x=rad_r,y expr=((2.394960-\thisrow{eta_EESPT})/2.394960)/((\thisrow{CPU_EESPT}-14)/14)] {capteur_effort_2D_radius_ratio.txt};
\end{axis}
\end{tikzpicture}
\hspace{2cm}
\begin{tikzpicture}[baseline]
\pgfplotsset{
xlabel near ticks,
ylabel near ticks,
tick label style={font=\footnotesize},
label style={font=\small},
legend style={font=\small},
try min ticks=7
}
\begin{axis}[
	width=0.35\textwidth,
	scaled x ticks=base 10:-3,
	scaled y ticks=base 10:2,
	scale only axis,
	xlabel=number of elements,
	ylabel=gain in $\eta$ / loss of $t$,
	legend pos=north east,
	legend entries={$\displaystyle{(g_{\eta} / l_t )}_{EET}$, $\displaystyle{(g_{\eta} / l_t )}_{EESPT}$},
	x filter/.code={
		\ifnum\coordindex<1
			\def\pgfmathresult{}
		\fi
	}
]
\addplot+[sharp plot,blue,solid,mark=*,mark options={blue}] table[x=nb_elem,y expr=((2.387400-\thisrow{eta_EET})/2.387400)/((\thisrow{CPU_EET}-12)/12)] {capteur_effort_2D_radius_ratio.txt};
\addplot+[sharp plot,red,solid,mark=triangle*,mark options={red}] table[x=nb_elem,y expr=((2.394960-\thisrow{eta_EESPT})/2.394960)/((\thisrow{CPU_EESPT}-14)/14)] {capteur_effort_2D_radius_ratio.txt};
\end{axis}
\end{tikzpicture}
\caption{Ratio between gain in effectivity index and loss of normalized CPU time for the error estimates EET and EESPT with respect to the radius ratio (left) and to the number of elements involved in the optimized procedure (right).}\label{fig1:capteur_effort_2D_radius_ratio_bis}
\end{figure*}

The graphs represented in \Figs{fig1:capteur_effort_2D_radius_ratio_bis}, \ref{fig1:capteur_effort_2D_edge_ratio_bis} and \ref{fig1:capteur_effort_2D_estimate_ratio_bis} show the evolution of the efficiency factors introduced in \Sect{6.1} with respect to each criterion (l.h.s. graphs) and the number of elements involved in the optimization (r.h.s. graphs) for both EET and EESPT methods. Results are similar to the previous two-dimensional case and proves that a local optimization involving only a small part of the elements of the initial FE mesh suffices to achieve sharper upper bounds.

\begin{figure*}
\centering
\begin{tikzpicture}[baseline]
\pgfplotsset{
xlabel near ticks,
ylabel near ticks,
tick label style={font=\footnotesize},
label style={font=\small},
legend style={font=\small},
try min ticks=7
}
\begin{axis}[
	width=0.35\textwidth,
	scaled y ticks=base 10:2,
	scale only axis,
	xlabel=edge ratio,
	ylabel=gain in $\eta$ / loss of $t$,
	x filter/.code={
		\ifnum\coordindex<1
			\def\pgfmathresult{}
		\fi
	}
]
\addplot+[sharp plot,blue,solid,mark=*,mark options={blue}] table[x=edg_r,y expr=((2.387400-\thisrow{eta_EET})/2.387400)/((\thisrow{CPU_EET}-12)/12)] {capteur_effort_2D_edge_ratio.txt};
\addplot+[sharp plot,red,solid,mark=triangle*,mark options={red}] table[x=edg_r,y expr=((2.394960-\thisrow{eta_EESPT})/2.394960)/((\thisrow{CPU_EESPT}-14)/14)] {capteur_effort_2D_edge_ratio.txt};
\end{axis}
\end{tikzpicture}
\hspace{2cm}
\begin{tikzpicture}[baseline]
\pgfplotsset{
xlabel near ticks,
ylabel near ticks,
tick label style={font=\footnotesize},
label style={font=\small},
legend style={font=\small},
try min ticks=7
}
\begin{axis}[
	width=0.35\textwidth,
	scaled x ticks=base 10:-3,
	scaled y ticks=base 10:2,
	scale only axis,
	xlabel=number of elements,
	ylabel=gain in $\eta$ / loss of $t$,
	legend pos=north east,
	legend entries={$\displaystyle{(g_{\eta} / l_t )}_{EET}$, $\displaystyle{(g_{\eta} / l_t )}_{EESPT}$},
	x filter/.code={
		\ifnum\coordindex<1
			\def\pgfmathresult{}
		\fi
	}
]
\addplot+[sharp plot,blue,solid,mark=*,mark options={blue}] table[x=nb_elem,y expr=((2.387400-\thisrow{eta_EET})/2.387400)/((\thisrow{CPU_EET}-12)/12)] {capteur_effort_2D_edge_ratio.txt};
\addplot+[sharp plot,red,solid,mark=triangle*,mark options={red}] table[x=nb_elem,y expr=((2.394960-\thisrow{eta_EESPT})/2.394960)/((\thisrow{CPU_EESPT}-14)/14)] {capteur_effort_2D_edge_ratio.txt};
\end{axis}
\end{tikzpicture}
\caption{Ratio between gain in effectivity index and loss of normalized CPU time for the error estimates EET and EESPT with respect to the edge ratio (left) and to the number of elements involved in the optimized procedure (right).}\label{fig1:capteur_effort_2D_edge_ratio_bis}
\end{figure*}
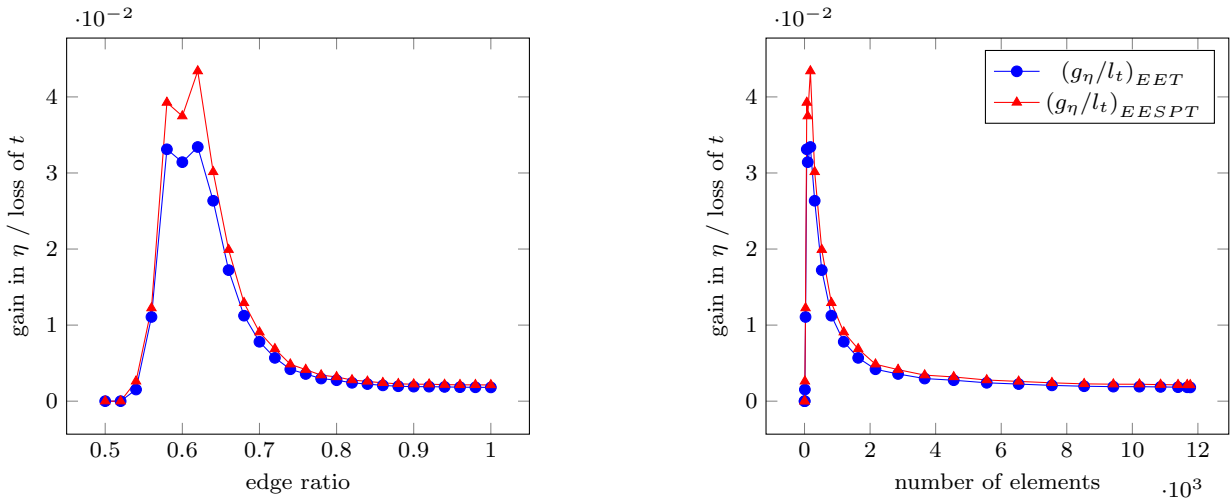

\begin{figure*}
\centering
\begin{tikzpicture}[baseline]
\pgfplotsset{
xlabel near ticks,
ylabel near ticks,
tick label style={font=\footnotesize},
label style={font=\small},
legend style={font=\small},
try min ticks=7
}
\begin{axis}[
	width=0.35\textwidth,
	scaled y ticks=base 10:2,
	scale only axis,
	ymax=0.16,
	xlabel=estimate ratio,
	x dir=reverse,
	ylabel=gain in $\eta$ / loss of $t$,
	x filter/.code={
		\ifnum\coordindex<1
			\def\pgfmathresult{}
		\fi
	}
]
\addplot+[sharp plot,blue,solid,mark=*,mark options={blue}] table[x=est_r,y expr=((2.387400-\thisrow{eta_EET})/2.387400)/((\thisrow{CPU_EET}-12)/12)] {capteur_effort_2D_estimate_ratio.txt};
\addplot+[sharp plot,red,solid,mark=triangle*,mark options={red}] table[x=est_r,y expr=((2.394960-\thisrow{eta_EESPT})/2.394960)/((\thisrow{CPU_EESPT}-14)/14)] {capteur_effort_2D_estimate_ratio.txt};
\end{axis}
\end{tikzpicture}
\hspace{2cm}
\begin{tikzpicture}[baseline]
\pgfplotsset{
xlabel near ticks,
ylabel near ticks,
tick label style={font=\footnotesize},
label style={font=\small},
legend style={font=\small},
try min ticks=7
}
\begin{axis}[
	width=0.35\textwidth,
	scaled x ticks=base 10:-3,
	scaled y ticks=base 10:2,
	scale only axis,
	ymax=0.16,
	xlabel=number of elements,
	ylabel=gain in $\eta$ / loss of $t$,
	legend pos=north east,
	legend entries={$\displaystyle{(g_{\eta} / l_t )}_{EET}$, $\displaystyle{(g_{\eta} / l_t )}_{EESPT}$},
	x filter/.code={
		\ifnum\coordindex<1
			\def\pgfmathresult{}
		\fi
	}
]
\addplot+[sharp plot,blue,solid,mark=*,mark options={blue}] table[x=nb_elem_EET,y expr=((2.387400-\thisrow{eta_EET})/2.387400)/((\thisrow{CPU_EET}-12)/12)] {capteur_effort_2D_estimate_ratio.txt};
\addplot+[sharp plot,red,solid,mark=triangle*,mark options={red}] table[x=nb_elem_EESPT,y expr=((2.394960-\thisrow{eta_EESPT})/2.394960)/((\thisrow{CPU_EESPT}-14)/14)] {capteur_effort_2D_estimate_ratio.txt};
\end{axis}
\end{tikzpicture}
\caption{Ratio between gain in effectivity index and loss of normalized CPU time for the error estimates EET and EESPT with respect to the estimate ratio (left) and to the number of elements involved in the optimized procedure (right).}\label{fig1:capteur_effort_2D_estimate_ratio_bis}
\end{figure*}
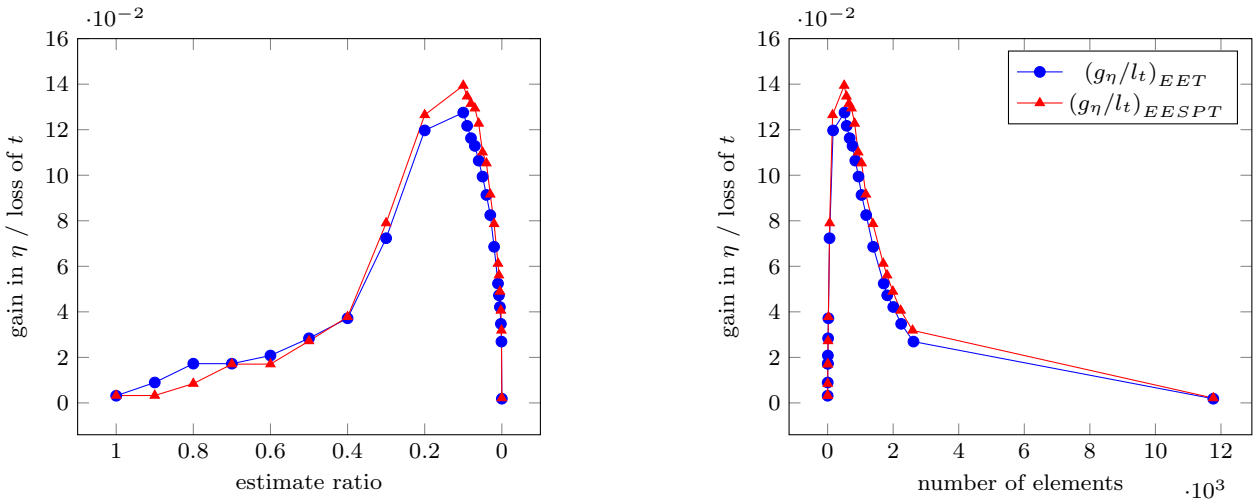

\subsection{Plate with a hole}\label{6.3}

\begin{figure*}
\centering\includegraphics[width=0.9\textwidth]{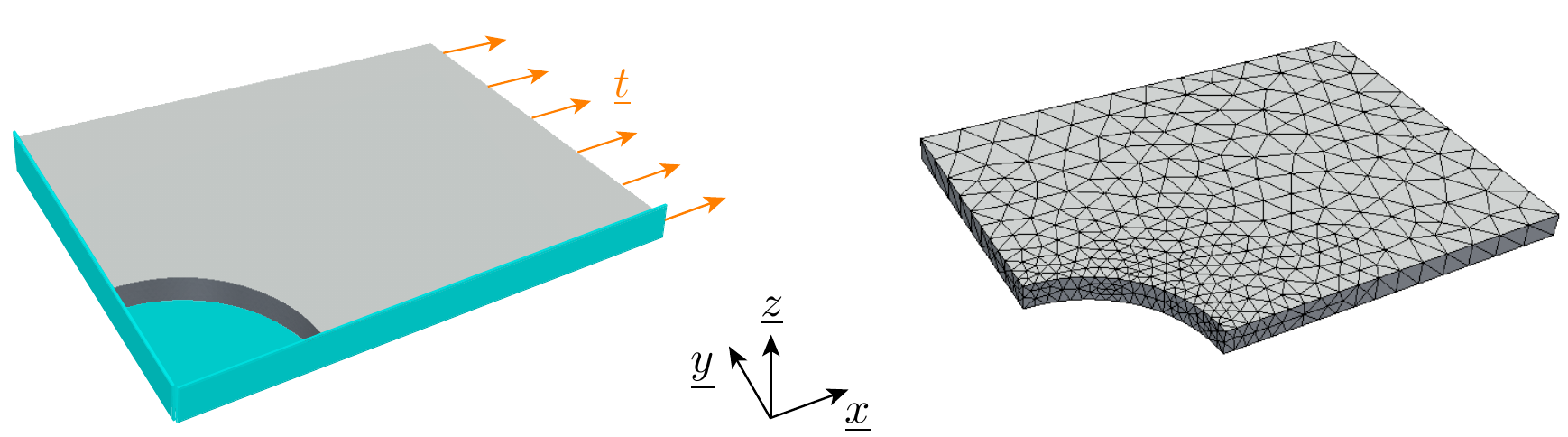}
\caption{Plate with a hole model problem (left) and associated finite element mesh (right). Light blue plans represent symmetry boundary conditions.}
\label{fig1:plaque_trouee_3D_geometry_coarse_mesh}
\end{figure*}

Let us now consider a three-dimensional open-hole specimen modeled by a perforated plate, represented in \Fig{fig1:plaque_trouee_3D_geometry_coarse_mesh}. The plate is $20$ mm long, $15$ mm large, $1$ mm high and presents a hole of radius $2.5$ mm. Due to symmetry, only one eighth of the structure is modeled. Symmetry boundary conditions are applied on the light blue surfaces represented in \Fig{fig1:plaque_trouee_3D_geometry_coarse_mesh}. The structure is loaded along the right side with a unit traction force density $\und{t} = +\und{x}$. The hole and the top side are traction-free boundaries. The FE mesh contains $2 \, 075$ linear tetrahedral elements and $724$ nodes (\ie $2 \, 172$ d.o.f.), is given in \Fig{fig1:plaque_trouee_3D_geometry_coarse_mesh}.
Once again, the reference solution is not the exact one, but the one given considering an approximation space of very large dimension. Indeed, the refined mesh obtained by subsplitting each tetrahedron into $512$ tetrahedra, thus consisting of $1 \, 062 \, 400$ linear tetrahedral elements and $199 \, 293$ nodes (\ie $597 \, 879$ d.o.f.).
\begin{figure}
\centering\includegraphics[scale = 0.32]{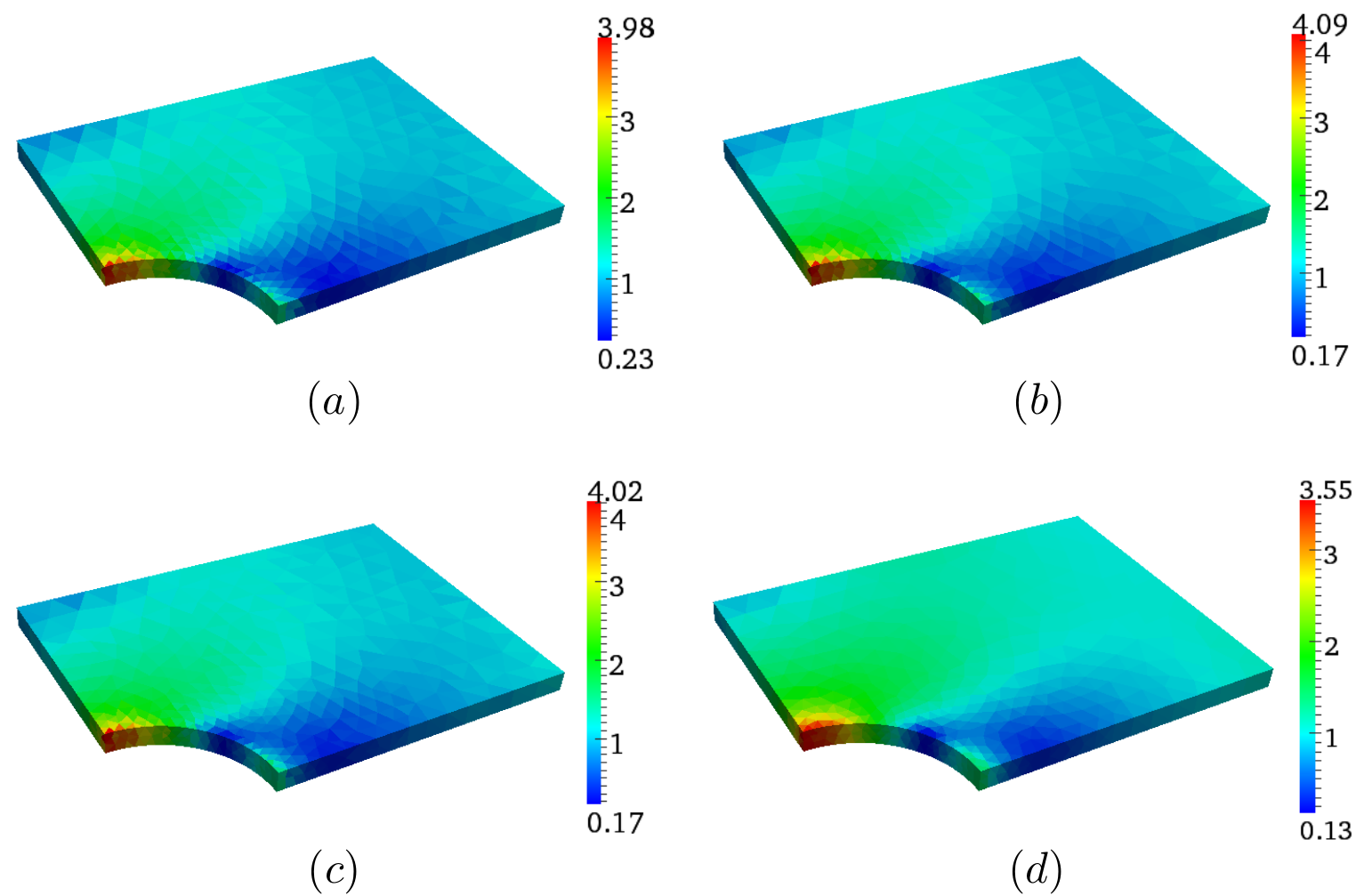}
\caption{Magnitude of the FE stress field (a) and the admissible stress field calculated using the standard versions of EET (b) and EESPT (c), the full enhanced version (d).}
\label{fig1:sigma_EF_hat_plaque_trouee_3D}
\end{figure}
\begin{figure}
\centering\includegraphics[scale = 0.32]{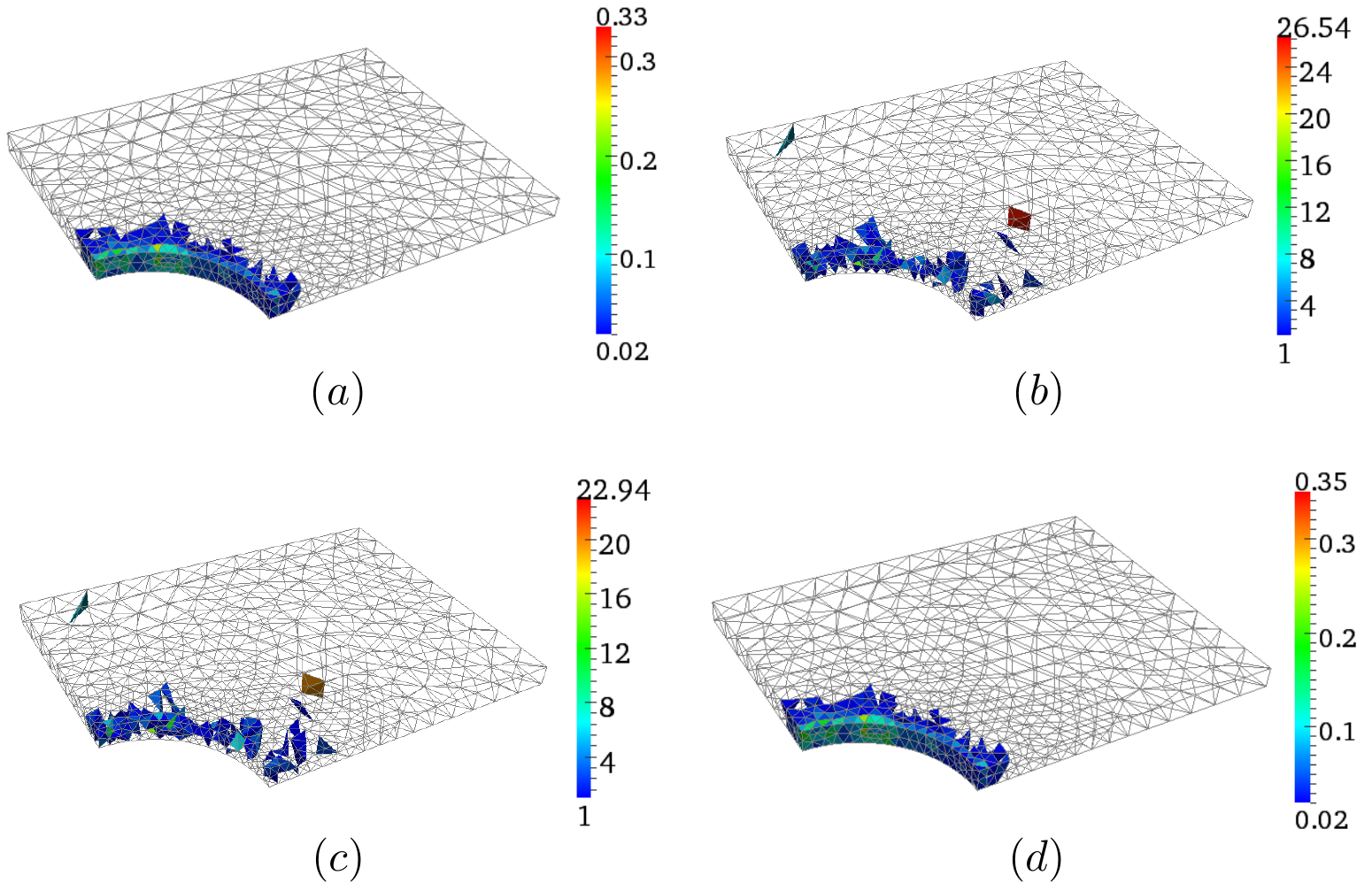}
\caption{Spatial distribution of relevant local contributions to the density of the energy norm of the reference error (a) and that of the contributions to the density of the error estimates calculated using the standard versions of the EET (b) and the EESPT (c), and the full enhanced version (d).}
\label{fig1:distribution_discretization_error_estimators_plaque_trouee_3D}
\end{figure}

\Fig{fig1:sigma_EF_hat_plaque_trouee_3D} represents the FE stress field and the admissible stress fields computed using the standard versions of EET and EESPT methods and the full enhanced version. \Fig{fig1:distribution_discretization_error_estimators_plaque_trouee_3D} shows only the elements which contribute the most to the density of the energy norm of the reference error and to the density of the error estimates. As expected, the region located in the neighborhood of the hole contains the major part of the reference and estimated error. Furthermore, the full enhanced version seems to be less affected by the distortion of the FE mesh than the standard versions of EET and EESPT methods. Indeed, one can see that some ill-shaped elements take a significant part in the error estimates calculated using the standard versions of EET and EESPT methods, whereas these are not involved in the main contributions to the energy norm of the reference error and in that of the error estimates obtained from the full enhanced version.

The calculation of the value of the energy norm of the reference error leads to:
\begin{equation}
\begin{aligned}
\lnorm{\und{e}_h}_{u, \Om} & = \sqrt{\lnorm{\uu}^2_{u, \Om} - \lnorm{\uu_h}^2_{u, \Om}} \\
& \simeq \sqrt{\lnorm{\uu_{ref}}^2_{u, \Om} - \lnorm{\uu_h}^2_{u, \Om}} \simeq 0.368999,
\end{aligned}
\end{equation}
and requires a CPU time of about $1$ s, whereas the calculation of the local contributions to $\lnorm{\und{e}_h}_{u, \Om}$ takes about $1$ hour and $25$ minutes.

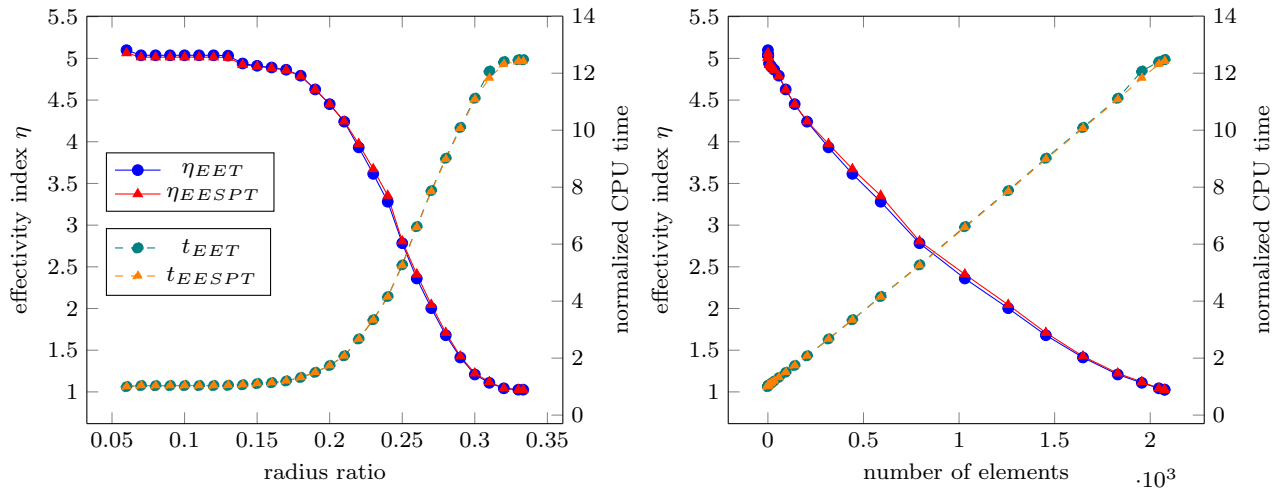
\begin{figure*}
\centering
\begin{tikzpicture}[baseline]
\pgfplotsset{
xlabel near ticks,
ylabel near ticks,
tick label style={font=\footnotesize},
label style={font=\small},
legend style={font=\small},
try min ticks=8
}
\begin{axis}[
	width=0.36\textwidth,
	scale only axis,
	xticklabel style={/pgf/number format/fixed},
	axis y line*=left,
	axis on top,
	xlabel=radius ratio,
	ylabel=effectivity index $\eta$,
	legend style={at={(0.04,0.52)},anchor=south west,legend columns=1},
	legend entries={$\eta_{EET}$,$\eta_{EESPT}$}
]
\addplot+[sharp plot,blue,solid,mark=*,mark options={blue}] table[x=rad_r,y=eta_EET] {plaque_trouee_3D_radius_ratio.txt};
\addplot+[sharp plot,red,solid,mark=triangle*,mark options={red}] table[x=rad_r,y=eta_EESPT] {plaque_trouee_3D_radius_ratio.txt};
\end{axis}
\begin{axis}[
	width=0.36\textwidth,
	scale only axis,
	ymax=14,
	axis x line=none,
	axis y line*=right,
	ylabel=normalized CPU time,
	legend style={at={(0.04,0.48)},anchor=north west,legend columns=1},
	legend entries={$t_{EET}$,$t_{EESPT}$}
]
\addplot+[sharp plot,teal,dashed,mark=*,mark options={teal}] table[x=rad_r,y expr=\thisrow{CPU_EET}/1199] {plaque_trouee_3D_radius_ratio.txt};
\addplot+[sharp plot,orange,dashed,mark=triangle*,mark options={orange}] table[x=rad_r,y expr=\thisrow{CPU_EESPT}/1205] {plaque_trouee_3D_radius_ratio.txt};
\end{axis}
\end{tikzpicture}
\begin{tikzpicture}[baseline]
\pgfplotsset{
xlabel near ticks,
ylabel near ticks,
tick label style={font=\footnotesize},
label style={font=\small},
legend style={font=\small},
try min ticks=8,
}
\begin{axis}[
	width=0.36\textwidth,
	scaled x ticks=base 10:-3,
	scale only axis,
	axis y line*=left,
	axis on top,
	xlabel=number of elements,
	ylabel=effectivity index $\eta$,
	xtick={0,500,1000,1500,2000},
	legend style={at={(0.96,0.96)},anchor=north east,legend columns=1}
]
\addplot+[sharp plot,blue,solid,mark=*,mark options={blue}] table[x=nb_elem,y=eta_EET] {plaque_trouee_3D_radius_ratio.txt};
\addplot+[sharp plot,red,solid,mark=triangle*,mark options={red}] table[x=nb_elem,y=eta_EESPT] {plaque_trouee_3D_radius_ratio.txt};
\end{axis}
\begin{axis}[
	width=0.36\textwidth,
	scale only axis,
	ymax=14,
	axis x line=none,
	axis y line*=right,
	ylabel=normalized CPU time,
	legend style={at={(0.04,0.48)},anchor=north west,legend columns=1}
]
\addplot+[sharp plot,teal,dashed,mark=*,mark options={teal}] table[x=nb_elem,y expr=\thisrow{CPU_EET}/1199] {plaque_trouee_3D_radius_ratio.txt};
\addplot+[sharp plot,orange,dashed,mark=triangle*,mark options={orange}] table[x=nb_elem,y expr=\thisrow{CPU_EESPT}/1205] {plaque_trouee_3D_radius_ratio.txt};
\end{axis}
\end{tikzpicture}
\caption{Effectivity indices and normalized CPU time for the error estimates EET and EESPT with respect to the radius ratio (left) and to the number of elements involved in the optimized procedure (right).}\label{fig1:plaque_trouee_3D_radius_ratio}
\end{figure*}

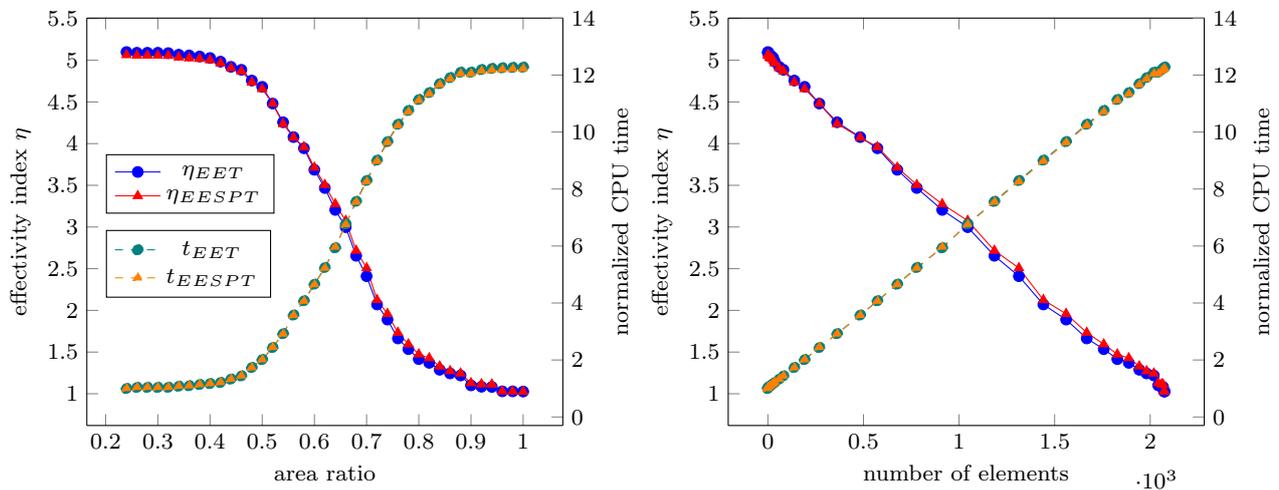
\begin{figure*}
\centering
\begin{tikzpicture}[baseline]
\pgfplotsset{
xlabel near ticks,
ylabel near ticks,
tick label style={font=\footnotesize},
label style={font=\small},
legend style={font=\small},
try min ticks=8
}
\begin{axis}[
	width=0.36\textwidth,
	scale only axis,
	axis y line*=left,
	axis on top,
	xlabel=area ratio,
	ylabel=effectivity index $\eta$,
	legend style={at={(0.04,0.52)},anchor=south west,legend columns=1},
	legend entries={$\eta_{EET}$,$\eta_{EESPT}$}
]
\addplot+[sharp plot,blue,solid,mark=*,mark options={blue}] table[x=are_r,y=eta_EET] {plaque_trouee_3D_area_ratio.txt};
\addplot+[sharp plot,red,solid,mark=triangle*,mark options={red}] table[x=are_r,y=eta_EESPT] {plaque_trouee_3D_area_ratio.txt};
\end{axis}
\begin{axis}[
	width=0.36\textwidth,
	scale only axis,
	ymax=14,
	axis x line=none,
	axis y line*=right,
	ylabel=normalized CPU time,
	legend style={at={(0.04,0.48)},anchor=north west,legend columns=1},
	legend entries={$t_{EET}$,$t_{EESPT}$}
]
\addplot+[sharp plot,teal,dashed,mark=*,mark options={teal}] table[x=are_r,y expr=\thisrow{CPU_EET}/1199] {plaque_trouee_3D_area_ratio.txt};
\addplot+[sharp plot,orange,dashed,mark=triangle*,mark options={orange}] table[x=are_r,y expr=\thisrow{CPU_EESPT}/1205] {plaque_trouee_3D_area_ratio.txt};
\end{axis}
\end{tikzpicture}
\begin{tikzpicture}[baseline]
\pgfplotsset{
xlabel near ticks,
ylabel near ticks,
tick label style={font=\footnotesize},
label style={font=\small},
legend style={font=\small},
try min ticks=8
}
\begin{axis}[
	width=0.36\textwidth,
	scaled x ticks=base 10:-3,
	scale only axis,
	axis y line*=left,
	axis on top,
	xlabel=number of elements,
	ylabel=effectivity index $\eta$,
	legend style={at={(0.96,0.96)},anchor=north east,legend columns=1}
]
\addplot+[sharp plot,blue,solid,mark=*,mark options={blue}] table[x=nb_elem,y=eta_EET] {plaque_trouee_3D_area_ratio.txt};
\addplot+[sharp plot,red,solid,mark=triangle*,mark options={red}] table[x=nb_elem,y=eta_EESPT] {plaque_trouee_3D_area_ratio.txt};
\end{axis}
\begin{axis}[
	width=0.36\textwidth,
	scale only axis,
	ymax=14,
	axis x line=none,
	axis y line*=right,
	ylabel=normalized CPU time,
	legend style={at={(0.04,0.48)},anchor=north west,legend columns=1}
]
\addplot+[sharp plot,teal,dashed,mark=*,mark options={teal}] table[x=nb_elem,y expr=\thisrow{CPU_EET}/1199] {plaque_trouee_3D_area_ratio.txt};
\addplot+[sharp plot,orange,dashed,mark=triangle*,mark options={orange}] table[x=nb_elem,y expr=\thisrow{CPU_EESPT}/1205] {plaque_trouee_3D_area_ratio.txt};
\end{axis}
\end{tikzpicture}
\caption{Effectivity indices and normalized CPU time for the error estimates EET and EESPT with respect to the area ratio (left) and to the number of elements involved in the optimized procedure (right).}\label{fig1:plaque_trouee_3D_area_ratio}
\end{figure*}

\begin{figure*}
\centering
\begin{tikzpicture}[baseline]
\pgfplotsset{
xlabel near ticks,
ylabel near ticks,
tick label style={font=\footnotesize},
label style={font=\small},
legend style={font=\small},
try min ticks=8
}
\begin{axis}[
	width=0.36\textwidth,
	scale only axis,
	axis y line*=left,
	axis on top,
	xlabel=estimate ratio,
	x dir=reverse,
	ylabel=effectivity index $\eta$,
	legend style={at={(0.04,0.52)},anchor=south west,legend columns=1},
	legend entries={$\eta_{EET}$,$\eta_{EESPT}$}
]
\addplot+[sharp plot,blue,solid,mark=*,mark options={blue}] table[x=est_r,y=eta_EET] {plaque_trouee_3D_estimate_ratio.txt};
\addplot+[sharp plot,red,solid,mark=triangle*,mark options={red}] table[x=est_r,y=eta_EESPT] {plaque_trouee_3D_estimate_ratio.txt};
\end{axis}
\begin{axis}[
	width=0.36\textwidth,
	scale only axis,
	ymax=14,
	axis x line=none,
	axis y line*=right,
	ylabel=normalized CPU time,
	x dir=reverse,
	legend style={at={(0.04,0.48)},anchor=north west,legend columns=1},
	legend entries={$t_{EET}$,$t_{EESPT}$}
]
\addplot+[sharp plot,teal,dashed,mark=*,mark options={teal}] table[x=est_r,y expr=\thisrow{CPU_EET}/1199] {plaque_trouee_3D_estimate_ratio.txt};
\addplot+[sharp plot,orange,dashed,mark=triangle*,mark options={orange}] table[x=est_r,y expr=\thisrow{CPU_EESPT}/1205] {plaque_trouee_3D_estimate_ratio.txt};
\end{axis}
\end{tikzpicture}
\begin{tikzpicture}[baseline]
\pgfplotsset{
xlabel near ticks,
ylabel near ticks,
tick label style={font=\footnotesize},
label style={font=\small},
legend style={font=\small},
try min ticks=8
}
\begin{axis}[
	width=0.36\textwidth,
	scaled x ticks=base 10:-3,
	scale only axis,
	axis y line*=left,
	axis on top,
	xlabel=number of elements,
	ylabel=effectivity index $\eta$,
	legend style={at={(0.96,0.96)},anchor=north east,legend columns=1}
]
\addplot+[sharp plot,blue,solid,mark=*,mark options={blue}] table[x=nb_elem_EET,y=eta_EET] {plaque_trouee_3D_estimate_ratio.txt};
\addplot+[sharp plot,red,solid,mark=triangle*,mark options={red}] table[x=nb_elem_EESPT,y=eta_EESPT] {plaque_trouee_3D_estimate_ratio.txt};
\end{axis}
\begin{axis}[
	width=0.36\textwidth,
	scale only axis,
	ymax=14,
	axis x line=none,
	axis y line*=right,
	ylabel=normalized CPU time,
	legend style={at={(0.04,0.48)},anchor=north west,legend columns=1}
]
\addplot+[sharp plot,teal,dashed,mark=*,mark options={teal}] table[x=nb_elem_EET,y expr=\thisrow{CPU_EET}/1199] {plaque_trouee_3D_estimate_ratio.txt};
\addplot+[sharp plot,orange,dashed,mark=triangle*,mark options={orange}] table[x=nb_elem_EESPT,y expr=\thisrow{CPU_EESPT}/1205] {plaque_trouee_3D_estimate_ratio.txt};
\end{axis}
\end{tikzpicture}
\caption{Effectivity indices and normalized CPU time for the error estimates EET and EESPT with respect to the estimate ratio (left) and to the number of elements involved in the optimized procedure (right).}\label{fig1:plaque_trouee_3D_estimate_ratio}
\end{figure*}
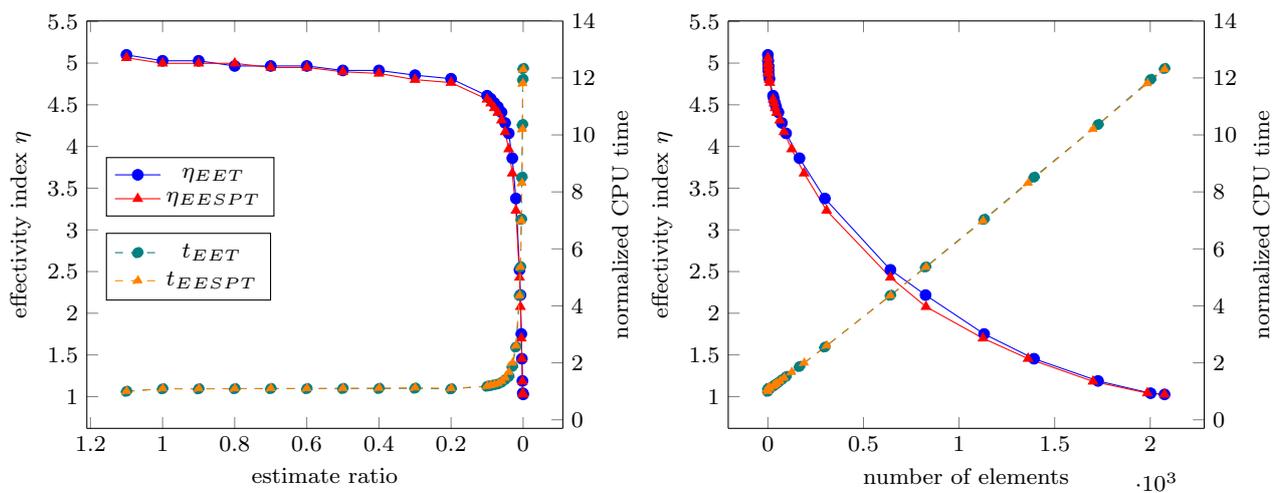

\Figs{fig1:plaque_trouee_3D_radius_ratio}, \ref{fig1:plaque_trouee_3D_area_ratio} and \ref{fig1:plaque_trouee_3D_estimate_ratio} show the effectivity indices and corresponding normalized computational cost as functions of radius ratio, area ratio and estimate ratio, respectively. The radius ratio, area ratio and estimate ratio range between $0.06689$ to $0.3298$, $0.2564$ to $0.9900$ and $1.77 \, 10^{-4}$ to $1.0$, respectively. Values of the effectivity indices and normalized CPU times corresponding to values $0.06$, $0.24$ and $1.1$ for the radius ratio, area ratio and estimate ratio, respectively, are that calculated using the standard versions of EET and EESPT methods. As regards both geometric criteria, results are similar and show a quasi-linear evolution of the effectivity indices with respect to the number of elements implicated in the optimized procedure. Concerning the error estimate criterion, values of the effectivity indices do not decline as strongly as for the two-dimensional cases, but however they experience a substantial decrease compared to the ones obtained using both geometric criteria.

The efficiency factors $g_{\eta} / l_t$, defined by relations (\ref{eq1:efficiencyfactor}) in \Sect{6.1}, computed with respect to the EET and EESPT methods, are plotted versus the radius ratio, the area ratio and the estimate ratio (l.h.s. graphs), and versus the number of elements involved in the enhanced procedure (r.h.s. graphs) respectively, in \Figs{fig1:plaque_trouee_3D_radius_ratio_bis}, \ref{fig1:plaque_trouee_3D_area_ratio_bis} and \ref{fig1:plaque_trouee_3D_estimate_ratio_bis}. Once again, upper bounds of better quality can be obtained without impairing too much the computational cost and with no need to perform any mesh refinement by using an error estimate criterion in order to optimize the tractions only in local zones in which the contributions to the estimate are greatest.

\begin{figure*}
\centering
\begin{tikzpicture}[baseline]
\pgfplotsset{
xlabel near ticks,
ylabel near ticks,
tick label style={font=\footnotesize},
label style={font=\small},
legend style={font=\small},
try min ticks=7
}
\begin{axis}[
	width=0.35\textwidth,
	scaled y ticks=base 10:2,
	xticklabel style={/pgf/number format/fixed},
	scale only axis,
	xlabel=radius ratio,
	ylabel=gain in $\eta$ / loss of $t$,
	x filter/.code={
		\ifnum\coordindex<1
			\def\pgfmathresult{}
		\fi
	}
]
\addplot+[sharp plot,blue,solid,mark=*,mark options={blue}] table[x=rad_r,y expr=((5.098370-\thisrow{eta_EET})/5.098370)/((\thisrow{CPU_EET}-1199)/1199)] {plaque_trouee_3D_radius_ratio.txt};
\addplot+[sharp plot,red,solid,mark=triangle*,mark options={red}] table[x=rad_r,y expr=((5.063770-\thisrow{eta_EESPT})/5.063770)/((\thisrow{CPU_EESPT}-1205)/1205)] {plaque_trouee_3D_radius_ratio.txt};
\end{axis}
\end{tikzpicture}
\hspace{2cm}
\begin{tikzpicture}[baseline]
\pgfplotsset{
xlabel near ticks,
ylabel near ticks,
tick label style={font=\footnotesize},
label style={font=\small},
legend style={font=\small},
try min ticks=7
}
\begin{axis}[
	width=0.35\textwidth,
	scaled x ticks=base 10:-3,
	scaled y ticks=base 10:2,
	scale only axis,
	xlabel=number of elements,
	ylabel=gain in $\eta$ / loss of $t$,
	legend pos=north east,
	legend entries={$\displaystyle{(g_{\eta} / l_t )}_{EET}$, $\displaystyle{(g_{\eta} / l_t )}_{EESPT}$},
	x filter/.code={
		\ifnum\coordindex<1
			\def\pgfmathresult{}
		\fi
	}
]
\addplot+[sharp plot,blue,solid,mark=*,mark options={blue}] table[x=nb_elem,y expr=((5.098370-\thisrow{eta_EET})/5.098370)/((\thisrow{CPU_EET}-1199)/1199)] {plaque_trouee_3D_radius_ratio.txt};
\addplot+[sharp plot,red,solid,mark=triangle*,mark options={red}] table[x=nb_elem,y expr=((5.063770-\thisrow{eta_EESPT})/5.063770)/((\thisrow{CPU_EESPT}-1205)/1205)] {plaque_trouee_3D_radius_ratio.txt};
\end{axis}
\end{tikzpicture}
\caption{Ratio between gain in effectivity index and loss of normalized CPU time for the error estimates EET and EESPT with respect to the radius ratio (left) and to the number of elements involved in the optimized procedure (right).}\label{fig1:plaque_trouee_3D_radius_ratio_bis}
\end{figure*}
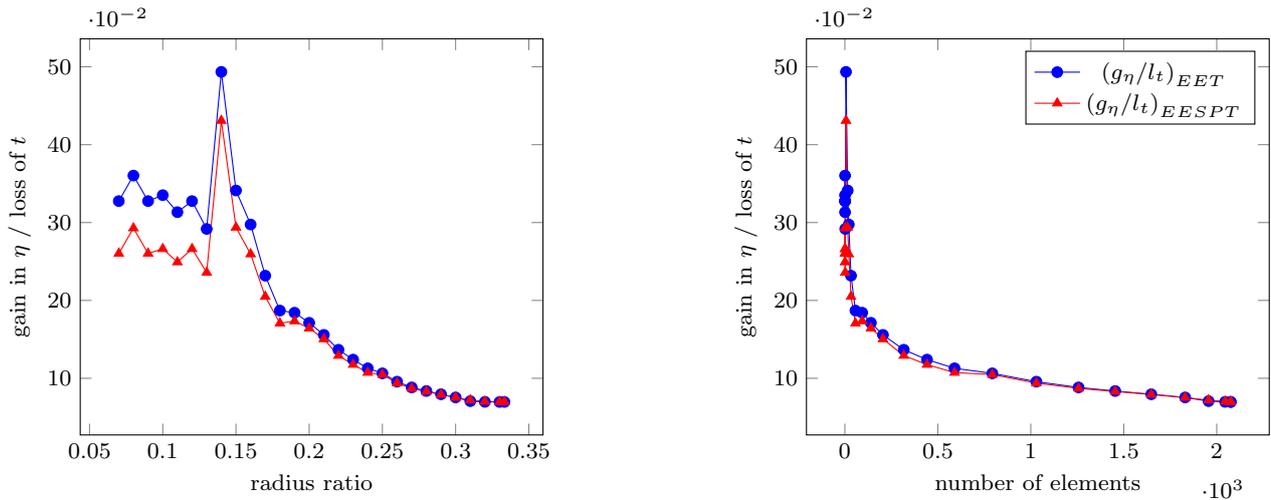

\begin{figure*}
\centering
\begin{tikzpicture}[baseline]
\pgfplotsset{
xlabel near ticks,
ylabel near ticks,
tick label style={font=\footnotesize},
label style={font=\small},
legend style={font=\small},
try min ticks=7
}
\begin{axis}[
	width=0.35\textwidth,
	scaled y ticks=base 10:2,
	scale only axis,
	xlabel=area ratio,
	ylabel=gain in $\eta$ / loss of $t$,
	x filter/.code={
		\ifnum\coordindex<1
			\def\pgfmathresult{}
		\fi
	}
]
\addplot+[sharp plot,blue,solid,mark=*,mark options={blue}] table[x=are_r,y expr=((5.098370-\thisrow{eta_EET})/5.098370)/((\thisrow{CPU_EET}-1199)/1199)] {plaque_trouee_3D_area_ratio.txt};
\addplot+[sharp plot,red,solid,mark=triangle*,mark options={red}] table[x=are_r,y expr=((5.063770-\thisrow{eta_EESPT})/5.063770)/((\thisrow{CPU_EESPT}-1205)/1205)] {plaque_trouee_3D_area_ratio.txt};
\end{axis}
\end{tikzpicture}
\hspace{2cm}
\begin{tikzpicture}[baseline]
\pgfplotsset{
xlabel near ticks,
ylabel near ticks,
tick label style={font=\footnotesize},
label style={font=\small},
legend style={font=\small},
try min ticks=7
}
\begin{axis}[
	width=0.35\textwidth,
	scaled x ticks=base 10:-3,
	scaled y ticks=base 10:2,
	scale only axis,
	xlabel=number of elements,
	ylabel=gain in $\eta$ / loss of $t$,
	legend pos=north east,
	legend entries={$\displaystyle{(g_{\eta} / l_t )}_{EET}$, $\displaystyle{(g_{\eta} / l_t )}_{EESPT}$},
	x filter/.code={
		\ifnum\coordindex<1
			\def\pgfmathresult{}
		\fi
	}
]
\addplot+[sharp plot,blue,solid,mark=*,mark options={blue}] table[x=nb_elem,y expr=((5.098370-\thisrow{eta_EET})/5.098370)/((\thisrow{CPU_EET}-1199)/1199)] {plaque_trouee_3D_area_ratio.txt};
\addplot+[sharp plot,red,solid,mark=triangle*,mark options={red}] table[x=nb_elem,y expr=((5.063770-\thisrow{eta_EESPT})/5.063770)/((\thisrow{CPU_EESPT}-1205)/1205)] {plaque_trouee_3D_area_ratio.txt};
\end{axis}
\end{tikzpicture}
\caption{Ratio between gain in effectivity index and loss of normalized CPU time for the error estimates EET and EESPT with respect to the area ratio (left) and to the number of elements involved in the optimized procedure (right).}\label{fig1:plaque_trouee_3D_area_ratio_bis}
\end{figure*}
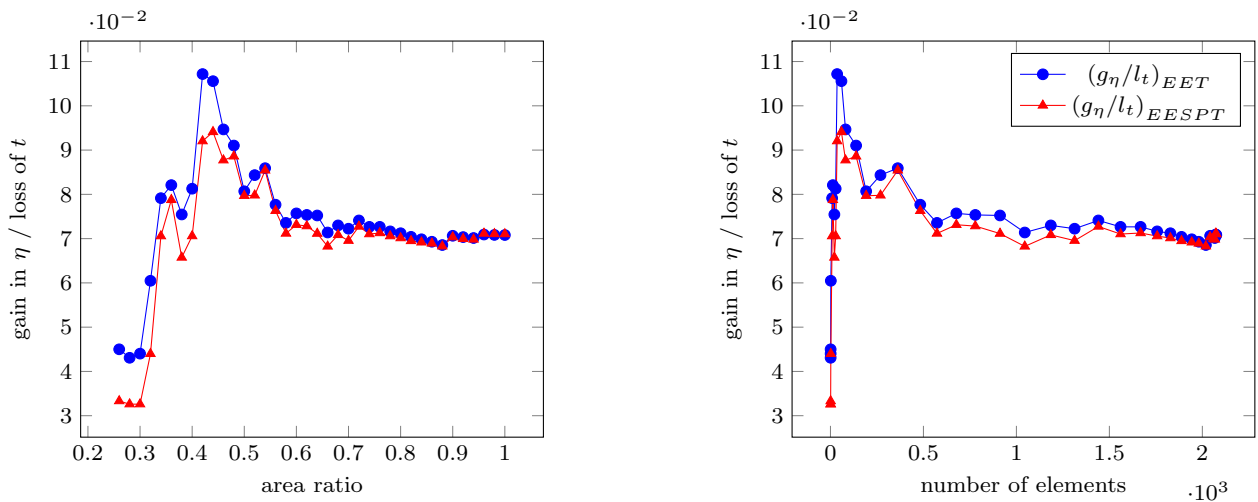

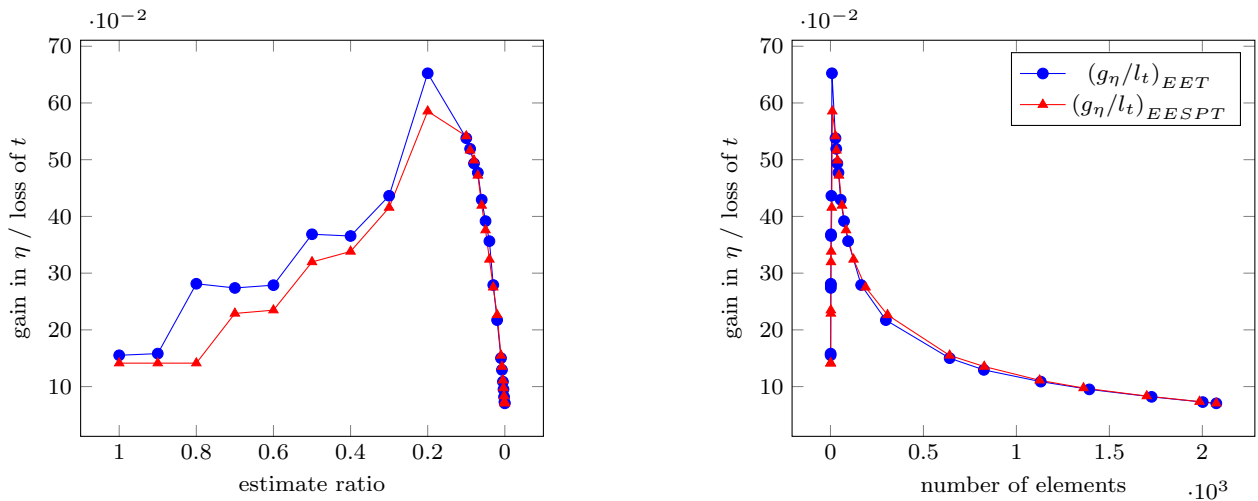
\begin{figure*}
\centering
\begin{tikzpicture}[baseline]
\pgfplotsset{
xlabel near ticks,
ylabel near ticks,
tick label style={font=\footnotesize},
label style={font=\small},
legend style={font=\small},
try min ticks=7
}
\begin{axis}[
	width=0.35\textwidth,
	scaled y ticks=base 10:2,
	scale only axis,
	xlabel=estimate ratio,
	x dir=reverse,
	ylabel=gain in $\eta$ / loss of $t$,
	x filter/.code={
		\ifnum\coordindex<1
			\def\pgfmathresult{}
		\fi
	}
]
\addplot+[sharp plot,blue,solid,mark=*,mark options={blue}] table[x=est_r,y expr=((5.098370-\thisrow{eta_EET})/5.098370)/((\thisrow{CPU_EET}-1199)/1199)] {plaque_trouee_3D_estimate_ratio.txt};
\addplot+[sharp plot,red,solid,mark=triangle*,mark options={red}] table[x=est_r,y expr=((5.063770-\thisrow{eta_EESPT})/5.063770)/((\thisrow{CPU_EESPT}-1205)/1205)] {plaque_trouee_3D_estimate_ratio.txt};
\end{axis}
\end{tikzpicture}
\hspace{2cm}
\begin{tikzpicture}[baseline]
\pgfplotsset{
xlabel near ticks,
ylabel near ticks,
tick label style={font=\footnotesize},
label style={font=\small},
legend style={font=\small},
try min ticks=7
}
\begin{axis}[
	width=0.35\textwidth,
	scaled x ticks=base 10:-3,
	scaled y ticks=base 10:2,
	scale only axis,
	xlabel=number of elements,
	ylabel=gain in $\eta$ / loss of $t$,
	legend pos=north east,
	legend entries={$\displaystyle{(g_{\eta} / l_t )}_{EET}$, $\displaystyle{(g_{\eta} / l_t )}_{EESPT}$},
	x filter/.code={
		\ifnum\coordindex<1
			\def\pgfmathresult{}
		\fi
	}
]
\addplot+[sharp plot,blue,solid,mark=*,mark options={blue}] table[x=nb_elem_EET,y expr=((5.098370-\thisrow{eta_EET})/5.098370)/((\thisrow{CPU_EET}-1199)/1199)] {plaque_trouee_3D_estimate_ratio.txt};
\addplot+[sharp plot,red,solid,mark=triangle*,mark options={red}] table[x=nb_elem_EESPT,y expr=((5.063770-\thisrow{eta_EESPT})/5.063770)/((\thisrow{CPU_EESPT}-1205)/1205)] {plaque_trouee_3D_estimate_ratio.txt};
\end{axis}
\end{tikzpicture}
\caption{Ratio between gain in effectivity index and loss of normalized CPU time for the error estimates EET and EESPT with respect to the estimate ratio (left) and to the number of elements involved in the optimized procedure (right).}\label{fig1:plaque_trouee_3D_estimate_ratio_bis}
\end{figure*}

\section{Conclusion and Prospects}\label{7}

We studied an enhanced version of the procedure used for the calculation of balanced tractions over element edges involved in both EET and EESPT methods. We pointed out that the optimized procedure offers greater flexibility. Therefore, the use of such a procedure in order to improve the quality of the tractions constructed achieves better efficiency than the standard version, although it requires higher computational cost. The novelty of this work is embodied in the different criteria used in the enhanced procedure. The analysis of the results reflects a downward trend of the global effectivity indices of each estimator. The optimized procedure concerning the sensitivity to error contribution leads to a better computational efficiency than the ones regarding the sensitivity to geometric parameters. Indeed the global effectivity indices for both estimators experience a sharp decrease by using error estimate criterion, while the use of geometric criteria yields a slight downturn in the effectivity index. Thus, this enhanced procedure is very effective and the improvement is particularly significant at affordable cost when it brings only the mostly concentrated error elements of the FE mesh into play. As regards robust error estimation, a way to achieve a sharper estimate while keeping a low-cost error estimation procedure and a given FE mesh is the use of the enhanced procedure with error estimate criterion. Therefore, it will be used in a forthcoming work dealing with global/goal-oriented error estimation.


\bibliographystyle{spmpsci}      
\bibliography{Biblio}   


\end{document}